\documentclass[leqno]{cornell}
\title{COHOMOLOGY OF $AUT(F_n)$}
\author{Craig A. Jensen}
\date{August 1998}

\usepackage{amsfonts}
\newcommand{\Z}{\mathbb{Z}}
\newcommand{\Q}{\mathbb{Q}}
\newcommand{\F}{\mathbb{F}}

\newcommand{\A}{\mathbb{A}}

\input amssym.def
\newsymbol\square 1003
\newsymbol\lozenge 1006
\newsymbol\surjarrow 1310
\newsymbol\injarrow 131A
\newsymbol\ltimes 226E
\newsymbol\rtimes 226F
\newsymbol\void 203F
\newsymbol\leftrightsquigarrow 1321
\def\END {$\square$}
\def\PF {\noindent{\bf Proof. }}
\def\proclaim #1. #2\par #3\par {\medbreak
\noindent{\bf#1.\enspace}{\sl#2}\par\medbreak
\noindent{\bf Proof.} #3 \par}
\def\proclaimb #1. #2\par {\medbreak
\noindent{\bf#1.\enspace}{\sl#2}\par\medbreak}
\font\esr = eusm10 scaled \magstep 0
\font\esrc = eusm10 scaled \magstep 2
\font\ecr = eurm10 scaled \magstep 0
\newtheorem{thm}[equation]{Theorem}
\newtheorem{lemma}[equation]{Lemma}
\newtheorem{prop}[equation]{Proposition}
\newtheorem{cor}[equation]{Corollary}
\newtheorem{defn}[equation]{Definition}
\newtheorem{claim}[equation]{Claim}
\newtheorem{remark}[equation]{Remark}
\newtheorem{fact}[equation]{Fact}
\begin{document}

\maketitle

\newpage

\makecopyright

\newpage

\begin{abstract}

For odd primes $p$, we examine the Farrell cohomology
$\hat H^*(Aut(F_{2(p-1)}); \Z_{(p)})$ of the group of automorphisms of
a free group $F_{2(p-1)}$ on $2(p-1)$ generators, with coefficients
in the integers localized at the prime $(p) \subset \Z$.  This extends
results in \cite{[G-M]} by Glover and Mislin, whose calculations
yield
$\hat H^*(Aut(F_n); \Z_{(p)})$ for $n \in \{ p-1,p\}$
and is concurrent with work by Chen in \cite{[C]}
where he calculates
$\hat H^*(Aut(F_n); \Z_{(p)})$ for $n \in \{ p+1,p+2\}.$
The main tools used are
Ken Brown's ``Normalizer spectral sequence'' from \cite{[B]},
a modification of Krstic and Vogtmann's proof of the contractibility
of fixed point sets for outer space in \cite{[K-V]}, and
a modification of the Degree Theorem of Hatcher and
Vogtmann in \cite{[H-V]}.

Other cohomological calculations in the paper yield that $H^5(Q_m; \Z)$
never stabilizes as $m \to \infty$,
where $Q_m$ is the quotient of the spine $X_m$ of ``auter space''
introduced in \cite{[H-V]} by Hatcher and Vogtmann.  This contrasts
with the theorems in
\cite{[H-V]} where various stability results are shown for
$H^n(Aut(F_m); \Z)$, $H^n(Aut(F_m); \Q)$, and $H^n(Q_m; \Q)$.

\end{abstract}

\begin{biosketch}

Born in Davis, California, in 1969, Craig has lived in
Davis; Livermore, California; Malvern, Pennsylvania;
Los Alamos, New Mexico; Logan, Utah;
Madison, Wisconsin; and Ithaca, New York.  He graduated Summa
Cum Laude from Utah State University in 1992 with a
B.S. in both Mathematics and Philosophy.  Following this,
he attended the University of Wisconsin--Madison and received an
M.A. in Mathematics in 1994.  His graduate education
was completed by earning a Ph.D. in Mathematics from
Cornell University in 1998.

Craig is an avid reader of books, and has read (perhaps
to the detriment of his mathematical education!) at
least three fiction books every week since he was
11 years old.  His other hobbies include
philosophy (in which he was sufficiently interested
in as a hobby to obtain a degree),
poetry (alas, he is better at reading it than writing it),
games (all types), and
computers (a reasonable hobby for someone in
a technical discipline.)  Overall, he has enjoyed
his time spent obtaining a Ph.D. in Mathematics,
but is glad finally to have completed this
facet of his education.  He looks forward to
his future position as a Zassenhaus
Assistant Professor at Ohio State University.

\end{biosketch}

\begin{dedication}

{\ecr Dedicated with love to my mother}

{\esrc J}{\esr OANNE }{\esrc R}{\esr EID }{\esrc J}{\esr ENSEN} 

{\ecr for helping me and being a great friend}

\bigskip

-- {\ecr and} --

\bigskip

{\ecr In loving memory of my father}

{\esrc C}{\esr ARL }{\esrc A J}{\esr ENSEN} 

{\ecr wishing that we had had more time}

\end{dedication}

\begin{acknowledgements}

First, I would like to thank Karen Vogtmann,
my thesis advisor, for all that she has done.  This
dissertation would not have been possible without her
help and advice.  She has done all of the things
that a good thesis advisor could do, including
offering advice about possible avenues
of research, suggesting references where I
could find information, giving encouragement
where necessary, commenting on the things I
have written, and sitting through 
hours of what must have been fairly tiresome
presentations of mine.  In short,
Prof. Vogtmann has helped me by doing various
things which are difficult, if not impossible,
for a graduate student to do on his own.
It has been a pleasure to know her, and I
am grateful for her aid and assistance.

Second, I would like to thank the other members of
my special committee, Allen Hatcher
and Marshall Cohen.  Professor Hatcher,
who at one time was chair of my committee,
was of immense help to me throughout my
time at Cornell.  He has
frequently offered advice,
suggested refinements to proofs,
helped me locate articles, and also listened
to a good number of presentations that
I have done.  Professor Cohen has been
very kind and encouraging to me while I
have attended Cornell.
I am glad to have known him.

Next I would like to thank two other people
who have helped me complete the process of
obtaining a dissertation and
completing graduate school,
Ken Brown and Alejandro Adem.  Both of
them have kindly answered 
questions (which I seem, like clockwork,
to have posed to each of
them once every two months)
about topics that have come up in the
cohomology of groups.  Their answers have
contained invaluable information,
such as references, summaries of topics, and guesses
about how difficult or complex they felt the
subjects were.  Ken Brown, in particular, deserves the
award for ``most helpful person at
Cornell who is not officially on my
committee'', as he has both written letters of
recommendation for me and served as a
proxy member of my committee.

While attending graduate school, I have
had the honor of being a Hertz
Fellow and have been supported by
a Hertz Fellowship.
The fellowship has been
of constant help for me, by allowing me
to devote more time and effort
to studying and researching
mathematics.
I would like to
thank the Fannie and John Hertz
Foundation for awarding me the
fellowship and funding me.
In particular, I would like to thank
Wilson Talley and Lowell Wood
for interviewing me and
recommending me for the fellowship.

Although I always wanted to be some
sort of scientist, it is
because of my high school math
teacher, Carletta Elich, that
I specifically became a mathematician.
In addition to posing
interesting problems and ideas
to me, she was genuinely kind
to me as a person and
concerned about my education.  I am
happy to call her a friend.

Finally, I would like to thank
my family.  My mother, Joanne,
has always encouraged me
to pursue the goal of
earning a Ph.D. in the sciences.
Throughout the past ten years
of university education, she
has bestowed priceless
love, support, and counsel.
Although it goes
without saying that I would not
have finished this process
without her assistance, 
it still needs to be said that
she has consistently
been there for me when I
needed love and friendship the most.
I would also like to thank my sister,
JanLee, for being such a good friend
and someone whose abilities
I have always admired.
Last of all, I would like to thank
my brother, Doug, for all the fun times
we have had and for being my
best friend ever.

\end{acknowledgements}

\contentspage

\figurelistpage

\preface
\doublespacing

Let $F_n$ denote the free group on $n$ letters and let
$Aut(F_n)$ and $Out(F_n)$ denote the automorphism group
and outer automorphism group, respectively, of $F_n$.
In \cite{[C-V]} Culler and Vogtmann defined a space on which
$Out(F_n)$ acts nicely called ``outer space''.  By
studying the action of $Out(F_n)$ on this space, various
people have been able to calculate the cohomology
of $Out(F_n)$ in specific cases.  More recently, Hatcher in
\cite{[H]} and Hatcher and Vogtmann in \cite{[H-V]} have defined
a space on which $Aut(F_n)$ acts nicely called ``auter space''
and have used this to calculate the cohomology of
$Aut(F_n)$ in specific cases.

Recall (see \cite{[L-S]}) that $Out(F_1)=\Z/2$
and $Out(F_2)=GL_2(\Z)$; in general, for all $n \geq 1$ the map
$Out(F_n) \to GL_n(\Z)$ has a torsion free kernel.
In \cite{[Br]}, Tom Brady calculated the integral cohomology of $Out(F_3)$.
This remains even today the only complete (nontrivial) calculation of
the integral cohomology of $Out(F_n)$ or $Aut(F_n)$.

Hatcher shows in \cite{[H]} that the integral cohomology of
the infinite symmetric group $\Sigma_\infty$ is a direct
summand of the integral cohomology of $Aut(F_\infty)$.
It is unknown whether or not the cohomology of $\Sigma_\infty$
equals the cohomology of $Aut(F_\infty)$.
In \cite{[H]} he also derives quadratic stability ranges
for the integral cohomology of $Aut(F_n)$ or $Out(F_n)$.

A ``Degree Theorem'' is introduced by Hatcher and Vogtmann
in \cite{[H-V]} which is a very useful tool for
simplifying cohomological calculations
concerning $Aut(F_n)$.
For example, they are able to derive linear stability
ranges for the integral cohomology of
$Aut(F_n)$.

Using the degree theorem and computer calculations,
progress toward calculating
the rational cohomology of $Aut(F_n)$ is made by
Hatcher and Vogtmann in \cite{[V]} where they show that
for $1 \leq t \leq 6$ and all $n \geq 1$,
$$\matrix{
\hfill H_t(Aut(F_n); \Q) &=& &\left\{\matrix{
\Q \hfill &\hbox{if } t = n = 4 \hfill \cr
0 \hfill &\hbox{otherwise} \hfill \cr} \right. \hfill \cr
}$$

Glover and Mislin \cite{[G-M]} calculated the
cohomology with coefficients in $\Z_{(p)}$ of $Out(F_n)$
for $n=p-1,p,p+1$.  In addition, Chen \cite{[Ch]} calculates
the integral cohomology of $Out(F_n)$ for $n=p+2$ and
$Aut(F_n)$ for $n=p+1,p+2$.  In each of the above cases,
the maximal $p$-subgroups of $Out(F_n)$ or $Aut(F_n)$ had $p$-rank one.
The case of $Aut(F_{2(p-1)})$ is the first one where
the maximal $p$-subgroups can have a higher $p$-rank, and it is
this case which we calculate here.

\noindent In this paper we show that for odd primes $p$, and
$n=2(p-1)$, that 
$\hat H^t(Aut(F_{n}); \Z_{(p)})$

\medskip 

$$\matrix{
\cong \left\{\matrix{
\Z/p^2 \oplus p(\Z/p) \hfill &t = 0 \hfill \cr
(p + [{3k \over 2}] - 1)\Z/p \hfill &|t| = kn \not = 0 \hfill \cr
\Z/p \hfill &t = 1 \hfill \cr
0 \hfill &t = kn + 1 > 1 \hfill \cr
([{3k \over 2}] - 1)\Z/p \hfill &t = -kn+1 < 0 \hfill \cr
H^{n-1}(\tilde Q_{p-1}; \Z/p) \oplus
([{3k \over 2}] - 1)\Z/p \hfill &t = kn-1 > 0 \hfill \cr
H^{n-1}(\tilde Q_{p-1}; \Z/p) \hfill &t = -kn-1 < 0 \hfill \cr
H^r(Q_{p}^\omega;\Z/p) \oplus
\sum_{i=0}^{p-1} H^r(\tilde Q_{i} \times Q_{p-1-i};\Z/p)
\hfill &t = kn + r, 2 \leq r \leq n-2 \hfill \cr
} \right. \hfill \cr}$$

\bigskip

\noindent or alternatively that
$\hat H^t(Aut(F_{n}); \Z_{(p)}) = $ 

\medbreak

$$\matrix{
\cong \left\{\matrix{
\Z/p^2 \oplus p(\Z/p) \hfill &t = 0 \hfill \cr
 \hfill & \hfill \cr
(p + [{3k \over 2}] - 1)\Z/p \hfill &|t| = kn \not = 0 \hfill \cr
 \hfill & \hfill \cr
\Z/p \hfill &t = 1 \hfill \cr
 \hfill & \hfill \cr
0 \hfill &t = kn + 1 > 1 \hfill \cr
 \hfill & \hfill \cr
([{3k \over 2}] - 1)\Z/p \hfill &t = -kn+1 < 0 \hfill \cr
 \hfill & \hfill \cr
([{3k \over 2}] - 1)\Z/p \hfill &t = kn-1 > 0 \hfill \cr
\oplus H^{n-1}(F_{p-1}
\rtimes Aut(F_{p-1}); \Z/p) \hfill & \hfill \cr
 \hfill & \hfill \cr
H^{n-1}(F_{p-1} \rtimes Aut(F_{p-1}); \Z/p)  \hfill &t = -kn-1 < 0 \hfill \cr
 \hfill & \hfill \cr
H^r((F_{p-2} \times F_{p-2})
\rtimes (\Z/2 \times Aut(F_{p-2}));\Z/p)  \oplus
\hfill &t = kn + r, \hfill \cr
\sum_{i=0}^{p-1} H^r((F_i \rtimes Aut(F_i)) \times Aut(F_{p-1-i});\Z/p)
\hfill &2 \leq r \leq n-2 \hfill \cr
} \right. \hfill \cr}$$

\bigskip

\noindent where $Q_{k}$, $\tilde Q_{k}$, and $Q_{p}^\omega$
are cellular complexes of dimensions $2k-2$, $2k-1$, and  $2p-4$,
respectively (provided $k \not = 0$; both $Q_0$ and
$\tilde Q_0$ can be taken to be points);
they are quotient complexes of various contractible
complexes by various groups and 
will be defined later in this paper.
Although the second formulation of 
$$\hat H^*(Aut(F_{2(p-1)}); \Z_{(p)})$$
above might seem more useful as it does not involve the
above quotient complexes (that is,  $Q_{k}$, $\tilde Q_{k}$, and 
$Q_{p}^\omega$), it (arguably) contains less
information than the first formulation.

A quick note about our notation is appropriate here.
In general, groups without any additional
structure will be written using multiplicative notation
(e.g., $\Z/p \times \Z/p \cong (\Z/p)^2$) but modules like cohomology
groups will be written using additive notation
(e.g., $\Z/p \oplus \Z/p \cong 2(\Z/p)$.)  Hence
the case $t=0$ of
our main result above should be read as stating that
$$\hat H^0(Aut(F_{2(p-1)}); \Z_{(p)})
\cong \Z/p^2 \oplus (\Z/p \oplus \cdots \oplus \Z/p)$$
where there are $p$ copies of $\Z/p$ in
$(\Z/p \oplus \cdots \oplus \Z/p)$.

For the case of the prime $p=3$,
the cohomology groups of all of the
relevant quotient spaces are calculated in
Appendix \ref{appen} of this paper.  This allows us
to state concisely what the above result gives us in the case
$p=3$.  This calculation was also (independently)
done by Glover and Henn for both $Out(F_4)$
and $Aut(F_4)$.  In the appendix, we show that
our main result simplifies to the following in
the case $p=3$:

\medskip

$$\hat H^t(Aut(F_{4}); \Z_{(3)}) = 
\left\{\matrix{
\Z/9 \oplus 3(\Z/3) \hfill &t = 0 \hfill \cr
([{3k \over 2}] + 2)\Z/3 \hfill &|t| = 4k \not = 0 \hfill \cr
\Z/3 \hfill &t = 1 \hfill \cr
0 \hfill &|t| \equiv 1,2  \hbox{ } (\hbox{mod } 4)
\hbox{ and } t \not = 1 \hfill \cr
([{3k \over 2}] - 1)\Z/3 \hfill &|t| = 4k-1 \hfill \cr
} \right.$$
\medskip

A couple more remarks about the statement of our main result are in
order.  First, you will note that all calculations
were done with $\Z_{(p)}$, the integers localized at the
prime ideal $(p)$, as the coefficient ring.  This was done
for two reasons.  One is that taking coefficients in
$\Z_{(p)}$ is like taking coefficients in $\Z/p$ in that it 
allows us (most of the time) just to look at simplices with ``p-symmetry''
in the various spectral sequences that we will use to calculate
the cohomology groups.  In this way, calculations are simplified.
Another reason is that by taking coefficients in
$\Z_{(p)}$ -- rather than $\Z/p$ -- we are still able to
distinguish torsion from non-torsion in the cohomology.
Although this is irrelevant when looking at Farrell
cohomology (as all Farrell cohomology groups are torsion),
it will still have some technical advantages (over taking $\Z/p$
coefficients) when we make statements about the usual
cohomology of $Aut(F_n)$ in chapter \ref{c6} and
for the calculations about moduli spaces of graphs
in part \ref{p4}.

Second, Glover and Mislin used the standard spectral sequence for
equivariant cohomology covered by Brown in \cite{[B]} on page 173 or
282.  It seems difficult to use that spectral sequence to calculate
$\hat H^*(Aut(F_{2(p-1)}); \Z_{(p)})$ and so instead we used the
``normalizer spectral sequence'' discussed by Brown in \cite{[B]} on
page 293.  The normalizer spectral sequence is also discussed by Henn
in \cite{[HN]}, where in addition work giving a ``centralizer spectral
sequence'' is done.

There are equivariant cohomology spectral sequences for
both the usual cohomology and the Farrell cohomology
of a group acting on a space.  Unfortunately, the normalizer spectral
sequence only applies to Farrell cohomology.  Hence our results
are stated in terms of Farrell cohomology.  In chapter \ref{c6}
we state what our results imply about the usual cohomology
of $Aut(F_n)$.  In general, it is a more difficult problem to
calculate the usual cohomology of a group than to calculate its
Farrell cohomology.  For example, although the results of
Glover and Mislin about $Out(F_p)$ can be used quickly to derive the
Farrell cohomology of $Aut(F_p)$, the usual cohomology of
$Aut(F_p)$ has yet to be calculated ``precisely''.  This problem
still arises in Chen's calculations with $Aut(F_{p+1})$
and $Aut(F_{p+2})$.

For the convenience of the reader, we now state these
previous (or concurrent, in the case of Chen's work)
calculations concerning the cohomology of $Aut(F_n)$.
From the work of Glover and Mislin in \cite{[G-M]}, it can be seen
that for odd primes $p$,
$$\matrix{
\hfill H^t(Aut(F_{p-1}); \Z_{(p)}) &=& &\left\{\matrix{
H^t(\Sigma_p; \Z_{(p)}) \hfill &t > 2p-4 \hfill \cr
H^t(Q_{p-1}; \Z_{(p)}) \hfill &t \leq 2p-4 \hfill \cr} \right. \hfill \cr
}$$
where $Q_m$ is the quotient of the spine $X_m$ of auter space
and $\Sigma_p$ is the symmetric group on $p$ letters.

Restating this in terms of Farrell cohomology,
we find that for all $t$,
$$\matrix{
\hfill \hat H^t(Aut(F_{p-1}); \Z_{(p)}) &=& \hat H^t(\Sigma_p;
\Z_{(p)}). \hfill \cr}$$
Similarly, it is clear from Glover and Mislin's work that
$$\matrix{
\hfill \hat H^t(Aut(F_{p}); \Z_{(p)}) &=& 3\hat H^t(\Sigma_p;
\Z_{(p)}), \hfill \cr
}$$
and Chen in \cite{[C]} shows that if $p \geq 5$ then
$$\matrix{
\hfill \hat H^t(Aut(F_{p+1}); \Z_{(p)}) &=& 4\hat H^t(\Sigma_p;
\Z_{(p)}) \hfill \cr
}$$
and that
$$\matrix{
\hfill \hat H^t(Aut(F_{p+2}); \Z_{(p)}) &=& 5\hat H^t(\Sigma_p;
\Z_{(p)}) \oplus \hat H^{t-4}(\Sigma_p; \Z_{(p)}). \hfill \cr
}$$

As a final note,
it is easy to check that $\hat H^t(\Sigma_p; \Z_{(p)})$ is periodic
with period $2(p-1)$ and that
$$\hat H^t(\Sigma_p; \Z_{(p)}) = \left\{\matrix{
\Z/p \hfill &t = 2(p-1)k \hfill \cr
0 \hfill &\hbox{otherwise } \hfill \cr} \right.$$

\newpage

\normalspacing
\pagenumbering{arabic}
\pagestyle{cornell}

\part{Cohomological calculations} \label{p1}

\chapter{The equivariant spectral sequence and the spine
of ``auter space''} \label{c1}

Let $G$ be a group acting cellularly on a
finite dimensional CW-complex $X$ such that the stabilizer
$stab_G(\delta)$ of every cell $\delta$ is finite and such
that the quotient of $X$ by $G$ is finite.  Further suppose that
for every cell $\delta$ of $X$, the group $stab_G(\delta)$ fixes
$\delta$ pointwise.  Let $M$ be a $G$-module. 
Recall (see \cite{[B]}) that the equivariant
cohomology groups of the $G$-complex $X$
with coefficients in $M$ are defined by
$$H_G^*(X; M) = H^*(G; C^*(X;M))$$
and that if in addition $X$ is contractible
(which will usually, but not always, be the case in this
paper) then
$$H_G^*(X; M) = H^*(G;M).$$

In \cite{[B]} a spectral sequence
\begin{equation} \label{e1}
\tilde E_1^{r,s} = \prod_{[\delta] \in \Delta_n^r}
H^s(stab(\delta); M) \Rightarrow H_G^{r+s}(X; M)
\end{equation}
is defined, where $[\delta]$ ranges over the set $\Delta_n^r$ of orbits
of $r$-simplices $\delta$ in $X$.

If $M$ is $\Z/p$ or $\Z_{(p)}$ then
a nice property should be noted about the spectral sequence \ref{e1}.
This property will greatly reduce the calculations we need to
go through, and in general will make concrete
computations possible.
Since each
group $stab(\delta)$ is finite, a standard restriction-transfer argument
in group cohomology yields that $|stab(\delta)|$ annihilates
$H^s(stab(\delta); M)$ for all $s > 0$.
(For examples of these sorts of
arguments see \cite{[A-M]} or \cite{[B]}.)
Since all primes not equal to
$p$ are divisible in $\Z/p$ or $\Z_{(p)}$,
this in turn shows that the $p$-part of
$|stab(\delta)|$ annihilates $H^s(stab(\delta); M)$ for
$s>0$.
In particular, if $p$ does not divide some $|stab(\delta)|$, then this
$[\delta]$ does not contribute anything to the spectral sequence
\ref{e1} except in the horizontal row $s=0$.  It follows that if our
coefficients are $\Z/p$ or $\Z_{(p)}$ then we are mainly just concerned
with the simplices $\delta$ which have ``$p$-symmetry''.

If $G$ is a group with finite virtual cohomological dimension (vcd)
and $M$ is a $G$-module, then Farrell cohomology groups
$$\hat H^*(G; M)$$
are defined in \cite{[F]} and \cite{[B]}.
For the basics about
Farrell cohomology, along with several useful properties,
see \cite{[B]} or \cite{[F]}.  Briefly, if $M$ is
one of the standard coefficient rings like $\Z$, $\Q$,
$\Z/p$, or $\Z_{(p)}$ then
these cohomology groups have
the property that:
$$\hat H^t(G; M) = H^*(G; M) \hbox{ if } t > \hbox{ vcd}(G)$$
and
$$\hat H^*(G; M) = 0 \hbox{ if } G \hbox{ is torsion free.}$$
Moreover, unlike the usual cohomology groups, they tend to be
nonzero for negative indices $t < 0$ as well as
for positive indices.

As before,
suppose $G$ acts cellularly on a
finite dimensional CW-complex $X$ such that the stabilizer
$stab_G(\delta)$ of every cell $\delta$ is finite and such
that the quotient of $X$ by $G$ is finite.  Further suppose that
for every cell $\delta$ of $X$, the group $stab_G(\delta)$ fixes
$\delta$ pointwise.
Just as in the usual case, the equivariant Farrell
cohomology groups of the $G$-complex $X$
with coefficients in $M$ are defined (see \cite{[B]}
for details) by
$$\hat H_G^*(X; M) = \hat H^*(G;C^*(X;M))$$
and if in addition $X$ is contractible
then
$$\hat H_G^*(X; M) = \hat H^*(G;M).$$

For Farrell cohomology, Brown defines in \cite{[B]}
an equivariant cohomology spectral sequence which is
analogous to \ref{e1}.
It is given by
\begin{equation} \label{e2}
E_1^{r,s} = \prod_{[\delta] \in \Delta_n^r}
\hat H^s(stab(\delta); M) \Rightarrow \hat H_G^{r+s}(X; M).
\end{equation}

If $M$ is $\Z/p$ or $\Z_{(p)}$ then as before only the
simplices with ``$p$-symmetry'' are relevant.
That is, we know that  $|stab(\delta)|$ annihilates
$\hat H^s(stab(\delta); \Z_{(p)})$ 
for all $s$ (not just for $s>0$ as before.)  
Since all primes not equal to
$p$ are divisible in M, this shows that the $p$-part of
$|stab(\delta)|$ annihilates $\hat H^s(stab(\delta); M)$ for
all $s$.
Accordingly, if $p$ does not divide some $|stab(\delta)|$, then this
$[\delta]$ does not contribute anything to the spectral sequence
\ref{e2}.

We now specialize to the case where
$G=Aut(F_n)$ and $X$ is the spine $X_n$ of ``auter space.''
First, we review some basic properties and definitions of $Aut(F_n)$ and
auter space.  Most of these can be found in \cite{[C-V]},
\cite{[H-V]}, \cite{[S-V1]}, and \cite{[S-V2]}.
Consider the automorphism group $Aut(F_n)$ of
a free group $F_n$ of rank $n$ (where $n=2(p-1)$ for our work.)
Let $(R_n,v_0)$ be the
$n$-leafed rose,
a wedge of $n$ circles. We say a basepointed graph $(G,x_0)$
is \emph{admissible} if it has no free edges, all vertices except the
basepoint
have valence at least three, and there is a basepoint-preserving continuous
map $\phi\colon R_n \to G$ which induces an isomorphism on $\pi_1$. The
triple $(\phi,G,x_0)$ is called a \emph{marked graph}.  Two marked graphs
$(\phi_i,G_i,x_i) \hbox{ for } i=0,1$  are \emph{equivalent}  if there is
a homeomorphism $\alpha \colon (G_0,x_0) \to (G_1,x_1)$ such that
$ (\alpha\circ\phi_0)_\# = (\phi_1)_\# : \pi_1(R_n,v_0) \to \pi_1(G_1,x_1)$.
Define a partial order on the set of all equivalence classes of
marked graphs by setting $(\phi_0,G_0,x_0) \leq
(\phi_1,G_1,x_1)$ if $G_1$ contains a \emph{forest}
(a disjoint union of trees in $G_1$ which contains all of the vertices
of $G_1$) such that collapsing each tree in
the forest to a point yields $G_0$, where the collapse is compatible with
the maps $\phi_0$ and $\phi_1$.

From \cite{[H]} and \cite{[H-V]} we have that $Aut(F_n)$ acts
with finite stabilizers on
a contractible space $X_n$.
The space $X_n$ is the geometric realization of the poset of
marked graphs that we defined above.
Let $Q_n$ be the quotient of $X_n$ by
$Aut(F_n)$.
Note that the CW-complex $Q_n$ is not necessarily a simplicial
complex.
Since $Aut(F_n)$ has a torsion free subgroup of finite index
\cite{[H]} and it acts on the contractible,
finite dimensional space
$X_n$ with finite stabilizers and finite quotient,
$Aut(F_n)$ has finite vcd.  Thus it makes sense to
talk about its Farrell cohomology, and to apply spectral
sequences like \ref{e2} to calculate that cohomology.
For $n=0$, we set
$X_0 = Q_0 = \{\hbox{a point}\}$
as a notational convenience.

Let $p$ be an odd prime number, and let $\Z_{(p)}$ be the localization of
$\Z$ at the prime ideal $(p)$.  Then we
can apply the spectral sequences \ref{e1} and \ref{e2} to get
\begin{equation} \label{er1}
\tilde E_1^{r,s} = \prod_{[\delta] \in \Delta_n^r}
H^s(stab(\delta); \Z_{(p)}) \Rightarrow H^{r+s}(Aut(F_n); \Z_{(p)})
\end{equation}
and
\begin{equation} \label{er2}
E_1^{r,s} = \prod_{[\delta] \in \Delta_n^r}
\hat H^s(stab(\delta); \Z_{(p)}) \Rightarrow \hat H^{r+s}(Aut(F_n); \Z_{(p)})
\end{equation}
where $[\delta]$ ranges over the set $\Delta_n^r$ of orbits
of $r$-simplices $\delta$ in $X_n$.

The spectral sequences \ref{er1} and \ref{er2}
requires as input the stabilizers $stab_{Aut(F_n)}(\delta)$ of
simplices $\delta$ in $X_n$.  Smillie and Vogtmann \cite{[S-V1]}
examined the structure of these stabilizers in
detail, and we list their results here.
Consider a given $r$-simplex
$$(\phi_r,G_r,x_r) > \cdots  > (\phi_1,G_1,x_1)
> (\phi_0,G_0,x_0)$$
with corresponding forest collapses
$$(H_r \subseteq G_r), \ldots, (H_2 \subseteq G_2), (H_1 \subseteq G_1).$$
For each $i \in {0,1, \ldots, r}$, let $F_i$ be the inverse image
under the map
$$G_r \to \cdots \to G_{i+1} \to G_i$$ of forest collapses,
of the
forest $H_i$.  That is, we
have
$$F_r \subseteq \cdots \subseteq F_2 \subseteq F_1 \subseteq G_r.$$
It is shown in \cite{[S-V1]}
that the stabilizer of the simplex under
consideration is isomorphic to the group
$Aut(G_r,F_1,\ldots,F_r,x_r)$ of basepointed automorphisms of the graph
$G_r$ that respect each of the forests $F_i$.  For example,
the stabilizer of a point $(\phi,G,x_0)$ in $X_n$ is isomorphic to
$Aut(G,x_0)$.

To sum things up, suppose we wanted to use
the spectral sequences \ref{er1} or \ref{er2} to
compute the cohomology of $Aut(F_n)$.
Since our coefficient ring is $\Z_{(p)}$,
we have already remarked that
for the terms in the spectral sequence above the
horizontal axis, we are concerned only with simplices
whose stabilizers are divisible by $p$.
In addition, the stabilizer of a simplex consists of graph automorphisms
that respect the forest collapses in the simplex.  Let us investigate
what simplices arise in the case $n=2(p-1)$.
In other words, we want to calculate which
graphs $G$ with a $\Z/p$
action on them have $\pi_1(G) \cong F_n$.
To simplify our calculations, we will consider only 
reduced $\Z/p$-graphs, where a
$\Z/p$-graph
$G$ is \emph{reduced} if it contains no $\Z/p$-invariant
subforests.

The spectral sequence \ref{er1} is the one used by Glover and
Mislin in \cite{[G-M]} for their calculations of
the cohomology of $Out(F_n)$ for $n=p-1, p, p+1$.
We are interested, however, in the cohomology of $Aut(F_n)$
for $n=2(p-1)$. 
The step of finding which simplices in $X_{2(p-1)}$ have
$p$-symmetry would be a necessary one if we actually
intended to use the spectral sequence \ref{er1} to
calculate the cohomology of $Aut(F_{2(p-1)})$.
Unfortunately, it seems that the spectral sequence \ref{er1}
becomes too cluttered and difficult to
compute with in our case, and so we will not be
using it to calculate our cohomology groups.
Instead, we will use Ken Brown's ``normalizer spectral sequence''
\cite{[B]} which will be introduced in chapter \ref{c2}.
Nevertheless, it will still prove instructive and useful
to examine which simplices in $X_{2(p-1)}$ have
$p$-symmetry.

\begin{figure}[here]
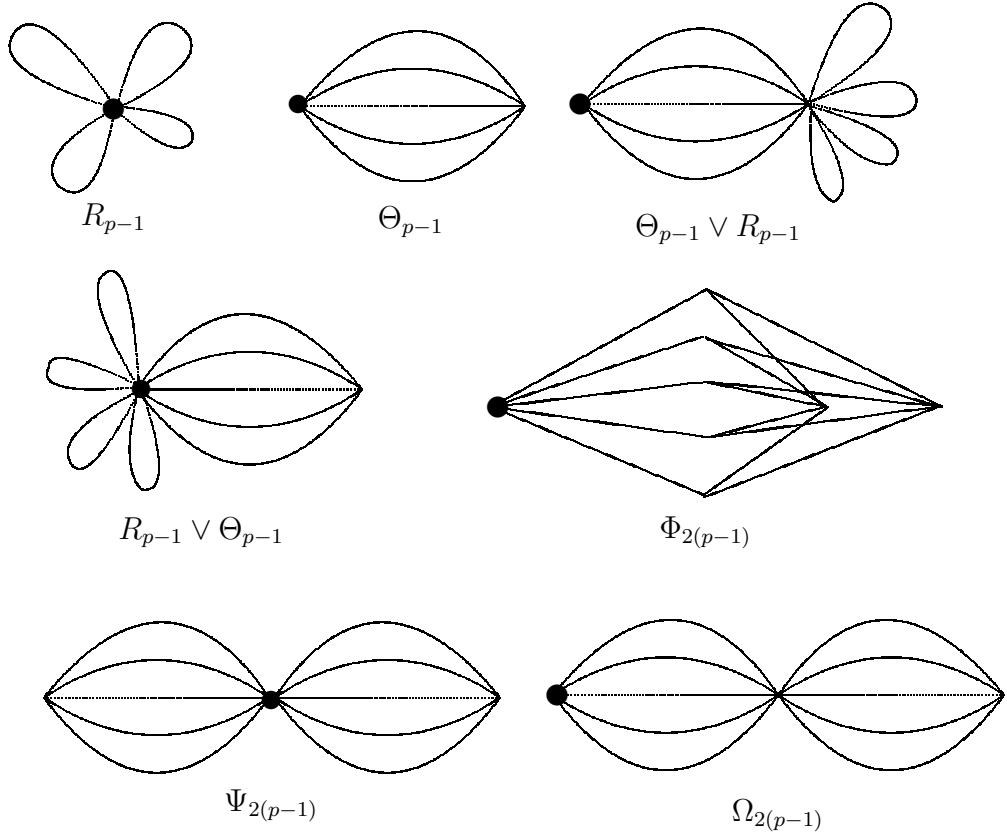

\caption{\label{fig1} Some graphs with $p$-symmetry}
\end{figure}

Some preliminary definitions of a few common graphs are in order.
Let $\Theta_{p-1}$ be the graph with two vertices and $p$ edges, each of
which goes from one vertex to the other (see Figure \ref{fig1}.)
Say the ``leftmost vertex''
of $\Theta_{p-1}$ is the basepoint.  Hence when we write
$\Theta_{p-1} \vee R_{p-1}$ then we are stipulating that the rose
$R_{p-1}$ is attached to the non-basepointed vertex of $\Theta_{p-1}$,
while when we write $R_{p-1} \vee \Theta_{p-1}$ then we are saying
that the rose is attached to the basepoint of $\Theta_{p-1}$.
Let $\Phi_{2(p-1)}$ be a graph with $3p$ edges $a_1, \ldots, a_p,$
$b_1, \ldots, b_p,$ $c_1, \ldots, c_p,$ and $p+3$ vertices
$v_1, \ldots, v_p, x, y, z.$  The basepoint is $x$ and each of the
edges $a_i$ begin at $x$ and end at $v_i$.  The edges $b_i$ and $c_i$
begin at $y$ and $z$, respectively, and end at $v_i$.  Note that there
are obvious actions of $\Z/p$ on $\Theta_{p-1}$ and $\Phi_{2(p-1)}$,
given by rotation, and that these actions are unique up to conjugacy.
Let $\Psi_{2(p-1)}$ be the graph obtained from $\Phi_{2(p-1)}$ by
collapsing all of the edges $a_i$ to a point.  Let $\Omega_{2(p-1)}$
be the graph obtained from $\Phi_{2(p-1)}$ by collapsing either the
edges $b_i$ or the edges $c_i$ (the resulting graphs are isomorphic)
to a point.  Note that the only difference between $\Psi_{2(p-1)}$ and
$\Omega_{2(p-1)}$ is where the basepoint is located.

\input{fig1.pic}

\begin{lemma} \label{tr1}
Let $p$ be an odd prime and $n=2(p-1)$,
and let $\Gamma$
be a reduced (basepointed) graph
with a
nontrivial $\Z/p$-action, where $\pi_1(\Gamma) \cong F_n$.
Let $e$ be an edge of $\Gamma$ which is moved
by the $\Z/p$ action.  Then the orbit of the edge $e$
forms either a rose $R_p$ or a $\Theta$-graph $\Theta_{p-1}$.
\end{lemma}

\PF Say $e_1, e_2, \ldots, e_p$ is the orbit of
the edge $e$ under the $\Z/p$-action.
(By an edge here we mean an unoriented one; that is, a 1-cell
in the 1-dimensional CW-complex $\Gamma$.)

Since $\Gamma$ is reduced,
the edges $e_i$ have some endpoints in common,
else they form an invariant subforest.  Similarly, they cannot have
just one endpoint in common, or we could collapse the star that
they would form.  Finally, they also cannot form a $p$-gon, else
we could take a minimal path from $e_1$ to the basepoint, consider
its orbit under $\Z/p$, and find in those edges a star of $p$ edges
that could be collapsed.  So either all of the endpoints of the
$e_i$ are the same and they form a rose, or there are two distinct
endpoints and they form a $\Theta$-graph. \END

\begin{prop} \label{t1}
Let $p$ be an odd prime and $n=2(p-1)$.
The only reduced (basepointed) graphs $\Gamma$
with a
nontrivial $\Z/p$-action
and $\pi_1(\Gamma) \cong F_n$ are $R_n$,
$R_k \vee \Theta_{p-1} \vee R_{p-1-k}$, $\Psi_n$, and $\Omega_n$.
\end{prop}

\PF First, suppose that $\Gamma$ has only one nontrivial $\Z/p$-orbit of
edges $$e_1, e_2, \ldots, e_p.$$
From the lemma, the $e_i$ either form a rose $R_p$
or a $\Theta$-graph $\Theta_{p-1}$.
In the first case,
$\Gamma$ must be $R_{2(p-1)}$ with $\Z/p$ rotating $p$ of the leaves
and leaving the other $p-2$ fixed.  In the second case, $\Gamma$ must
be of the form $R_k \vee \Theta_{p-1} \vee R_{p-1-k}$
where $\Z/p$ rotates the edges of the $\Theta$-graphs and leaves
the roses at either end of the $\Theta$-graph fixed.

Second, we examine the case where $\Gamma$ has more than one
nontrivial $\Z/p$-orbit of unoriented edges.  Take two
distinct orbits $e_1, e_2, \ldots, e_p$ and
$f_1, f_2, \ldots, f_p$.  As in the previous paragraph,
we can use the fact that $\Gamma$ is reduced
and Lemma \ref{tr1} to get that the $e_i$
and $f_i$ either form roses or $\Theta$-graphs.  If either one
of them is a rose, then the rank of $\pi_1(\Gamma)$ is at
least $2p-1$, which is a contradiction.  So both form
$\Theta$-graphs.  Since $\pi_1(\Gamma) \cong F_{2(p-1)}$, $\Gamma$ doesn't
have any other edges aside from the $e_i$ and $f_i$.
It follows that $\Gamma$ is either $\Psi_n$ or $\Omega_n$,
depending upon where the basepoint is. \END

\chapter{The normalizer spectral sequence} \label{c2}

Ken Brown \cite{[B]} introduces another spectral sequence that can be
used to calculate $\hat H^*(Aut(F_n); \Z_{(p)})$.  It involves normalizers
of elementary abelian $p$-subgroups of the group under consideration
($Aut(F_n)$ in our case) and hence will
be called the ``normalizer spectral sequence.''  This is not the only
alternative available
to the standard spectral sequence for equivariant
cohomology, as Hans-Werner Henn \cite{[HN]} has created a
``centralizer spectral sequence'' that involves centralizers of
elementary abelian $p$-subgroups;
however, Brown's spectral sequence appears
to be the easiest to apply in our situation.

Let $p$ be a prime,
$G$ be a group with finite virtual cohomological dimension,
$M$ be a $G$-module,
and
$\mathcal{A}$ be the poset of nontrivial elementary abelian
$p$-subgroups of $G$.
Let $G$ act on $\mathcal{A}$ by conjugation,
so that the simplicial realization $|\mathcal{A}|$ is
a G-complex.
Finally, denote by
$\hat H^*(G; M)_p$
the $p$-primary part of $\hat H^*(G; M)$.
Brown shows that there is an isomorphism
\begin{equation} \label{er3}
\hat H^*(G; M)_{p} \approx \hat H^*_G(\mathcal{A};
M)_{p}.
\end{equation}
Now let the coefficient module $M$ be either $\Z/p$ or
$\Z_{(p)}$ so that \ref{er3} becomes just
$$\hat H^*(G; M) \approx \hat H^*_G(\mathcal{A};M).$$
Then use the equivariant spectral sequence
\ref{e2} to compute $\hat H^*_G(\mathcal{A};M).$
This spectral sequence involves
stabilizers of simplices of $|\mathcal{A}|$.
If $$(A_0 \subset \cdots \subset A_r)$$ is an
$r$-simplex of $|\mathcal{A}|$ then its stabilizer is
the intersection of normalizers
$$\bigcap_{i=0}^r N_{G}(A_i)$$
because $G$ acts on $\mathcal{A}$ by conjugation.
Hence the equivariant spectral sequence
used to calculate $\hat H^*_G(\mathcal{A};M)$
becomes the following
``normalizer spectral sequence'':
\begin{equation} \label{er4}
E_1^{r,s} = \prod_{(A_0 \subset \cdots \subset A_r) \in |\mathcal{B}|_r}
\hat H^s( \bigcap_{i=0}^r N_{G}(A_i); M)
\Rightarrow \hat H^{r+s}(G; M)
\end{equation}
where $\mathcal{B}$ denotes the poset of conjugacy classes of
nontrivial elementary
abelian $p$-subgroups of $G$, and $|\mathcal{B}|_r$ is the
set of $r$-simplices in the realization $|\mathcal{B}|$.

For our case (i.e., the case $G=Aut(F_n)$),
we set $M=\Z_{(p)}$ and note that from
\cite{[H]} we know that $Aut(F_n)$ has finite vcd.
So we can apply \ref{er4} to obtain
\begin{equation} \label{e3}
E_1^{r,s} = \prod_{(A_0 \subset \cdots \subset A_r) \in |\mathcal{B}|_r}
\hat H^s( \bigcap_{i=0}^r N_{Aut(F_n)}(A_i); \Z_{(p)}) \Rightarrow \hat H^{r+s}(Aut(F_n); \Z_{(p)})
\end{equation}
where $\mathcal{B}$ denotes the poset of conjugacy classes of
nontrivial elementary
abelian $p$-subgroups of $Aut(F_n)$, and $|\mathcal{B}|_r$ is the
set of $r$-simplices in the realization $|\mathcal{B}|$.

It is important to note, however, that although we will not be using a
spectral sequence for equivariant cohomology directly to
compute the cohomology of $Aut(F_n)$, we will be computing the cohomology
groups of normalizers of elementary abelian $p$-subgroups
using spectral sequences for equivariant cohomology.
The output of these spectral sequences will then be entered in as
input into the $E_1$-page of the normalizer spectral sequence
\ref{e3}.

In this chapter, we will compute the elementary abelian
$p$-subgroups of $Aut(F_n)$.  In the next chapter, we will adapt
methods of Krstic to compute the normalizers of these
subgroups, in the chapter after that we will
define certain contractible subcomplexes of
auter space that these normalizers act on,
and in the following chapter we will use
the equivariant cohomology spectral sequence associated
to the action of the normalizers on these
contractible subcomplexes to compute
the cohomology of the normalizers.

From \cite{[B-T]} and \cite{[M]},
we see that $p^2$ is an upper bound for the order of any $p$-subgroup
of $Aut(F_n)$.  (Recall that we set $n=2(p-1)$ earlier.)
Every nontrivial
$p$-subgroup of $Aut(F_n)$ is isomorphic to either $\Z/p$ or
$\Z/p \times \Z/p$
(see, for example, Smillie and Vogtmann in \cite{[S-V2]}.)

A subgroup $G$ of $Aut(F_n)$ is said to be {\em realized} by
a specific marked graph $\eta : R_n \to \Gamma$
if it is contained in the stabilizer of that marked graph.
Since (from \cite{[S-V1]})
the stabilizer of $\eta : R_n \to \Gamma$ is isomorphic
to the group $Aut(\Gamma,*)$ of basepoint-preserving graph
automorphisms of $\Gamma$, this means that we can
think of $\Gamma$ as a graph with an action of $G$ on it.
If two subgroups
$G_1$ and $G_2$ are realized by marked graphs
$\eta_1: R_n \to \Gamma$ and $\eta_2: R_n \to \Gamma$,
respectively, and if both $G_1$ and $G_2$ give the same
action on the graph $\Gamma$, then $G_1$ and $G_2$
are conjugate.

We want to find all conjugacy classes of (elementary abelian)
$p$-subgroups of $Aut(F_n)$.
We can do this by analyzing marked
graphs with $p$-symmetry, so that the stabilizers of these marked
graphs have elements of order $p$.  Since we are only interested in
conjugacy classes of subgroups, we will just examine the underlying
graphs of the marked graphs.

We now define several $p$-subgroups
$A$, $B_k$, $C$, $D$, and $E$ of $Aut(F_n)$.  Our goal is to
show that these are a complete listing of 
the distinct conjugacy classes of
$p$-subgroups of $Aut(F_n)$.

\begin{itemize}
\item There is an action of $\Z/p$ on the rose $R_{2(p-1)}$ given by
rotating the first $p$ leaves of the rose.  By looking at the
stabilizer of a marked graph with underlying graph $R_{2(p-1)}$, this
action gives us a subgroup $A \cong \Z/p$ of $Aut(F_{2(p-1)})$.  This
subgroup is a maximal $p$-subgroup, in the sense that no other
$p$-subgroup properly contains it.
\item For each $k \in \{0, \ldots, p-1 \}$, there is an action of
$\Z/p$ on $R_k \vee \Theta_{p-1} \vee R_{p-1-k}$ given by rotating the
edges in $\Theta_{p-1}$.  This action gives us a subgroup
$B_k \cong \Z/p$ of $Aut(F_{2(p-1)})$.  If $k \in \{1, \ldots, p-2
\}$, then $B_k$ is a maximal $p$-subgroup.
\item There is an action of $\Z/p$ on $\Phi_{2(p-1)}$ which gives us a
(non-maximal) $p$-subgroup $C \cong \Z/p$ of $Aut(F_{2(p-1)})$.
\item There is an action of $\Z/p \times \Z/p$ on $\Omega_{2(p-1)}$
given by having the first $\Z/p$ rotate one of the $\Theta$-graphs in
$\Omega_{2(p-1)}$ and having the second $\Z/p$ rotate the other
$\Theta$-graph.  This action gives us a subgroup $D \cong \Z/p \times
\Z/p$ of $Aut(F_{2(p-1)})$.  This subgroup is maximal among
$p$-subgroups, and it contains $B_0$, $B_{p-1}$, and $C$.
\item There is an action of $\Z/p \times \Z/p$ on $\Psi_{2(p-1)}$
given by having the first $\Z/p$ rotate one of the $\Theta$-graphs in
$\Psi_{2(p-1)}$ and having the second $\Z/p$ rotate the other
$\Theta$-graph.  This action gives us a subgroup $E \cong \Z/p \times
\Z/p$ of $Aut(F_{2(p-1)})$.  This subgroup is maximal among
$p$-subgroups, and it contains $B_{p-1}$ and $C$.
\end{itemize}

\begin{prop} \label{t2}
Every $p$-subgroup of $Aut(F_{2(p-1)})$
is conjugate to one of
$$A, B_0, \ldots , B_{p-1}, C, D, \hbox{ or } E.$$
\end{prop}

\PF As we asserted earlier, every nontrivial 
$p$-subgroup $P$ is either $\Z/p$ or $\Z/p \times \Z/p$.
From Culler, this subgroup is realized by an action on
a basepointed marked graph $(\eta,\Gamma,*)$.  Since we are only
trying to calculate conjugacy classes of elementary
abelian $p$-subgroups, we need only think of $P$ as
acting on the underlying basepointed graph $\Gamma$.  In addition,
we might as well collapse invariant forests to get an action
of $P$ on a reduced graph $\Gamma$.

If $P=\Z/p$, then Proposition \ref{t1} gives us that $\Gamma$ is one of
$R_n$, $R_k \vee \Theta_{p-1} \vee R_{p-1-k}$,
$\Psi_n$, or $\Omega_n$.  If $\Gamma$ is $R_n$ then $P$ is conjugate
to $A$.  If $\Gamma$ is $R_k \vee \Theta_{p-1} \vee R_{p-1-k}$ then
$P$ is conjugate to $B_k$.  Finally, note that
by collapsing different invariant forests the action
of $\Z/p$ on $\Phi_n$ gives a diagonal action of $\Z/p$ on
both $\Psi_n$ and $\Omega_n$.  Hence if $\Gamma$ is either
$\Psi_n$ or $\Omega_n$ then $P$ is conjugate to $C$.

Next, suppose $P=\Z/p \times \Z/p = (\alpha) \times (\beta)$.
The first cyclic summand must rotate $p$ edges
$e_1, e_2, \ldots, e_p$ of $\Gamma$.  Without loss of
generality, we may assume
$(\beta)$ fixes these edges
(by replacing
$\beta$ with $\beta - \alpha^j$ for some $j$
if necessary.)  Now $\beta$
must rotate $p$ other edges $f_1, f_2, \ldots, f_p$.  Again,
we can also assume without loss of generality that $\alpha$
fixes the edges $f_i$.  From Lemma \ref{tr1},
the $e_i$ and $f_i$
must form either roses or $\Theta$-graphs.  Neither can form
a rose or the rank of $\pi_1(\Gamma)$ is at least $2p-1$.  So both
form $\Theta$-graphs.  Hence $\Gamma$ is either $\Psi_n$ or
$\Omega_n$.  In the former case $P$ is conjugate to $E$ while
in the latter case it is conjugate to $D$. \END

In order to determine whether or not two different graphs
give us the same elementary abelian $p$-subgroup,
we need to use work of Krstic in \cite{[K]} involving
{\em Nielsen transformations.}  More will be said on this
in the next chapter, but for now we briefly 
recall Krstic's definition of
Nielsen transformations and state the result of
Krstic that we need.

\begin{defn}[Nielsen transformation] \label{tr2}
Let $G$ be a finite subgroup of $Aut(F_n)$, which is realized by an
action of $G$ on a reduced, basepointed graph $\Gamma$.  
Let $V$ and $E$ be the vertex and oriented edge sets,
respectively, of $\Gamma$.  Finally, let $\iota, \tau : E \to V$ be
maps which give the initial and terminal points, respectively, of
oriented edges.
If there are two edges $e$ and $f$ of
$\Gamma$ such that:
\begin{itemize}
\item $e$ and $f$ are in different orbits (i.e., $f \not
\in Ge \cup G\bar e$),
\item $\tau e = \tau f$, and
\item $stab(e) \subseteq stab(f)$.
\end{itemize}
then there is an \emph{admissible Nielsen transformation}
$<e,f>$ from $\Gamma$ to
a new graph $<e,f>\Gamma$.
The graph $\Gamma' = <e,f>\Gamma$ has the same vertex and
edge sets $V$ and $E$ as $\Gamma$; however, the map $\tau' : E \to V$
which gives the terminal point of an edge is changed as follows.  For
edges $h$ not in the orbit of $e$, set $\tau ' h=\tau h$; but set
$\tau' ge = \iota g f$ for $g \in G$.
\end{defn}

If a sequence of Nielsen transformations can change
a $G$-graph $\Gamma_1$ into a $G$-graph $\Gamma_2$,
then $\Gamma_1$ is said to be {\em Nielsen equivalent} to $\Gamma_2$.

We need the following theorem of Krstic, Theorem 2 from
\cite{[K]}:

\begin{thm}[Krstic] \label{tr4}
Let $\Gamma_1$ and $\Gamma_2$ realize the same subgroup
$G$ of
$Aut(F_n)$.  If $\Gamma_1$ and $\Gamma_2$ are reduced as
$G$-graphs then they are Nielsen equivalent,
up to an equivariant isomorphism (a basepoint preserving
isomorphism.)
\end{thm}

\begin{prop} \label{tr5}
The subgroups
$A$, $B_0, \ldots , B_{p-1}$, $C$, $D$, $E$
are in distinct conjugacy classes.  The diagram of
subgroups up to conjugacy is
\begin{equation} \label{e4}
\matrix{ & & B_0 & & \cr
           & & \downarrow & & \cr
           & & D & & \cr
           & \nearrow & & \nwarrow & \cr
           B_{p-1} & & & & C \cr
           & \searrow & & \swarrow & \cr
           & & E & & \cr}
\end{equation}
\end{prop}

\PF We apply Theorem \ref{tr4}.
Any two reduced graphs realized by the same subgroup
of $Aut(F_n)$ can be connected up by a sequence of Nielsen
transformations.  It is easy to see that the only case
in which this could occur for two different graphs listed in
Proposition \ref{t1} is when $P=\Z/p$ and the two graphs are
$\Psi_n$ and $\Omega_n$.  We already know, however, that in
this case $P$ is conjugate to the subgroup $C$ of $Aut(F_n)$. \END

\chapter{Normalizers of $p$-subgroups of $Aut(F_n)$} \label{c3}

The structure of centralizers of finite subgroups of $Aut(F_n)$ is
given by Krstic in \cite{[K]} in some detail.
He shows that, for a finite subgroup $G$ of $Aut(F_n)$, an element
of the centralizer $C_{Aut(F_n)}(G)$ is a product of
Nielsen transformations followed by a centralizing
graph isomorphism. Moreover, one of the
main propositions of his paper can be used to yield information about
the structure of normalizers of finite subgroups, which is what
we are interested in here.

We now recall some definitions and theorems from Krstic \cite{[K]}.

\begin{defn}[Nielsen isomorphism] \label{tr3}
Let $<e,f>$ be a Nielsen transformation from $\Gamma$ to $\Gamma'$
(see Definition \ref{tr2}.)
A \emph{Nielsen
isomorphism} is the isomorphism between fundamental groupoids:
$$<e,f> : \Pi(\Gamma) \to \Pi(\Gamma')$$
given by $<e,f>h=h$ if h is not in the orbit of $e$ and
$<e,f>ge=gef$ for $g \in G$.
\end{defn}

Sometimes
this notation is abused and we say
Nielsen transformation when we mean Nielsen
isomorphism.  This abuse will tend to happen when
we are implicitly dealing with maps
between fundamental groupoids even though we only
state things on the level of
deforming graphs by applying Nielsen transformations
to them.

The following is Proposition $4'$ from \cite{[K]}:

\begin{prop}[Krstic] \label{tr6}
If $\Gamma_1$ and
$\Gamma_2$ are finite basepointed, reduced
 $G$-graphs and
$$F : \Pi(\Gamma_1) \to \Pi(\Gamma_2)$$
is an equivariant groupoid isomorphism
which preserves the base vertex, then there
exists a product $T$ of Nielsen transformations and an equivariant
isomorphism of basepointed graphs
$$H: T\Gamma_1 \to \Gamma_2$$ such that
$F=HT$.
\end{prop}

Recall that $\pi_1(\Gamma_i)$ is a sub-groupoid of
$\Pi(\Gamma_i)$. Krstic uses the above theorem to
get information about maps between fundamental groups.
He shows the following, which is Corollary $1'$ in \cite{[K]}:

\begin{cor}[Krstic] \label{tr7}
Let $\Gamma_1$ and $\Gamma_2$ be reduced pointed $G$-graphs of
rank $\geq 2$.  Then every equivariant isomorphism
$$\pi_1(\Gamma_1,*) \to \pi_1(\Gamma_2,*)$$
is the restriction of an equivariant isomorphism
$$\Pi(\Gamma_1) \to \Pi(\Gamma_2).$$
\end{cor}

For a $G$-graph $\Gamma$ and an automorphism $\phi : G \to G$, let
$\phi(\Gamma)$ be the $G$-graph with underlying graph $\Gamma$ and
$G$-action given by $gx := \phi(g)x$ where $g \in P$, $x \in V(\Gamma)
\cup E(\Gamma)$, and the latter multiplication is given by the
standard action on $\Gamma$.  In other words, $\phi(\Gamma)$
is just the graph $\Gamma$ with the $G$-action ``twisted''
by $\phi$.

Realize the finite subgroup $G$ of $Aut(F_n)$ by a
marked graph $\eta : R_n \to \Gamma$ where
the induced action of $G$ on $\Gamma$
is reduced.

\begin{prop} \label{t3}
Every element $\alpha$ of
$N_{Aut(F_n)}(G)$ is realized by
a $G$-equivariant isomorphism
$$\pi_1(\Gamma,*) \to \pi_1(\phi(\Gamma),*)$$
for some automorphism $\phi$ of $G$.
\end{prop}

Note that the automorphism $\phi$ can change
as we consider different elements of $N_{Aut(F_n)}(G)$.

\PF
Define $N_{Aut(F_n)}'(G)$ to be the set of all
equivalence classes of pairs
$(\phi,\psi)$
where $\phi \in Aut(G)$ and
$\psi : \Gamma \to \Gamma$ is a basepoint preserving, continuous
surjection of graphs such that
$$\psi_\# : \pi_1(\Gamma,*) \to \pi_1(\phi(\Gamma),*)$$
is a $G$-equivariant group isomorphism.
Two such pairs $(\phi_1,\psi_1)$ and $(\phi_2,\psi_2)$
are equivalent if $\phi_1=\phi_2$ and
$$(\psi_1)_\# = (\psi_2)_\# : \pi_1(\Gamma,*) \to \pi_1(\phi(\Gamma),*)$$
on the level of fundamental groups.

Define a group operation {\em composition} in
$N_{Aut(F_n)}'(G)$ in the obvious way by letting
$$(\phi_2,\psi_2) \circ (\phi_1,\psi_1) =
(\phi_2 \circ \phi_1, \psi_2 \circ \psi_1).$$

To prove the proposition, it suffices to show that
there is a group isomorphism
$$\xi : N_{Aut(F_n)}(G) \to N_{Aut(F_n)}'(G).$$

Choose a
fixed homotopy inverse $\bar \eta : \Gamma \to R_n$; i.e., a map such
that $$(\eta \bar \eta)_\# = 1 \in Aut(\pi_1(\Gamma))$$ and
$$(\bar \eta \eta)_\# = 1 \in Aut(\pi_1(R_n)).$$

\begin{enumerate}
\item In this step, we define the map $\xi$ and show that it is
a well defined map between sets.  (We delay showing that is
a well defined map between groups, i.e. a homomorphism,
until step 3.)

If $\alpha \in N_{Aut(F_n)}(G)$, then $\alpha$ induces an
automorphism $\phi$ of $G$ via conjugation: $\phi(g)=\alpha g
\alpha^{-1}$.  In addition, since $\alpha \in Aut(F_n)$ it
corresponds to a
map $\tilde \alpha : R_n \to R_n$ which gives a marked graph $\alpha
\eta$ by precomposition:
$$R_n \buildrel \tilde \alpha \over \to R_n
\buildrel \eta \over \to \Gamma.$$
Similarly, $\alpha$ induces a map $\psi = \eta \tilde \alpha
\bar \eta$:
$$\Gamma \buildrel \bar \eta \over \to
R_n \buildrel \tilde \alpha \over \to
R_n \buildrel \eta \over \to \Gamma$$
which in turn gives a map
$$\psi_\# : \pi_1(\Gamma,*) \to \pi_1(\Gamma,*).$$

It remains to show that $\psi_\#$ is
an equivariant map in the suitable sense.
Now for any $x \in \pi_1(\Gamma,*)$ and $g \in G$,
$$\matrix{\hfill \psi_\#(gx) &= (\psi_\# g \psi_\#^{-1})\psi_\#(x) \hfill \cr
\hfill &= (\eta \tilde \alpha g \tilde \alpha^{-1} \bar \eta)_\# \psi_\#(x) \hfill \cr
\hfill &= (\eta \phi(g) \bar \eta)_\# \psi_\#(x) \hfill \cr
\hfill &= \phi(g) \psi_\#(x). \hfill \cr }$$

At the risk of being overly pedantic, we explain the above
equalities in detail.
Both the first and third equalities are trivial.
 For the second equality, one can reason as follows.
As $\Gamma$ is a $G$-graph, one at first thinks of
$g$ as a map from $\Gamma$ to $\Gamma$.  If one wants to
translate this to thinking of $g$ as an element of
$Aut(F_n)$ realized on the rose, one forms
$$R_n \buildrel \eta \over \to
\Gamma \buildrel g \over \to
\Gamma \buildrel \bar \eta \over \to R_n.$$
We know that this is the correct way to think of $g$
as realized on the rose, because if we precompose
the marked graph $(\eta,\Gamma)$ with the above element
then we get $(\eta g \bar \eta)(\eta)=\eta g$, or the
result of the graph automorphism $g$ acting on
the marked graph.
So
$$\matrix{
\hfill (\psi_\# g \psi_\#^{-1}) &= (\psi g \psi^{-1})_\# \hfill \cr
\hfill &= (\eta \tilde \alpha \bar \eta g \eta \tilde \alpha^{-1} \bar \eta)_\# \hfill \cr
\hfill &= (\eta \tilde \alpha  g \tilde \alpha^{-1} \bar \eta)_\#. \hfill \cr}$$
This justifies the second equality.

For the fourth equality, note that
$$R_n \buildrel \phi(g) \over \to
R_n \buildrel \eta \over \to \Gamma$$
is the same as
$$R_n \buildrel \eta \over \to
\Gamma \buildrel \phi(g) \over \to \Gamma$$
where in the first case $\phi(g) \in G$
is thought of as an element of $Aut(F_n)$
realized on the rose and in the second case
$\phi(g) \in G$ is thought of as
giving a graph automorphism of the
$G$-graph $\Gamma$.

Thus we have shown that if we think of the map $\psi_\#$
as
$$\psi_\# : \pi_1(\Gamma,*) \to \pi_1(\Gamma,*)$$
then
$\psi_\#(gx) = \phi_(g) \psi_\#(x).$
So if instead we think of $\psi_\#$ as a map
$$\psi_\# : \pi_1(\Gamma,*) \to \pi_1(\phi(\Gamma),*)$$
then it is $G$-equivariant.

\item In this step, we define a set-theoretic inverse map
$$\xi^{-1} : N_{Aut(F_n)}'(G) \to N_{Aut(F_n)}(G).$$

Now suppose we are given a pair $(\phi,\psi)$.  Then $\psi$
gives us an element $\tilde \alpha : R_n \to R_n$ via
$\bar \eta \psi \eta$:
$$R_n \buildrel \eta \over \to
\Gamma \buildrel \psi \over \to
\Gamma \buildrel \bar \eta \over \to R_n.$$
The element $\tilde \alpha$ induces a map
$\alpha \in Aut(F_n)$.

Recall that the the map
$$\psi_\# : \pi_1(\Gamma,*) \to \pi_1(\phi(\Gamma),*)$$
is $G$-equivariant.  If we instead think of it as a
map
$$\psi_\# : \pi_1(\Gamma,*) \to \pi_1(\Gamma,*)$$
then $$\psi_\#(gx) = \phi(g) \psi_\#(x)$$ for all
$g \in G$ and $x \in \pi_1(\Gamma,*)$.
So $$gx = \psi_\#^{-1} \phi(g) \psi_\#(x)$$
for all $x \in \pi_1(\Gamma,*)$.  Now similar arguments to
those used in $1.$ yield that $g=\alpha^{-1} \phi(g) \alpha$.
Consequently $$\alpha g \alpha^{-1} = \phi(G) \in G$$
for all $g \in G$.  This implies that
$\alpha \in N_{Aut(F_n)}(G)$.

That the map we are calling
$\xi^{-1}$ actually is a set-theoretic inverse
to $\xi$ can easily be verified by looking at their
definitions.  For example, if $(\phi,\psi) \in N_{Aut(F_n)}'(G)$
then
$$\matrix{
\hfill \xi(\xi^{-1}(\phi,\psi)) &= \xi( (\bar \eta \psi \eta)_\# ) \hfill \cr
\hfill &= (g \mapsto (\bar \eta \psi \eta)_\# g (\bar \eta \psi \eta)_\#^{-1},
\eta (\bar \eta \psi \eta) \bar \eta \hfill) \cr
\hfill &= (g \mapsto \phi(g), \psi). \hfill \cr}$$
Hence $\xi \xi^{-1} = 1$.
The equality $\xi^{-1} \xi = 1$ can
be shown similarly.

\item Finally, note that $\xi$ as defined is clearly a
group homomorphism because
$$\matrix{
\hfill \xi(\alpha_2) \circ \xi(\alpha_1) &=
(g \mapsto \alpha_2 g \alpha_2^{-1}, \eta \tilde \alpha_2 \bar \eta)
\circ (g \mapsto \alpha_1 g \alpha_1^{-1}, \eta \tilde \alpha_1 \bar \eta) \hfill \cr
\hfill &= (g \mapsto (\alpha_2 \circ \alpha_1) g (\alpha_2 \circ \alpha_1)^{-1},
\eta \tilde \alpha_2 \circ \tilde \alpha_1 \bar \eta) \hfill \cr
\hfill &= \xi(\alpha_2 \circ \alpha_1). \square \hfill \cr}$$
\end{enumerate}

As a result of Proposition \ref{t3}, we obtain the following:

\begin{prop} \label{t4} An element of $N_{Aut(F_n)}(G)$ is
a product of Nielsen transformations followed by a normalizing
graph isomorphism.
\end{prop}

\PF Any element $(\phi, \psi)$ of
$N_{Aut(F_n)}'(G)$ induces a $G$-equivariant map
$$\pi_1(\Gamma,*) \buildrel \psi_\# \over \to
\pi_1(\phi(\Gamma),*).$$
From Corollary \ref{tr7}, $\psi_\#$ is the restriction of
a $G$-equivariant map $F$ between fundamental groupoids:
$$\Pi(\Gamma) \buildrel F \over \to \Pi(\phi(\Gamma)).$$
Now from Proposition \ref{tr6} we know that there
exists a product $T$ of Nielsen transformations
starting with the graph $\Gamma$, and a
basepoint preserving graph isomorphism
$H : T\Gamma \to \phi(\Gamma)$ such that $F$ is the map
induced by $HT$: 
$$F: \Pi(\Gamma) \buildrel T \over \to
\Pi(T\Gamma) \buildrel H \over \to \Pi(\phi(\Gamma)).$$
By restricting the above map of fundamental groupoids to
one of fundamental groups, we have that $\psi_\#$ is
also the map induced by $HT$:
$$\psi_\#: \pi_1(\Gamma,*) \buildrel T \over \to
\pi_1(T\Gamma,*) \buildrel H \over \to
\pi_1(\phi(\Gamma),*).$$
It is in the above sense that we mean an element of
$N_{Aut(F_n)}(G)$ is a product of Nielsen transformations followed
by a normalizing graph isomorphism.
\END

\bigskip

We now use Proposition \ref{t4} to calculate the normalizers
of the subgroups
$A$, $B_0, \ldots , B_{p-1}$, $D$, and $E$
listed in Proposition \ref{tr5}.  We do this so that their
cohomology groups will be easier to calculate in a later chapter.
We will not need
to calculate the normalizer of $C$ in this explicit manner,
and later will find its cohomology through more geometric means.

\begin{lemma} \label{tr8}
$$N_{Aut(F_n)}(A) \cong N_{\Sigma_p}(\Z/p) \times \bigl( (F_{p-2} \times F_{p-2})
\rtimes (\Z/2 \times Aut(F_{p-2}))  \bigr)$$
where $\Z/2$ acts by exchanging the two copies of $F_{p-2}$ and
$Aut(F_{p-2})$ acts diagonally on the two copies of $F_{p-2}$.
\end{lemma}

\PF Recall that $A$ comes from the action of $\Z/p$ on the first $p$
petals of the rose $R_{2(p-1)}$.  From Proposition \ref{t4}
any element of $N_{Aut(F_n)}(A)$ is induced by a
product of Nielsen transformations $T$ followed by a
graph isomorphism
$$H: TR_{2(p-1)} \to \phi(R_{2(p-1)})$$
for some $\phi \in Aut(A)$.  The only types of Nielsen transformations
that are possible are
\begin{itemize}
\item Nielsen transformations obtained by pulling
either the front ends or the back ends of the first $p$ petals
of the rose uniformly around paths in the last $p-2$ petals.  This
subgroup is isomorphic to $F_{p-2} \times F_{p-2}$.  The first copy
of $F_{p-2}$ comes from pulling the front ends of the first $p$
petals, and the second copy comes from pulling the back ends of the
first $p$ petals.
\item Nielsen transformations involving just the last $p-2$ petals
of the rose.  In other words, Nielsen transformations contained
entirely in the copy of $Aut(F_{p-2})$ corresponding to
graph automorphisms and Nielsen transformations involving the
last $p-2$ petals.
\end{itemize}

Note that for a single Nielsen transformation $T$ in either of the above
cases, the $A$-graph $TR_{2(p-1)}$ is exactly the same
as the $A$-graph $R_{2(p-1)}$.  It follows that for any product
of Nielsen transformations $T$ the $A$-graph $TR_{2(p-1)}$ is exactly
the same as the $A$-graph $R_{2(p-1)}$ (although, of course,
the product of Nielsen isomorphisms induced by the product of Nielsen
transformations will probably not be trivial.)
Hence from Proposition \ref{t4}, we
know that
any element of $N_{Aut(F_n)}(A)$ is induced by a
product of Nielsen transformations $T$ followed by a
graph isomorphism
$$H: TR_{2(p-1)} = R_{2(p-1)} \to \phi(R_{2(p-1)})$$
for some $\phi \in Aut(A)$.
Thus the only ``normalizing graph isomorphisms'' that we need
to consider are ones which go from $R_{2(p-1)}$ to $\phi(R_{2(p-1)})$
for some $\phi \in Aut(A)$.  In other words, we need only examine
automorphisms of one particular graph $R_{2(p-1)}$ and see which are
in the normalizer $N_{Aut(F_n)}(A)$.

The normalizing graph automorphisms are:
\begin{itemize}
\item Those involving just the last $p-2$ petals of the rose.
As the action of $A$ on these petals is trivial,
the normalizer of $A$ contains any graph automorphism
involving just those last $p-2$ petals.  All of these graph automorphisms
are contained in 
the copy of $Aut(F_{p-2})$ obtained from graph automorphisms and
Nielsen transformations involving the last $p-2$ petals.
\item Normalizing graph automorphisms of the first $p$
petals of the rose.  This gives a subgroup of $N_{Aut(F_n)}(A)$ that is
isomorphic to $(\Z/p \rtimes \Z/(p-1)) \times \Z/2$ where the $\Z/2$
comes from inverting all of the $p$-petals at the same time, and
$$(\Z/p \rtimes \Z/(p-1)) \cong N_{\Sigma_p}(\Z/p)
= N_{\Sigma_p}((1 2 \ldots p)).$$
In other words, the generator of $\Z/(p-1)$ acts on the generator of
$\Z/p$ by conjugating it to its $s$-th power for some generator $s$ of
$\F_p^\times$ and
$$N_{\Sigma_p}(\Z/p) \cong
\langle g,h | g^p=1, h^{p-1}=1, hgh^{-1} = g^s \rangle$$
where $s$ is a primitive root modulo $p$.

We can verify in
a very concrete ``nuts and bolts'' manner
that the normalizing graph automorphisms of the
first $p$ petals of the rose take the above form.  We do this now
for the sake of providing at least one concrete example.

The full automorphism
group of the first $p$ petals is
$$(\Z/2)^p \rtimes \Sigma_p$$
where $\Sigma_p$ acts on $(\Z/2)^p$ by permuting the
$p$ copies of $\Z/2$.  Each of the $p$ copies of $\Z/2$ corresponds to
flipping one of the first $p$ petals of the rose.  The
$\Sigma_p$ corresponds to permuting the first $p$ petals
of the rose among themselves.  The subgroup $A$
corresponds to the subgroup generated by the permutation
$(1 2 \ldots p)$ in $\Sigma_p$ and rotates the
first $p$ petals sequentially.
We must see which elements of
$$(\Z/2)^p \rtimes \Sigma_p$$
are in the normalizer of $A$.  Say that the element
$$(a_1, \ldots, a_p) \cdot b$$
is in the normalizer, where
$(a_1, \ldots, a_p) \in (\Z/2)^p$
and $b \in \Sigma_p$.
Then
$$(a_1, \ldots, a_p) \cdot b \cdot (1 2 \ldots p) \cdot b^{-1}
\cdot (a_1, \ldots, a_p) = (1 2 \ldots p)^r$$
for some $r \in \{0, 1, \ldots, p-1\}$.
Now $$b \cdot (1 2 \ldots p) \cdot b^{-1}$$ is a
permutation of order $p$ in $\Sigma_p$ and so has the form
$$(c_1 c_2 \ldots c_p)$$
where the integers $c_i$ are all distinct.
Accordingly,
$$(a_1, \ldots, a_p) \cdot (c_1 c_2 \ldots c_p)
\cdot (a_1, \ldots, a_p) = (1 2 \ldots p)^r.$$
Recall that each $a_i \in \Z/2$.  We want to
show that if some $a_j=1$ then all of the $a_i$ are $1$.
Now if $a_j=1$, the graph automorphism
$$(a_1, \ldots, a_p) \cdot (c_1 c_2 \ldots c_p)
\cdot (a_1, \ldots, a_p)$$
first flips the $j$th petal, then sends that
flipped petal to the $(c_1 c_2 \ldots c_p)j$ petal, and then
finally flips that last petal back.  It must end by flipping
the $(c_1 c_2 \ldots c_p)j$ petal back because the graph automorphism
$(1 2 \ldots p)^r$ does not flip any petals.
So $a_{(c_1 c_2 \ldots c_p)j}$ is also equal to $1$.
We then iterate the above reasoning to get that
$a_{(c_1 c_2 \ldots c_p)^{2}j}$ equals $1$, etc.
Hence if some $a_j=1$ then all of the $a_i$ are $1$.
In other words, $$a_1 = a_2 = \ldots = a_p.$$
It now follows that
$$(a_1, \ldots, a_p) \cdot (c_1 c_2 \ldots c_p)
\cdot (a_1, \ldots, a_p) = (c_1 c_2 \ldots c_p)$$
and so
$$b \cdot (1 2 \ldots p) \cdot b^{-1} = (1 2 \ldots p)^r.$$
Thus $b$ is in $N_{\Sigma_p}((1 2 \ldots p))$
which is (see, for example, \cite{[B1]} page 68)
$$\langle g,h | g^p=1, h^{p-1}=1, hgh^{-1} = g^s \rangle.$$

Hence the normalizing graph automorphisms of the first $p$
petals form a subgroup of $N_{Aut(F_n)}(A)$ which is
isomorphic to
$$(\Z/p \rtimes \Z/(p-1)) \times \Z/2$$
as we claimed.
\end{itemize}

Putting all of this together, we see that
$$N_{Aut(F_n)}(A) \cong N_{\Sigma_p}(\Z/p) \times \bigl( (F_{p-2} \times F_{p-2})
\rtimes (\Z/2 \times Aut(F_{p-2}))  \bigr)$$
where $\Z/2$ acts by exchanging the two copies of $F_{p-2}$ and
$Aut(F_{p-2})$ acts diagonally on the two copies of $F_{p-2}$.
The $(F_{p-2} \times F_{p-2})$ corresponds to Nielsen
transformations obtained by pulling the
front or back ends of the first $p$ petals
uniformly around some path in the last $p-2$ petals;
the $\Z/2$ corresponds to flipping all of the first $p$
petals over at the same time; and the
$Aut(F_{p-2})$ corresponds to Nielsen transformations and
graph automorphisms involving the last $p-2$ petals.
We leave it as an exercise for the reader to show that
the various subgroups of $N_{Aut(F_n)}(A)$ fit together as described above. \END

A final remark about the structure
of the subgroup $N_{Aut(F_n)}(A)$, described
above in the proof of Lemma \ref{tr8},
is appropriate here:

\begin{remark} \label{tr9}
Consider a subgroup $\langle \omega \rangle \cong \Z/2$ of $Aut(F_{p})$
corresponding to the action of $\Z/2$ on $R_p$
given by switching the
just the first two petals of the
rose.  Note that
$$N_{Aut(F_{p})}(\omega) = C_{Aut(F_{p})}(\omega)=
\Z/2 \times \bigl( (F_{p-2} \times F_{p-2})
\rtimes (\Z/2 \times Aut(F_{p-2}))  \bigr)$$
where the action of $\Z/2 \times Aut(F_{p-2})$
on $F_{p-2} \times F_{p-2}$ in the semidirect
product is the same as before (that is, as in $N_{Aut(F_n)}(A)$.)
Consequently, we see that
$$N_{\Sigma_p}(\Z/p) \times N_{Aut(F_p)}(\omega) \cong \Z/2 \times N_{Aut(F_n)}(A)$$
and hence
$$\hat H^*(N_{\Sigma_p}(\Z/p) \times N_{Aut(F_p)}(\omega); \Z_{(p)})
= \hat H^*(N_{Aut(F_n)}(A); \Z_{(p)})$$
because $p \geq 3$ and so the first summand
$\Z/2$ in $\Z/2 \times N_{Aut(F_n)}(A)$ can be ignored
when taking Farrell cohomology with $\Z_{(p)}$
coefficients.
\end{remark}

\begin{lemma} \label{tr10}
$$N_{Aut(F_n)}(B_k) \cong N_{\Sigma_p}(\Z/p) \times (F_k \rtimes Aut(F_k))
\times Aut(F_{p-1-k})$$
\end{lemma}

\PF Recall that $B_k$ comes from the action of $\Z/p$ on the $\Theta$-graph in
$R_k \vee \Theta_{p-1} \vee R_{p-1-k}$.
For just this proof,
denote the middle $p$ edges of the $\Theta$ graph by $e_1, \ldots, e_p$,
the $k$ petals of the rose $R_k$ by $f_1, \ldots, f_k$,
and the $p-1-k$ petals of the rose $R_{p-1-k}$
by $g_1, \ldots, g_{p-1-k}$.  Say that the edges $e_i$ are oriented
such that they go from the rose $R_{p-1-k}$ to the rose
$R_k$.
Also, just for this proof, let $\Gamma$ denote the
$B_k$-graph $R_k \vee \Theta_{p-1} \vee R_{p-1-k}$.

The normalizer
$N_{Aut(F_n)}(B_k)$ contains four types of Nielsen
transformations:
\begin{itemize}
\item Ones obtained from Nielsen transformations involving
just the $k$ petals of the rose $R_k$.
As the action of $B_k$ on this rose is trivial, any Nielsen
transformation of those petals is admissible.
These are all contained
in a subgroup $Aut(F_k)$ of the normalizer obtained from
taking graph automorphisms and Nielsen transformations of
the rose $R_k$.
\item Ones obtained from Nielsen transformations involving
just the $p-1-k$ petals of the rose $R_{p-1-k}$.
The action of $B_k$ on this rose is trivial.
These Nielsen transformations are all contained
in a subgroup $Aut(F_{p-1-k})$ of the normalizer obtained from
taking graph automorphisms and Nielsen transformations of
the rose $R_{p-1-k}$.
\item Nielsen transformations obtained by pulling the back edges of the
$\Theta$-graph uniformly around the $p-1-k$ petals of the rose
$R_{p-1-k}$ on the right.  Note that the Nielsen transformations
around the last $p-1-k$ petals induce 
maps that are similar to conjugation maps
of the last $p-1-k$ petals of the rose,
on the level of fundamental groups
$$\pi_1(R_k \vee \Theta_{p-1} \vee R_{p-1-k},*)
\to \pi_1(R_k \vee \Theta_{p-1} \vee R_{p-1-k},*).$$
That is, if we pull the back ends of all of the $e_i$
uniformly around some $g_j$, then on fundamental groups
the result is:
\begin{enumerate}
\item The loops $\bar e_{i_1} e_{i_2}$ are all sent back to themselves,
for $i_1, i_2 \in \{1, \ldots, p\}$.
This is because $\bar e_{i_1} e_{i_2}$ is sent to the
loop $(\bar e_{i_1} g_j) (\bar g_j e_{i_2}) = \bar e_{i_1} e_{i_2}$.
\item The loops $f_i$ are all sent back to themselves.
\item The loops $\bar e_r g_i e_s$ are sent to
$\bar e_r g_j g_i \bar g_j e_s$ for all
$r,s \in \{1, \ldots, p\}$, $i \in \{1, \ldots, p-1-k\}$.
In other words, paths in the rose $R_{p-1-k}$
are ``conjugated'' by $g_j$, where the word conjugated is
in quotes because we have first to get to the
path in the rose $R_{p-1-k}$ by travelling along some $\bar e_r$ and
we have to leave the the path in the
rose by travelling back along some $e_s$.
\end{enumerate}
We can thus disregard
these Nielsen transformations as they are already included in
the Nielsen transformations involving
just the rose $R_{p-1-k}$, which was covered in the previous case.
\item Nielsen transformations obtained by pulling the front edges of the
$\Theta$-graph uniformly around the $k$ petals of the rose
$R_k$ on the left.
Say we perform the Nielsen transformation given by
pulling the front ends of all of the $e_i$ around the edge $f_j$.
Then the resulting graph is the same as the one we started with, but the
Nielsen isomorphism on fundamental groupoids sends
$$\matrix{
e_i & \mapsto & e_i f_j & \hbox{ for } i \in \{1, \ldots, p\} \cr
f_i & \mapsto & f_i & \hbox{ for } i \in \{1, \ldots, k\} \cr
g_i & \mapsto & g_i & \hbox{ for } i \in \{1, \ldots, p-1-k\} \cr
}$$
On the level of fundamental groups, this has the result of
conjugating paths that start with some $\bar e_{i_1}$ and end with some
$e_{i_2}$, by the letter $f_j$.  Direct examination
reveals that none of the Nielsen transformations covered in the
previous cases do this.
We could also take the inverse
$<e_1,\bar f_j>$ of the above Nielsen
transformation, namely pull the edges $e_i$ uniformly around some edge
$\bar f_j$.
This would have the result of
conjugating paths that start with some $\bar e_{i_1}$ and end with some
$e_{i_2}$, by the letter $\bar f_j$.
Hence all of these Nielsen transformations 
form a subgroup
of the normalizer that is isomorphic to $F_k$,
where the $j$th generator of the free group $F_k$ corresponds
to the Nielsen transformation $<e_1, f_j>$.
Moreover,
this subgroup $F_k$ is acted upon by
the subgroup $Aut(F_k)$ of $N_{Aut(F_n)}(B_k)$ corresponding to 
graph automorphisms and Nielsen transformations of the 
first $k$ petals.  You can see this by noting how the subgroup
$Aut(F_k)$ affects the above conjugation maps.
Also note that the subgroup $Aut(F_{p-1-k})$ of
$N_{Aut(F_n)}(B_k)$, corresponding to 
graph automorphisms and Nielsen transformations of the 
last $p-1-k$ petals, commutes with the above Nielsen
transformations.
\end{itemize}

As was the case for Lemma \ref{tr8} above,
for any product $T$ of any of the above
types of Nielsen transformations,
the $B_k$-graph $T\Gamma$ is exactly the
same as the $B_k$-graph $\Gamma$
(although, again,
the product of Nielsen isomorphisms induced by the product of Nielsen
transformations will probably not be trivial.)
From Proposition \ref{t4}, we
know that
any element of $N_{Aut(F_n)}(B_k)$ is induced by a
product of Nielsen transformations $T$ followed by a
graph isomorphism
$$H: T\Gamma = \Gamma \to \phi(\Gamma)$$
for some $\phi \in Aut(B_k)$.
The only ``normalizing graph isomorphisms'' that we need
to consider are ones which go from $\Gamma$ to $\phi(\Gamma)$
for some $\phi \in Aut(B_k)$.

The normalizing graph automorphisms take one of three forms:

\begin{itemize}
\item Those automorphisms involved only with the
petals of the rose $R_k$.  These are contained in the
subgroup $Aut(F_k)$ 
obtained from graph automorphisms and Nielsen transformations
involving just the rose $R_k$.
\item Those automorphisms involved only with the
petals of the rose $R_{p-1-k}$.  These are contained in the
subgroup $Aut(F_{p-1-k})$ 
obtained from graph automorphisms and Nielsen transformations
involving just the rose $R_{p-1-k}$.
\item Normalizing graph automorphisms of the $p$ edges $e_i$ in
$\Theta_{p-1}$.  A graph automorphism can only permute the
edges $e_i$ among themselves (and can not, of course, flip
them over because of where the basepoint is located.)  Hence we can identify
the graph automorphisms which move the edges $e_i$ with the
symmetric group $\Sigma_p$.  The group $B_k$ corresponds
to the subgroup of $\Sigma_p$ generated by the
permutation $(1 2 \ldots p)$.  In addition, the
graph automorphisms of $\Sigma_p$ which are
normalizing are exactly those which are in the normalizer
of the subgroup generated by $(1 2 \ldots p)$.
This gives a subgroup of $N_{Aut(F_n)}(B_k)$ that is
isomorphic to
$$(\Z/p \rtimes \Z/(p-1)) \cong N_{\Sigma_p}(\Z/p).$$
\end{itemize}

Putting all of this together,
it is not a difficult exercise for the reader to
verify that
$$N_{Aut(F_n)}(B_k) \cong N_{\Sigma_p}(\Z/p) \times (F_k \rtimes Aut(F_k))
\times Aut(F_{p-1-k}).$$
We have already remarked that $Aut(F_k)$ acts on $F_k$
as described and that \newline
$Aut(F_{p-1-k})$ commutes
with $F_k$.  Moreover, the subgroup $Aut(F_k)$ commutes
with $Aut(F_{p-1-k})$ as the relevant Nielsen transformations
and graph automorphisms of the two groups have nothing to
do with each other.
Finally, the subgroup
$$(F_k \rtimes Aut(F_k)) \times Aut(F_{p-1-k})$$
commutes with $N_{\Sigma_p}(\Z/p)$ for two reasons.
First, all of the Nielsen transformations in
$$(F_k \rtimes Aut(F_k)) \times Aut(F_{p-1-k})$$
commute with
$N_{\Sigma_p}(\Z/p)$ because if any one of these
Nielsen transformations moved some edge $e_j$
then it did this by pulling all of the edges
$e_i$ (which are what the subgroup $N_{\Sigma_p}(\Z/p)$
permutes) uniformly together around some path.
Second, it is clear that graph
automorphisms of the $\Theta$ graph in the
middle of $\Gamma$ commute with graph automorphisms of
either rose in $\Gamma$. \END

\begin{lemma} \label{tr11}
$$N_{Aut(F_n)}(D) = N_{\Sigma_p}(\Z/p) \times N_{\Sigma_p}(\Z/p).$$
$$N_{Aut(F_n)}(E) =
(N_{\Sigma_p}(\Z/p) \times N_{\Sigma_p}(\Z/p)) \rtimes \Z/2.$$
\end{lemma}

\PF The groups
$N_{Aut(F_n)}(D)$ and $N_{Aut(F_n)}(E)$ are the easiest to calculate of
all of the normalizers.
There are no admissible Nielsen transformations of the
reduced $(\Z/p \times \Z/p)$-graphs $\Omega_n$ or $\Psi_n$.
It follows that $N_{Aut(F_n)}(D)$ and $N_{Aut(F_n)}(E)$ are just
finite groups consisting of
normalizing graph automorphisms of $\Omega_n$ and
$\Psi_n$, respectively.  So we need only examine the
graphs $\Omega_n$ or $\Psi_n$ and see which graph
automorphisms are in the normalizers of the
respective $(\Z/p \times \Z/p)$-actions.

We will do this for $E \cong \Z/p \times \Z/p$
acting on $\Psi_n$. (The argument for $D$ acting on
$\Omega_n$ is similar.)  The automorphism group of the
graph $\Psi_n$ is
$$(\Sigma_p \times \Sigma_p) \rtimes \Z/2$$
where the first $\Sigma_p$ corresponds to permuting the
edges of first $\Theta$-graph in $\Psi_n$, the second
$\Sigma_p$ corresponds to permuting the edges of
the second $\Theta$-graph, and the $\Z/2$ corresponds to
flipping the graph over and thus exchanging the two $\Theta$-graphs.
Now $E \cong \Z/p \times \Z/p$ and the
$E$-action on the graph is given by the first
$\Z/p$ rotating the $p$ edges of the first $\Theta$-graph
while the latter $\Z/p$ rotates the $p$-edges of the
latter $\Theta$-graph.
Hence we need only compute the normalizer of
$$\Z/p \times \Z/p \subset
(\Sigma_p \times \Sigma_p) \rtimes \Z/2.$$
The subgroup
$$N_{\Sigma_p}(\Z/p) \times N_{\Sigma_p}(\Z/p)$$
is clearly in the normalizer.

Moreover, the flip $\Z/2$
is also in the normalizer because if
$$(a,b) \in \Z/p \times \Z/p \cong E$$
and $c$ is the flip in $\Z/2$ then
$$c \cdot (a,b) \cdot c^{-1} = (b,a) \in \Z/p \times \Z/p \cong E.$$
Further examination reveals that these are the only
normalizing graph automorphisms, so that
$$N_{Aut(F_n)}(E) =
(N_{\Sigma_p}(\Z/p) \times N_{\Sigma_p}(\Z/p)) \rtimes \Z/2. \square$$

\chapter{The normalizers, their fixed point sets, and related fixed point sets} \label{c8}

Let $G$ be a finite subgroup of $Aut(F_n)$, realized by a reduced
graph $\Gamma$ as in the previous chapter.  Define the
{\em fixed point subcomplex} $X_n^G$
of $G$ in the
spine $X_n$ of auter space by
$$X_n^G = \{x \in X_n : xg = x \hbox{ for all } g \in G\}.$$
The normalizer
$N_{Aut(F_n)}(G)$ acts on the fixed point set $X_n^G$ since if
$n \in N_{Aut(F_n)}(G)$ and $x \in X_n^G$ then for all $g \in G$,
$$(xn)g=(xngn^{-1})n=xn.$$

We will delay proving some of the more difficult results that we need in
this chapter until Part \ref{p3}.  For now, we will just
state the results that we need as facts.

\begin{fact} \label{tr12}
The space $X_n^G$ is a contractible, finite-dimensional
complex and \newline
$N_{Aut(F_n)}(G)$ acts on it with finite quotient and
finite stabilizers.
\end{fact}

Note that the poset structure on $X_n$ is such that a simplex $\delta$
in $X_n$ is fixed by $G$ if and only if it is fixed pointwise.
The space $X_n^G$ is finite dimensional since $X_n$ is
finite dimensional by \cite{[H]}.
That the normalizer acts on the space with finite quotient and
finite stabilizers follows from work by Krstic and Vogtmann in
\cite{[K-V]}.  Theorem \ref{tr13} in Part \ref{p3} proves that
$X_n^G$ is contractible.

A few spaces related to fixed point subcomplexes will come up so
frequently in our work that we will give them special names,
$X_p^\omega$, $Q_p^\omega$, $\tilde X_m$, and $\tilde Q_m$.
We now define these spaces.

Recall from Krstic and Vogtmann \cite{[K-V]} that an edge of a reduced
$G$-graph $\Gamma$ is {\em inessential} if it is contained in
every maximal $G$-invariant forest in $\Gamma$,
and that there is an $N_{Aut(F_p)}(G)$-equivariant deformation
retract of $X_n^G$ obtained by collapsing all inessential
edges of marked graphs in $X_n^G$.

\begin{defn}[$X_p^\omega$ and $Q_p^\omega$] \label{tr17}
Let $\omega$ be the automorphism of $F_p$ defined by
interchanging the first two basis elements of $F_p$
(cf. Remark \ref{tr9}),
and let $X_p^{\langle \omega \rangle}$ be the subcomplex of
$X_p$ fixed by the subgroup $\langle \omega \rangle \cong \Z/2$.
Define $X_p^\omega$ to be the associated
$N_{Aut(F_p)}(\omega)$-equivariant deformation retract
of $X_p^{\langle \omega \rangle}$ and let $Q_p^\omega$
be the quotient of $X_p^\omega$ by the action of
$N_{Aut(F_p)}(\omega)$.
\end{defn}

In Part \ref{p3} of this paper we will show (Remark \ref{tr29})
that
$$dim(X_p^\omega) = dim(Q_p^\omega) = 2p-4.$$

We will use $X_p^\omega$ to study the cohomology of the normalizer
of $A$.  In order to do this, we need to define an action of
$N_{Aut(F_n)}(A)$ on $X_p^\omega$.
Recall from Remark \ref{tr9} that
$$N_{\Sigma_p}(\Z/p) \times N_{Aut(F_p)}(\omega) \cong
\langle \omega \rangle \times N_{Aut(F_n)}(A).$$

The rest of this chapter will be devoted to
defining the space $\tilde X_m$, finding
an alternative characterization of this space,
and explaining how the space relates to
to the normalizer
$N_{Aut(F_n)}(B_m)$.

\begin{defn}[$\tilde X_m$ and $\tilde Q_m$] \label{tr18}
Let $m$ be a positive integer. Let $\alpha$ be the automorphism
of $F_{m+2}$ defined by:
$$\left\{ \matrix{
x_i \mapsto x_i \hbox{ for } i \leq m, \hfill \cr
x_{m+1} \mapsto x_{m+2} \hfill \cr
x_{m+2} \mapsto x_{m+2}^{-1} x_{m+1}^{-1}. \hfill \cr
} \right.$$
The subgroup $\mathcal{Q}$ generated by $\alpha$ has order 3,
and is realized on the graph $R_m \vee \Theta_2$ by rotating the
edges of $\Theta_2$ cyclically.  Define $\tilde X_m$ to be the fixed point set
$X_{m+2}^\mathcal{Q}$, and let $\tilde Q_m$ be the quotient of
$\tilde X_m$ by $N_{Aut(F_{m+2})}(\mathcal{Q})$.

For $m=0$, we set
$\tilde X_0 = \tilde Q_0 = \{\hbox{a point}\}$
as a notational convenience.
\end{defn}

Let $\tilde \Gamma$ be the graph $R_m \vee \Theta_{2}$,
where the basepoint of
the resulting graph $\tilde \Gamma$ is the center of the rose
$R_m$.
Label the $3$ edges of the $\Theta$-graph in $\tilde \Gamma$
as $e_1, e_2, e_3$ and orient them so that they begin
at the basepoint of $\tilde \Gamma$.
Construct a specific marked graph
$$\tilde \eta: R_{m+2} = R_{m} \vee R_{2} \to R_{m} \vee \Theta_{p-1}$$
by sending the first rose $R_{m}$ in $R_{m} \vee R_{2}$
to $R_{m}$ in $R_{m} \vee \Theta_{2}$ via the identity map.
Then send the second rose $R_{2}$ in $R_{m} \vee R_{2}$
to $\Theta_{2}$ in $R_{m} \vee \Theta_{2}$ by
sending the $i$th petal to $e_i * \bar e_{i+1}$.
As noted above in Definition \ref{tr18},
the subgroup $\mathcal{Q}$ acts on this marked graph by rotating the
edges of $\Theta_2$ cyclically.

We will
now study the structure of
$\tilde X_m$ and $\tilde Q_m$
in some detail.
Recall from Lemma \ref{tr10} that
$$N_{Aut(F_n)}(B_k) \cong N_{\Sigma_p}(\Z/p) \times (F_k \rtimes Aut(F_k))
\times Aut(F_{p-1-k}).$$

In a similar way, we could use
Proposition \ref{t4} just as Lemma \ref{tr10} did to
get that
$$N_{Aut(F_{m+2})}(\mathcal{Q}) \cong \Sigma_3
\times (F_m \rtimes Aut(F_m)).$$
The $\Sigma_3 \cong N_{\Sigma_3}(\Z/3)$ comes from normalizing
graph automorphisms of the $\Theta_{2}$, the
$F_{m}$ comes from Nielsen moves done by pulling
the back ends of all of
the edges $e_i$ uniformly around some loop in the rose
$R_{m}$, and the $Aut(F_{m})$ comes from Nielsen
transformations and graph automorphisms concerning just
the petals of the rose $R_{m}$.
  
Let us examine the Nielsen transformations corresponding
to the $F_{m}$ in a little more detail.
Say the petals of the rose are $r_1, \ldots, r_{m}$.
The group
$F_{m}$ is generated by the $m$ Nielsen moves
$a_k := <\bar e_1,\bar r_k>$ where $<\bar e_1,\bar r_k>$
fixes the rose and sends
$e_j$ to $r_k * e_j$ for all $j \in \{1, 2, 3\}$.
What effect does this move have on the marked
graph $\tilde \eta$?  First, observe that the new marked graph
sends $R_{m}$ identically to
$R_{m}$ as before.  However, the $i$th petal of the
second rose is now sent to $r_k * e_i * \bar e_{i+1}
* \bar r_k$.   Hence some word
$$w=a_{k_1}, \ldots, a_{k_s} \in
F_m \subset N_{Aut(F_{m+2})}(\mathcal{Q})$$
in the letters $a_1, \ldots, a_{m}$ acts on the
marked graph $\tilde \eta$ as follows.  The first rose $R_{m}$
in $R_{m} \vee R_{2}$ is still sent identically
to the rose $R_{m}$ in $R_{m} \vee \Theta_{2}$.
On the other hand, the second rose $R_{2}$ now
has its $i$th petal sent to (recalling
that the action of $N_{Aut(F_{m+2})}(\mathcal{Q})$ on marked
graphs is a right action)
$$r_{k_1} * \cdots * r_{k_s} * e_i * \bar e_{i+1}
* \bar r_{k_s} * \cdots * \bar r_{k_1}.$$
In other words, the entire second rose $R_{2}$
in $R_{m} \vee R_{2}$ is mapped to
$R_{m} \vee \Theta_{2}$ by conjugating the old
way it was mapped by the path
$r_{k_1} * \cdots * r_{k_s}$ representing $w$.

\begin{prop} \label{t11} The fixed point
space $\tilde X_{m}$ can be characterized
as the realization of the poset of equivalence classes of
pairs $(\alpha,f)$, where
$\alpha : R_{m} \to \Gamma_{m}$ is
a basepointed marked graph whose underlying graph $\Gamma_{m}$
has a special (possibly valence 2)
vertex which is designated as $\circ$,
$\circ$ may equal the basepoint $*$ of $\Gamma_{m}$,
and $f : I \to \Gamma_{m}$ is a homotopy
class (rel endpoints) of maps from
$*$ to $\circ$ in $\Gamma_{m}$. 
\end{prop}

\PF By definition,
$\tilde X_{m}$ consists
of simplices in the spine $X_{m+2}$ that are
fixed by $\mathcal{Q}$.  It is the subcomplex
generated by marked graphs (i.e., vertices of $X_{m+2}$) which
realize the finite subgroup $\mathcal{Q} \subset Aut(F_{m+2})$.
From Theorem \ref{tr4}, any
two vertices in $\tilde X_{m}$ corresponding
to reduced marked graphs, are connected to
each other by Nielsen moves.  In other words,
they are connected to each other by Nielsen moves that
can be represented by
elements of
$$N_{Aut(F_{m+2})}(\mathcal{Q}) \cong \Sigma_3
\times (F_m \rtimes Aut(F_m)),$$
since the normalizer contains all of the relevant Nielsen
transformations.
In particular, the Nielsen moves come from
$<a_1, \ldots, a_{m}> \cong F_{m}$ and the
Nielsen moves in
$Aut(F_{m})$ involving only the petals of the rose
$R_{m}$ in $R_{m} \vee \Theta_{2}$.

Hence the reduced marked graphs representing vertices
of $\tilde X_{m}$ are of the form
$$\psi = \alpha \vee \beta :
R_{m} \vee R_{2} \to R_{m} \vee \Theta_{2}$$
where $\alpha : R_{m} \to R_{m}$
corresponds to any reduced marked graph representing
a vertex of $X_{m}$ and
$\beta : R_{2} \to R_{m} \vee \Theta_{2}$
sends the  $i$th petal to 
$$r_{k,1} * \cdots * r_{k,m} * e_i * \bar e_{i+1}
* \bar r_{k,m} * \cdots * \bar r_{k,1},$$
for some path
$r_{k_1} * \cdots * r_{k_s}$ corresponding
to a word in $F_{m}$.
We could thus represent the reduced
marked graph $\psi$ more compactly as a pair
$(\alpha,f)$ where $\alpha: R_{m} \to R_{m}$
is any reduced marked graph and $f:I=[0,1] \to R_{m}$
is a homotopy class (rel the endpoints $f(0)=f(1)=*$) of maps
representing a path $r_{k_1} * \cdots * r_{k_s}$
in the rose $R_{m}$.

Now that we have a characterization of reduced
marked graphs representing vertices of $\tilde X_{m}$,
we will obtain a characterization of all
marked graphs representing vertices of
$\tilde X_{m}$.  We can do this by considering
the stars, in $\tilde X_{m}$, of reduced
marked graphs
$(\alpha,f): R_{m+2} \to R_{m} \vee \Theta_{2}$.
Recall that $\mathcal{Q}$ rotates the edges of $\Theta_{2}$
and leaves everything else fixed.

\begin{defn} \label{tr22} Let $G$ be a finite
subgroup of $Aut(F_r)$ for some integer $r$.
A marked graph
$$\eta^1 : R_{r} \to \Gamma^1$$
is a {\em $G$-equivariant blowup in the fixed point space
$X_r^G$}
of a marked graph
$$\eta^2 : R_{r} \to \Gamma^2$$
if there is a $1$-simplex $\eta^1 > \eta^2$ in $X_r^G$.
In other words, this happens exactly when we can
collapse a $G$-invariant forest in $\Gamma^1$
to obtain $\Gamma^2$. 
We often abbreviate this and
just say that $\eta_1$ is a {\em blowup} of $\eta_2$ or
that $\eta_2$ can be {\em blown up} to get $\eta_1$.
If the forest that we collapsed was just a tree, we
say that we {\em blew up} the corresponding vertex
of $\eta_2$ to get $\eta_1$.  
\end{defn}


\begin{claim} \label{tr23}
Any equivariant blowup of $\Gamma$ is a blowup of $R_m$,
with $\Theta_2$ attached at a point.
\end{claim}

\PF
Suppose $\Gamma$ is obtained from $\Gamma'$ by collapsing an invariant forest
$F$.  It suffices to show that $F$ is fixed by the action of $\mathcal{Q}$, since then the
initial (resp. terminal) vertices of the edges mapping to $e_i$ must be the same.

Suppose $F$ is not fixed by $\mathcal{Q}$, and let $F_0$
be the union of edges of $F$ with non-trivial orbits.  Let $v$ be 
a terminal vertex of $F_0$, and let $e$ and $f$ be edges of $\Gamma'$
terminating at $v$ which are not in $F_0$.  If $e$ or $f$ is fixed by 
$\mathcal{Q}$, then $v$ is fixed; in particular, if $e$ or $f$ is in
$F-F_0$, then $v$ is fixed by $\mathcal{Q}$.  If $e$ and $f$ are not 
fixed, then since they are not in $F$, they must map to edges $e_i$ and $e_j$ of 
$\Gamma$, with $i \not = j$.  Since the map is equivariant, some element of $\mathcal{Q}$
takes $e$ to $f$.  But this element must then fix their common vertex, $v$, so
in all cases $v$ must be fixed by $\mathcal{Q}$.  Since all terminal vertices of
$F_0$ are fixed, $F_0$ must be fixed, contradicting the definition of $F_0$. \END

Hence the only way to blow up $(\alpha,f)$ is
by blowing up the $R_{m}$ part of
$R_{m} \vee \Theta_{2}$.  When you do this,
you can think of the resulting marked graph as
some $(\hat \alpha, \hat f)$ where
$$\hat \alpha : R_{m} \to \Gamma_{m}$$ is any
marked graph, except that
the underlying graph $\Gamma_{m}$ has
one extra distinguished vertex, which we will call $\circ$
and which might have valence $2$,
aside from the basepoint $*$, 
and where
$$\hat f : I \to \Gamma_{m}$$
is a homotopy class of
paths from $*$ to $\circ$ in $\Gamma_{m}$.
We allow the possibility that $\circ$ might just be $*$.
The pair $(\hat \alpha, \hat f)$ is really
representing a marked
graph $R_{m+2} \to \Gamma_{m} \vee \Theta_{2}$
as follows.  The wedge $\vee$ connects the
point $\circ$ of $\Gamma_{m}$ to the left
hand vertex of $\Theta_{2}$.  The basepoint of
the resulting graph
$\Gamma_{m} \vee \Theta_{2}$
is whatever the old basepoint $*$ of $\Gamma_{m}$
was.  The first $m$ petals of
$R_{m+2}$ map to $\Gamma_{m}$ via $\hat \alpha$
and the $(m+i)$-th petal of $R_{m+2}$
first goes around the image of $\hat f$ from
$*$ to $\circ$, then goes around $e_i * \bar e_{i+1}$
of the $\Theta$-graph, and then finally
goes back around the image of $\hat f$ in
reverse from $\circ$ to $*$.

For $(\hat \alpha, \hat f)$ to be in the
star of $(\alpha,f)$ it has to be
the case that when the graph
$\Gamma_{m} \vee \Theta_{2}$
is collapsed to $R_{m} \vee \Theta_{2}$
the marking $\hat \alpha$ collapses
to $\alpha$ and the homotopy class
(rel endpoints) of maps $\hat f$ collapses
to $f$. This concludes the proof of
Proposition \ref{t11}. \END

As in \cite{[C-V]}, there is an obvious
definition of when two of the marked graphs described in
Proposition \ref{t11} 
are equivalent. 

\begin{defn} \label{tr24}
Marked graphs
$(\alpha^1,f^1)$ and $(\alpha^2,f^2)$
are {\em equivalent} 
if there
is a homeomorphism $h$ from 
$\Gamma_{m}^1$ to $\Gamma_{m}^2$
which sends $*$ to $*$, $\circ$ to $\circ$,
such that
$$(h \alpha^1)_\# = (\alpha^2)_\# : \pi_1(R_m,*) \to \pi_1(\Gamma_m^2)$$
and such that the paths
$$ hf_1, f_2 : I \to \Gamma_m^2$$
are homotopic rel endpoints.
\end{defn}

We are now in a position to observe 
how 
$N_{Aut(F_{m+2})}(\mathcal{Q})$ acts on a marked graph
in $\tilde X_{p-1}$.  Namely,

\begin{prop} \label{t12} The group
$$N_{Aut(F_{m+2})}(\mathcal{Q}) \cong
\Sigma_3 \times (F_{m} \rtimes Aut(F_{m}))$$
acts on a marked graph
$$\matrix{\hfill (\alpha,f): R_{m} \coprod I \to \Gamma_{m}\hfill \cr}$$
in $\tilde X_{m}$ as follows.
The subgroup $\Sigma_3$ 
permutes the edges of the $\Theta$ graph attached at $\circ$,
giving a marked graph which is equivalent to the
original one.  An element $a_{k_1} a_{k_2} \ldots a_{k_s}$
in $F_{m} = <a_1, \ldots, a_{m}>$ doesn't change
$\alpha$ at all, but sends the path $f$ to the
path 
$$\alpha(r_{k_1}) * \alpha(r_{k_2}) * \cdots * \alpha(r_{k_s}) * f$$
where $r_i$ is the $i$th petal of the rose $R_{m}$
in the domain of $\alpha$.
Lastly, an element $\phi \in Aut(F_{m})$
does not change $f$ at all and acts on
$\alpha: R_{m} \to \Gamma_{m}$
by precomposition:
$$(\alpha,f) \cdot \phi = (\alpha \circ \phi, f).$$
\end{prop}

\PF This follows in a direct manner from recalling that
we know how \newline
$N_{Aut(F_{m+2})}(\mathcal{Q})$ acts on reduced marked
graphs and just looking at $(\alpha,f)$ as being a
marked graph
$$R_{m+2} \to \Gamma_{m} \vee \Theta_{2}$$
in the star of some reduced marked graph
$$R_{m+2} \to R_{m} \vee \Theta_{2}. \hbox{ \END}$$

Just as Proposition \ref{t11} gives us
a simple characterization of $\tilde X_m$, the following remark
provides a nice characterization of $\tilde Q_m$:

\begin{remark} \label{tr25} The quotient
space $\tilde Q_{m}$ of $\tilde X_m$
by $N_{Aut(F_{m+2})}(\mathcal{Q})$
can be characterized
as the realization of the poset of equivalence classes of
basepointed graphs $\Gamma_{m}$
which have 
a special (possibly valence 2)
vertex designated as $\circ$,
which may equal the basepoint $*$.
Two such graphs $\Gamma_{m}^1$ and $\Gamma_{m}^2$
are {\em equivalent} if there 
is a homeomorphism
$$h: \Gamma_{m}^1 \to \Gamma_{m}^2$$
such that $h(*)=*$ and $h(\circ)=\circ$.
Define the poset structure
on these graphs by
forest collapses.  That is,
$\Gamma_{m}^1 > \Gamma_{m}^2$
if there is a simplicial map
$g: \Gamma_{m}^1 \to \Gamma_{m}^2$
such that $g(*)=*$, $g(\circ)=\circ$,
and $g^{-1}(\hbox{vertices of }\Gamma_{m}^2)$ is
a subforest of $\Gamma_{m}^1$.
\end{remark}

Note that if we are concerned with the quotient
$\tilde Q_{m}$ of $\tilde X_{m}$
by $N_{Aut(F_{m+2})}(\mathcal{Q})$ then the situation
simplifies even further
from that of Proposition \ref{t11}.  When you
forget the marking of $(\alpha,f)$
on $\Gamma_{m} \vee \Theta_{m}$,
it becomes irrelevant where
$\alpha$ sends the rose $R_{m}$
and where
$f$ sends the interval $I$.
In other words, all that matters is the
graph $\Gamma_{m}$ with two distinguished
points, $*$ and $\circ$.  In this manner, the
quotient $\tilde Q_{m}$ is just the ``moduli space
of unmarked graphs $\Gamma_{m}$ with
$\pi_1(\Gamma_{m}) \cong F_{m}$
and where $\Gamma_{m}$ has two distinguished points.''

Recall from \cite{[H-V]} that the spine $X_m$ of auter space is a
deformation retraction of auter space $\A_m$.  Similarly, we can
construct a space $\tilde \A_{m}$ which deformation retracts
to $\tilde X_m$.  We can then think of $\tilde X_m$ as
being the ``spine'' of $\tilde \A_m$.

\begin{defn}[$\tilde \A_{m}$] \label{tr26}
Construct an
analog $\tilde \A_{m}$ of auter space
for $N_{Aut(F_{m+2})}(\mathcal{Q})$ by considering
markings
$$\matrix{\hfill (\alpha,f): R_{m} \coprod I \to \Gamma_{m}\hfill \cr}$$
where
the edges of $\Gamma_{m}$ are assigned
lengths which must sum to 1.  Just as
in \cite{[C-V]}, the space $\tilde \A_{m}$
deformation retracts to its spine
$\tilde X_{m}$.
\end{defn}

Now from $\cite{[H-V]}$ we have
$dim(X_{m}) = dim(Q_{m}) = 2m-2$
and $dim(\A_{m}) = 3m-3$.  As the graph
$\Gamma_{m}$ in a particular marked graph
has possibly one extra vertex
$\circ$ of valence $2$, we see that
$dim(\tilde X_{m}) = dim(\tilde Q_{m}) = 2m-1$
and $dim(\tilde \A_{m}) = 3m-2$.

Now that we are more familiar with the
structure of our spaces
$\tilde X_m$ and $\tilde Q_m$,
we can proceed with describing how they are
related to the normalizers
$N_{Aut(F_n)}(B_k)$.

\begin{defn}[Action of $N_{Aut(F_n)}(B_k)$ on $\tilde X_k$] \label{tr20}
For $k \in \{0, \ldots, p-1\}$,
the map from
$$N_{Aut(F_n)}(B_k) \cong N_{\Sigma_p}(\Z/p) \times (F_k \rtimes Aut(F_k))
\times Aut(F_{p-1-k})$$
to
$$N_{Aut(F_{m+2})}(\mathcal{Q}) \cong \Sigma_3
\times (F_m \rtimes Aut(F_m))$$
given by projection to the second factor followed by
inclusion induces an action of 
$N_{Aut(F_n)}(B_k)$ on $\tilde X_k$.
Since $N_{\Sigma_p}(\Z/p) \subset N_{Aut(F_n)}(B_k)$ acts
trivially, the quotient of this action is $\tilde Q_k$. 
\end{defn}

We can also define an action of $N_{Aut(F_n)}(B_k)$
on $X_{p-1-k}$:

\begin{defn}[Action of $N_{Aut(F_n)}(B_k)$ on $X_{p-1-k}$] \label{tr21}
For $k \in \{0, \ldots, p-1\}$,
the map from
$$N_{Aut(F_n)}(B_k) \cong N_{\Sigma_p}(\Z/p) \times (F_k \rtimes Aut(F_k))
\times Aut(F_{p-1-k})$$
to
$$Aut(F_{p-1-k})$$
given by projection on the third factor induces an action of
$N_{Aut(F_n)}(B_k)$ on $X_{p-1-k}$,
with quotient $Q_{p-1-k}$.
\end{defn}

From Definition \ref{tr20} and Definition \ref{tr21},
there is an induced action of \newline
$N_{Aut(F_n)}(B_k)$ on
$\tilde X_k \times X_{p-1-k}$.
Note that the space
$\tilde X_k \times X_{p-1-k}$
is a product of posets and can be given a poset
structure by saying that
the pair of marked graphs
$$(\alpha^1,f^1) \times \beta^1$$
is greater than or equal to the pair
$$(\alpha^2,f^2) \times \beta^2$$
if both
$$(\alpha^1,f^1) \geq (\alpha^2,f^2)$$
and
$$\beta^1 \geq \beta^2.$$  This gives
$\tilde X_k \times X_{p-1-k}$
a simplicial structure.  Often we will
not use this simplicial structure, however,
and just use the cellular structure that comes
from products of simplices in
$\tilde X_k$ and in $X_{p-1-k}$.

\begin{thm} \label{tr28}
The above induced action of
$N_{Aut(F_n)}(B_k)$ on the contractible
space $\tilde X_k \times X_{p-1-k}$
has finite stabilizers and
quotient
$\tilde Q_k \times Q_{p-1-k}$.
\end{thm}

\PF The space $\tilde X_k$ 
is contractible from Fact \ref{tr12}.  The space
$X_{p-1-k}$ is contractible by
Hatcher and Vogtmann in \cite{[H-V]}.
It follows from the assertions in
Definition \ref{tr20} and Definition \ref{tr21}
that the quotient of
$\tilde X_k \times X_{p-1-k}$
by
$N_{Aut(F_n)}(B_k)$
is
$\tilde Q_k \times Q_{p-1-k}$.

Consider the product of simplices
$$\delta_1 \times \delta_2 \subset \tilde X_k \times X_{p-1-k}$$ 
where $\delta_1$ and $\delta_2$
are simplices of
$\tilde X_k$ and $X_{p-1-k}$, respectively.
Let
$$\matrix{\hfill (\alpha,f): R_{k} \coprod I \to \Gamma_{k}^1\hfill \cr}$$
be the top marked graph in the simplex $\delta_1$
and let
$$\beta: R_{p-1-k} \to \Gamma_{p-1-k}^2$$
be the top marked graph in the simplex $\delta_2$.
Then the stabilizer of
$\delta_1 \times \delta_2$
under the action of
$$N_{Aut(F_n)}(B_k) \cong N_{\Sigma_p}(\Z/p) \times (F_k \rtimes Aut(F_k))
\times Aut(F_{p-1-k})$$
is a subgroup of
$$N_{\Sigma_p}(\Z/p) \times Aut(\Gamma_k^{1}) \times Aut(\Gamma_{p-1-k}^2),$$
which is a finite group.
Hence the action of
$N_{Aut(F_n)}(B_k)$ on 
$\tilde X_k \times X_{p-1-k}$
has finite stabilizers. \END

It is interesting to note that the space
$\tilde X_k \times X_{p-1-k}$ is in fact
homeomorphic to the fixed point subcomplex
$X_n^{B_k}$.  We leave this as a
straightforward exercise for the
reader, as this fact will not be used in this dissertation.

\chapter{The cohomology of the normalizers} \label{c4}

In this chapter, we will compute the cohomology of the normalizers
of the subgroups
$$A, B_k, C, D, \hbox{ and } E$$
listed in
Proposition \ref{tr5}.
First, we list some helpful facts that will allow us to
compute these cohomology groups.

Associated to any group extension is a ``five term
exact sequence in cohomology'' (see
section $3.5$ of \cite{[B2]} for a proof and details):

\begin{prop} \label{tr14}
Let $G$ be a group, $N$ be a normal subgroup of $G$, and $M$
be a $G$-module.  Then the following sequence is exact:
$$\matrix{
0 & \to & H^1(G/N; M^N) & \buildrel inf \over \to
& H^1(G; M)     & \buildrel res \over \to & H^1(N; M)^{G/N} 
& \buildrel d_2 \over \to  &  \cr
&     &               &     &   &   &  H^2(G/N; M^N)
& \buildrel inf \over \to & H^2(G; M), \cr
}$$
where
{\em inf} refers to the inflation map induced by the
quotient map $G \surjarrow G/N$, {\em res} refers to the restriction map
induced by the inclusion $N \injarrow G$, and $d_2$ is the differential
on the $E_2$-page of a spectral sequence defined in section
$3.5$ of \cite{[B2]}.
\end{prop}

One also has, from section $3.5$ of \cite{[B2]},
a ``five term exact sequence in homology'':

\begin{prop} \label{tr34}
Let $G$
be a group, $N$ be a normal subgroup of $G$, and $M$
be a $G$-module.  Then the following sequence is exact:
$$\matrix{
H_2(G; M)       & \to & H_2(G/N; M_N) & \to &  &  &  & &  \cr
 & & H_1(N; M)_{G/N} & \to & H_1(G; M)     & \to & H_1(G/N; M_N) & \to & 0.  \cr
}$$
\end{prop}

We will also need a result of Hatcher 
in \cite{[H]}:

\begin{thm}[Hatcher] \label{tr15}
The homology of the infinite
symmetric group \newline
$H_t(\Sigma_\infty; \Z)$ is a direct
summand of $H_t(Aut(F_\infty); \Z)$.  For $t=1,2$
the complementary summand is zero.
In addition,
$H_1(Aut(F_\infty); \Z) = H_1(Aut(F_m); \Z)$
for $m \geq 2$.
\end{thm}

It will also prove useful to cite a theorem of Hatcher and Vogtmann
in \cite{[H-V]} which gives a linear stability range for
the homology of $Aut(F_n)$:

\begin{thm}[Hatcher-Vogtmann] \label{tr16}
The natural stabilization map
$$H_t(Aut(F_m); \Z) \to H_t(Aut(F_{m+1}); \Z)$$
is an isomorphism for $n \geq 2t+3$.
\end{thm}

The following theorem of Swan's
(see \cite{[S]} or \cite{[A-M]}) is a standard
tool for computing the cohomology of groups:
\begin{thm}[Swan] \label{tr30}
If $G$ is a finite group with a $p$-Sylow
subgroup $P$ that is abelian, then
$$H^*(G; \Z_{(p)}) = H^*(P; \Z_{(p)})^{N_{G}(P)}.$$
\end{thm}

We begin our computations with
a few preliminary lemmas
about the first cohomology groups of the
quotient spaces $Q_m$ and $\tilde Q_m$.

\begin{lemma} \label{tr32}
\begin{enumerate}
\item For all integers $1 \leq m \leq p-1$ and for $t=1,2$,
$$H^t(Q_m; \Z/p) = H^t(Aut(F_m); \Z/p).$$
\item For all integers $1 \leq m \leq p-1$,
$$H^1(Aut(F_m); \Z/p)=0.$$
\end{enumerate}
\end{lemma}

\PF First, we show that
$H^t(Q_m; \Z/p) = H^t(Aut(F_m); \Z/p)$
for $t=1,2$.
Apply the equivariant cohomology spectral
sequence \ref{e1} corresponding to
$Aut(F_m)$ acting on $X_m$.  This requires us
to look at the cohomology groups
$$H^*(stab(\delta); \Z/p)$$
of stabilizers of simplices $\delta$
in $X_m$.

If $m < p-1$, none of the
simplices $\delta$ have any $p$-symmetry.
So $H^t(stab(\delta); \Z/p)$ is $0$ for
all $t \not = 0$ in those cases.
For $t=0$, the cohomology group is just
$\Z/p$.  Hence the spectral sequence \ref{e1}
applied to this situation
yields that
$H^t(Q_m; \Z/p) = H^t(Aut(F_m); \Z/p)$
for all $t$.  Specializing this to
$t=1,2$ gives us our desired result.

If $m=p-1$, then the only simplices
$\delta$ with $p$-symmetry
are vertices (i.e., marked graphs)
whose underlying graph is
$\Theta_{p-1}$.  For a such a $\delta$,
$$H^t(stab(\delta); \Z/p) = H^t(\Sigma_p; \Z/p) =
\left\{\matrix{
\Z/p \hfill &t \geq 0, t \equiv 0 \hbox{ } (\hbox{mod } n) \hfill \cr
\Z/p \hfill &t > 0, t \equiv -1 \hbox{ } (\hbox{mod } n) \hfill \cr
0 \hfill &\hbox{otherwise } \hfill \cr} \right.$$
Since $p \geq 3$, the three rows $s=1,2,3$ of
the $E_1$-page of spectral sequence \ref{e1}
are all zero. Consequently, the spectral sequence gives us that
$H^t(Q_{p-1}; \Z/p) = H^t(Aut(F_{p-1}); \Z/p)$
for $t=1,2$.

Hence if $t=1,2$ then
$H^t(Q_m; \Z/p) = H^t(Aut(F_m); \Z/p)$
in both cases and $1.$ is established.

It remains to show that
$H^1(Aut(F_m); \Z/p) = 0$.
If $m=1$ then 
$$H^1(Aut(F_{m}); \Z/p) = H^1(Aut(\Z); \Z/p) = H^1(\Z/2; \Z/p) = 0$$
and we are done.
On the other hand, if $m \geq 2$ then Theorem \ref{tr15}
and dual universal coefficients yield
$$H^1(Aut(F_{m}); \Z/p) = H^1(Aut(F_{\infty}); \Z/p) =
H^1(\Sigma_\infty; \Z/p) = 0. \square$$

If $m \geq 5$ then Theorem \ref{tr15}
and Theorem \ref{tr16}, along with dual
universal coefficients, yield
$$H^2(Aut(F_{m}); \Z/p) = H^2(Aut(F_{\infty}); \Z/p) =
H^2(\Sigma_\infty; \Z/p) = 0.$$

If $p \geq 7$, this gives
$$H^2(Aut(F_{p-1}); \Z/p)=0.$$
We will also need this group to vanish for $p=3$ and $p=5$,
in the course of our proof of Lemma \ref{t6}.
Although the theorems of Hatcher and Vogtmann do not as
stated yield the result for $p=3,5$, it
is an easy application of their work to obtain this:

\begin{fact} \label{tr35}
$H^2(Aut(F_{p-1}); \Z/p)=0.$
\end{fact}

For a proof of the above fact, refer to Proposition \ref{tr36}
in Chapter \ref{c10}. 
 
\begin{lemma} \label{tr33}
For all integers $1 \leq m \leq p-1$,
\begin{enumerate} 
\item $|H_1(F_{m} \rtimes Aut(F_{m}); \Z)| \in \Z^{+}$ is a
power of $2$.
\item $H^1(\tilde Q_m; \Z/p) =
H^1(F_{m} \rtimes Aut(F_{m}); \Z/p) = 0.$
\end{enumerate}
\end{lemma}

\PF For just this proof, set
$$G = F_{m} \rtimes Aut(F_{m})$$
as a notational convenience.
First, we show 1.
From the five term exact sequence \ref{tr34}
corresponding to the group
extension $G$, the following is exact:
$$H_1(F_{m}; \Z)_{Aut(F_{m})} \to H_1(G; \Z)
\to H_1(Aut(F_{m}); \Z_{F_{m}}).$$

Now $G$ acts trivially on $\Z$ so that $\Z_{F_{m}}=\Z$
and from Theorem \ref{tr15} we see that
$$H_1(Aut(F_{m}); \Z) = H_1(Aut(F_\infty); \Z) =
H_1(\Sigma_\infty; \Z) = \Z/2.$$

In addition,
$$H_1(F_{m}; \Z) = \Z \oplus \ldots \oplus \Z = m(\Z)$$
and $Aut(F_{m})$ acts on it by first projecting to
$GL_{m}(\Z)$ and then acting in the usual manner.
Let $\xi$ be the element of $Aut(F_m)$ sending each generator to 
its inverse.  Then $\xi$ acts on $m(\Z) = \Z \oplus \ldots \oplus \Z$
by sending $(x_1, \ldots, x_{m})$ to $(-x_1, \ldots, -x_{m})$.
Thus in the module of co-invariants
$(\Z \oplus \ldots \oplus \Z)_{Aut(F_{m})}$
we have
$$\left[(x_1, \ldots, x_{m})\right]=\left[(-x_1, \ldots, -x_{m})\right]$$
or $2\left[(x_1, \ldots, x_{m})\right]=0.$
Accordingly, we have that
$(\Z \oplus \ldots \oplus \Z)_{Aut(F_{m})}$
is all $2$-torsion.

From the three terms of the five term exact sequence
that we printed above,
$H_1(G; \Z)$ is a finite group whose order is a power of $2$.
This proves statement 1.

To prove the second statement, we use the equivariant
cohomology spectral sequence \ref{e1} for
$$G \subset N_{Aut(F_{m+2})}(\mathcal{Q})$$
acting on $\tilde Q_m$.
If $m < p-1$, none of the
simplices $\delta$ of $\tilde Q_m$ have any $p$-symmetry.
So $H^t(stab(\delta); \Z/p)$ is $0$ for
all $t \not = 0$ in those cases.
For $t=0$, the cohomology group is just
$\Z/p$.  Hence the spectral sequence \ref{e1}
applied to this situation
yields that
$H^t(\tilde Q_m; \Z/p) = H^t(G; \Z/p)$
for all $t$.  Specializing this to
$t=1$ gives us that
$H^1(\tilde Q_m; \Z/p) = H^1(G; \Z/p)$.

On the other hand, if $m=p-1$, then the only simplices
$\delta$ with $p$-symmetry
are vertices (i.e., marked graphs) whose
underlying graph is
$\Theta_{p-1}$.  (In other words, there are
really two possibilities for underlying graphs here.
One is the graph $\Theta_{p-1}^1$
where $*$=$\circ$ are both the same vertex on one side of
the $\Theta$-graph, and the other is
the graph $\Theta_{p-1}^2$
where
$*$ is the vertex on one side of the $\Theta$-graph
and $\circ$ is the vertex on the other side.)
For a such a $\delta$,
$$H^t(stab(\delta); \Z/p) = H^t(\Sigma_p; \Z/p) =
\left\{\matrix{
\Z/p \hfill &t \geq 0, t \equiv 0 \hbox{ } (\hbox{mod } n) \hfill \cr
\Z/p \hfill &t > 0, t \equiv -1 \hbox{ } (\hbox{mod } n) \hfill \cr
0 \hfill &\hbox{otherwise } \hfill \cr} \right.$$
Since $p \geq 3$, the three rows $s=1,2,3$ of
the $E_1$-page of spectral sequence \ref{e1}
are all zero. Consequently, the spectral sequence gives us that
$H^1(\tilde Q_m; \Z/p) = H^1(G; \Z/p)$.

Hence $H^1(\tilde Q_m; \Z/p) = H^1(G; \Z/p)$
in both cases.  By part $1.$ and dual
universal coefficients, we see that
$$H^1(\tilde Q_m; \Z/p) = H^1(G; \Z/p) = 0. \square$$

Finally, we will need the following fact,
which is proven later in Proposition \ref{t15}
and Proposition \ref{t16} of
Part \ref{p2}:

\begin{fact} \label{tr31}
$H^2(\tilde Q_{p-1}; \Z/p) = 0$.
\end{fact}

We can now calculate the cohomology of the
normalizers of the various subgroups $A$, $B_k$, $C$,
$D$, and $E$.

\begin{lemma} \label{t5} $$\hat H^t(N_{Aut(F_n)}(A); \Z_{(p)})
\cong \left\{\matrix{
\Z/p \hfill &t \equiv 0  \hbox{ } (\hbox{mod } n) \hfill \cr
0 \hfill &t \equiv \pm 1  \hbox{ } (\hbox{mod } n) \hfill \cr
H^r(Q_{p}^\omega; \Z/p) \hfill &t \equiv r  \hbox{ } (\hbox{mod } n), \hfill \cr
 &2 \leq r \leq n-2 \hfill \cr
} \right.$$
\end{lemma}

\PF 
We now define an action of $N_{\Sigma_p}(\Z/p) \times N_{Aut(F_p)}(\omega)$
on the space $X_{p}^\omega$ by
stipulating that $N_{\Sigma_p}(\Z/p)$ acts trivially and
that $N_{Aut(F_p)}(\omega)$ acts in the usual manner on $X_{p}^\omega$.
This in turn defines an action of $N_{Aut(F_n)}(A)$
on the contractible space $X_p^\omega$.  This action has finite
stabilizers and quotient $Q_p^\omega$.

The equivariant cohomology spectral sequence \ref{e2}
for this action is
\begin{equation} \label{e5}
E_1^{r,s} = \prod_{[\delta] \in \Delta_{p-1}^r}
\hat H^s(stab_{N_{Aut(F_n)}(A)}(\delta); \Z_{(p)}) \Rightarrow
\hat H^{r+s}(N_{Aut(F_n)}(A)); \Z_{(p)})
\end{equation}
where $[\delta]$ ranges over the set $\Delta_{p-1}^r$ of orbits
of $r$-simplices $\delta$ in $X_{p}^\omega$.

We claim
that each $stab_{N_{Aut(F_n)}(A)}(\delta)$ is the direct sum of
$N_{\Sigma_p}(\Z/p)$ with a finite subgroup of
$N_{Aut(F_p)}(\omega)$ that does
not have any $p$-torsion.  This is because
Glover and Mislin showed in \cite{[G-M]} that
the only $p$-torsion in
$Aut(F_{p})$ comes from stabilizers of marked graphs with underlying
graphs $\Theta_{p}$, $\Theta_{p-1} \vee R_1$,
$\Xi_p$, or $R_p$. (The graph $\Xi_p$, defined to be the
$1$-skeleton of the
cone over a $p$-gon, is illustrated in
Figure \ref{fig2}.)  Since $\omega$ acts by switching the
first two petals of the rose $R_p$,
$X_{p}^\omega$ does not contain any marked
graphs with underlying graph
$\Theta_{p}$, $\Theta_{p-1} \vee R_1$, or $\Xi_p$.
Although $X_{p}^\omega$ obviously does contain marked
graphs with underlying graph $R_p$, these will not worry us
as $p$ does not divide the orders of the stabilizers
-- under the action of
just $N_{Aut(F_p)}(\omega)$ -- of such marked
graphs.  For example,
$N_{Aut(F_p)}(\omega)$
does not contain
the permutation $(1 2 \ldots p)$ that
rotates the petals of the rose $R_p$
because
$$(1 2 \ldots p) \circ \omega \circ (p  \ldots 2 1)$$
is the permutation
$(2 3)$, which is not equal to $(1)$ or
$\omega = (1 2)$.

Thus for every $[\delta]$, we have
$$\hat H^s(stab_{N_{Aut(F_n)}(A)}(\delta); \Z_{(p)}) =
\hat H^s(N_{\Sigma_p}(\Z/p); \Z_{(p)})
= \hat H^s(\Sigma_p; \Z_{(p)}).$$
The $E_1^{r,s}$-page of the spectral sequence is $0$ in the rows
where $s \not = kn$.  On the rows where $s = kn$, each $[\delta]$ of
dimension $r$ contributes exactly one copy of $\Z/p$ to $E_1^{r,s}$.
Recall that the quotient of $X_{p}^\omega$ by $N_{Aut(F_{p})}(\omega)$
is $Q_{p}^\omega$.  Now each row where $s=kn$ is simply a copy
of the cellular cochain complex with $\Z/p$-coefficients of the
$(n-2)$-dimensional complex $Q_{p}^\omega$:
$$E_2^{r,s} = \left\{\matrix{
\Z/p \hfill &r=0 \hbox{ and } s=kn \hfill \cr
H^r(Q_{p}^\omega; \Z/p) \hfill &1 \leq r \leq n-2
\hbox{ and } s=kn \hfill \cr
0 \hfill &\hbox{otherwise } \hfill \cr} \right.$$
Hence we see that the spectral sequence converges at the $E_2$-page.

It remains to show that $H^1(Q_{p}^\omega; \Z/p) =0.$
From the equivariant spectral sequence \ref{e1} of the
group
$$\bigl( (F_{p-2} \times F_{p-2})
\rtimes (\Z/2 \times Aut(F_{p-2}))  \bigr)
\subset N_{Aut(F_n)}(A)$$
acting on the space $X_p^\omega$, we have that
$$H^1(Q_{p}^\omega; \Z/p) =  H^1\bigl( (F_{p-2} \times F_{p-2})
\rtimes (\Z/2 \times Aut(F_{p-2})); \Z/p  \bigr).$$

From the five term sequence of a group extension \ref{tr14},
we have that the sequence
$$\matrix{
\hfill H^1(\Z/2 \times Aut(F_{p-2}); \Z/p) \to \hfill \cr
\hfill H^1\bigl( (F_{p-2} \times F_{p-2})
\rtimes (\Z/2 \times Aut(F_{p-2})); \Z/p  \bigr) \to \hfill \cr
\hfill H^1(F_{p-2} \times F_{p-2};\Z/p)^{(\Z/2 \times Aut(F_{p-2}))} \hfill \cr
}$$
is exact.

First off,
$$H^1(\Z/2 \times Aut(F_{p-2}); \Z/p)
= H^1(Aut(F_{p-2}); \Z/p) = 0$$
from Lemma \ref{tr32}.

Next, we must show that
$$H^1(F_{p-2} \times F_{p-2};\Z/p)^{(\Z/2 \times Aut(F_{p-2}))}=0.$$
Basic calculations in the cohomology of groups yield that
$$H^1(F_{p-2} \times F_{p-2};\Z/p)^{(\Z/2 \times Aut(F_{p-2}))}
=
\bigl( (p-2)\Z/p \oplus (p-2)\Z/p \bigr)^{(\Z/2 \times
Aut(F_{p-2}))},$$
where the $\Z/2$ acts by exchanging the two copies of
$(p-2)\Z/p$ 
and the $Aut(F_{p-2})$ acts on the
$(p-2)\Z/p \oplus (p-2)\Z/p$ by first
projecting to $GL_{p-2}(\Z/p)$ and then
acting diagonally.

Suppose that
$$(x_1, \ldots, x_{p-2}) \oplus (y_1, \ldots, y_{p-2})$$
is an element of
$$\bigl( (p-2)\Z/p \oplus (p-2)\Z/p \bigr)^{(\Z/2 \times
Aut(F_{p-2}))}.$$
Define an element $\xi$ of $Aut(F_{p-2})$
on the generators $a_i$ of $F_{p-2}$
by $\xi(a_i)=a_i^{-1}$
for all $i$.
The automorphism $\xi$ sends
$$(x_1, \ldots, x_{p-2}) \oplus (y_1, \ldots, y_{p-2})$$
to
$$(-x_1, \ldots, -x_{p-2}) \oplus (-y_1, \ldots, -y_{p-2}).$$
So $x_i = -x_i$ and $y_j = -y_j$ for all $i,j$.
Recall that each $x_i$ and $y_j$ is in the additive group $\Z/p$.
As $p$ is odd, the only possibility is that $x_i=0$
for all $i$ and $y_j=0$ for all $j$.
Consequently, we have that
$$\bigl( (p-2)\Z/p \oplus (p-2)\Z/p \bigr)^{(\Z/2 \times
Aut(F_{p-2}))}=0.$$

The five term exact sequence thus gives that
$$H^1(Q_{p}^\omega; \Z/p) =  H^1\bigl( (F_{p-2} \times F_{p-2})
\rtimes (\Z/2 \times Aut(F_{p-2})); \Z/p  \bigr)=0.$$
This completes the proof of the lemma. \END

Note that the proof of this lemma was
complicated by the fact that there is no readily defined
K\"unneth formula for Farrell
cohomology.  Despite the fact that
$$\Z/2 \times N_{Aut(F_n)}(A)
\cong N_{\Sigma_p}(\Z/p) \times N_{Aut(F_p)}(\omega)$$
and that
a spectral sequence argument easily yields that
$\hat H^t(N_{Aut(F_p)}(\omega); \Z_{(p)})=0$,
we could not from then go on to
claim that
$$\hat H^t(N_{Aut(F_n)}(A); \Z_{(p)})=\hat H^t(\Sigma_p; \Z_{(p)}).$$

\begin{lemma} \label{t6} 
$$\matrix{
\hfill \cr
\hat H^t(N_{Aut(F_n)}(B_0); \Z_{(p)}) \cong
\left\{\matrix{
\Z/p^2 \hfill &t = 0 \hfill \cr
(j+1) \Z/p \hfill &|t| = nj \not = 0 \hfill \cr
(j-1) \Z/p \hfill &|t| = nj-1 \hfill \cr
0 \hfill &|t| \equiv 1  \hbox{ } (\hbox{mod } n) \hfill \cr
0 \hfill &t \equiv 2  \hbox{ } (\hbox{mod } n) \hfill \cr
H^r(Q_{p-1};\Z/p) \hfill
&t \equiv r  \hbox{ } (\hbox{mod } n), \hfill \cr
 \hfill &3 \leq r \leq n-2 \hfill \cr
}\right.
\hfill \cr
\hfill & \hfill \cr
}$$
$$\matrix{
\hbox{For $k \in \{1,\ldots,p-2\},$} \hfill \cr
\hfill \cr
\hat H^t(N_{Aut(F_n)}(B_k); \Z_{(p)}) \cong
\left\{\matrix{
\Z/p \hfill &t \equiv 0  \hbox{ } (\hbox{mod } n) \hfill \cr
0 \hfill &t \equiv \pm 1, -2,  \hbox{ } (\hbox{mod } n) \hfill \cr
H^r(\tilde Q_{k} \times Q_{p-1-k};\Z/p) \hfill
&t \equiv r  \hbox{ } (\hbox{mod } n), \hfill \cr
\hfill &2 \leq r \leq n-3 \hfill \cr
} \right.
\hfill \cr
\hfill & \hfill \cr}$$
$$\matrix{
\hfill \cr
\hat H^t(N_{Aut(F_n)}(B_{p-1}); \Z_{(p)}) \cong
\left\{\matrix{
\Z/p^2 \oplus \Z/p \hfill &t = 0 \hfill \cr
(2j+1) \Z/p \hfill &|t| = nj \not = 0 \hfill \cr
0 \hfill &t = nj + 1 > 0 \hfill \cr
(2j-2) \Z/p \hfill &t = -nj+1 < 0 \hfill \cr
(2j-2) \Z/p \oplus H^{n-1}(\tilde Q_{p-1}; \Z/p) \hfill
&t = nj-1 > 0 \hfill \cr
H^{n-1}(\tilde Q_{p-1}; \Z/p) \hfill &t = -nj-1 < 0 \hfill \cr
0 \hfill &t \equiv 2  \hbox{ } (\hbox{mod } n) \hfill \cr
H^r(\tilde Q_{p-1}; \Z/p) \hfill &t \equiv r  \hbox{ } (\hbox{mod } n), \hfill \cr
 &3 \leq r \leq n-2 \hfill \cr
} \right.
\hfill \cr
\hfill \cr}$$
\end{lemma}


\PF We use the equivariant cohomology spectral
sequence \ref{e2} corresponding to the action,
defined in Definition \ref{tr20} and Definition \ref{tr21},
of
$$N_{Aut(F_n)}(B_k)
\cong N_{\Sigma_p}(\Z/p) \times (F_k \rtimes Aut(F_k))
\times Aut(F_{p-1-k})$$
on the space $\tilde X_k \times X_{p-1-k}$.
From Theorem \ref{tr28}
the action of
$N_{Aut(F_n)}(B_k)$ on the contractible
space $\tilde X_k \times X_{p-1-k}$
has finite stabilizers and
quotient
$\tilde Q_k \times Q_{p-1-k}$.

Applying the spectral sequence \ref{e2} we obtain
the following $E_1$-page:
\begin{equation} \label{e6}
E_1^{r,s} = \prod_{[\delta] \in \Delta^r}
\hat H^s(stab_{N_{Aut(F_n)}(B_k)}(\delta); \Z_{(p)}) \Rightarrow
\hat H^{r+s}(N_{Aut(F_n)}(B_k)); \Z_{(p)})
\end{equation}
where $[\delta]$ ranges over the set $\Delta^r$ of orbits
of $r$-simplices $\delta$ in $\tilde X_k \times X_{p-1-k}$.
Vertices in $\Delta^0$ are
pairs of unmarked graphs 
$(\Gamma_1,\Gamma_2)$ where $\pi_1(\Gamma_1) \cong F_k$,
$\pi_1(\Gamma_2) \cong F_{p-1-k}$, the graph $\Gamma_1$ has
two distinguished points $*$ and $\circ$,
and the graph $\Gamma_2$ has one distinguished point $*$.

The stabilizer of this vertex of $\Delta^0$
under the action of
$$N_{Aut(F_n)}(B_k)
\cong N_{\Sigma_p}(\Z/p) \times (F_k \rtimes Aut(F_k))
\times Aut(F_{p-1-k})$$
is isomorphic to
$$N_{\Sigma_p}(\Z/p) \times Aut(\Gamma_1) \times Aut(\Gamma_2)$$
where by $Aut(\Gamma_1)$ we mean graph automorphisms of
$\Gamma_1$ that preserve
both distinguished points, and $Aut(\Gamma_2)$ is the group of graph
automorphisms of $\Gamma_2$ that preserve its distinguished point.

There are three cases, depending upon what value
$k$ takes.
\begin{enumerate}
\item $k \in \{1,\ldots,p-2\}.$
Consider an $r$-simplex $\delta$
$$((\alpha^0,f^0),\beta^0) > ((\alpha^1,f^1),\beta^1) >
\cdots > ((\alpha^r,f^r),\beta^r)$$
of $\tilde X_k \times X_{p-1-k}$.
Let $\Gamma^1$ be the underlying graph of
$(\alpha^0,f^0)$ and let $\Gamma^2$ be the underlying
graph of $\beta^0$.
The stabilizer of $\delta$ is a subgroup of the group
$$N_{\Sigma_p}(\Z/p) \times Aut(\Gamma_1) \times Aut(\Gamma_2).$$
The finite group $Aut(\Gamma_1)$
has no $p$-torsion
since $k < p-1$ and so none of the
underlying graphs of
marked graphs in $\tilde X_k$ have any $p$-symmetry.
Similarly,
the finite group $Aut(\Gamma_2)$
has no $p$-torsion since $k > 0$ and so none of the
underlying graphs of
marked graphs in $X_{p-1-k}$ have any $p$-symmetry.
Thus we have that
$$\hat H^*(stab(\delta); \Z_{(p)}) \cong
\hat H^*(N_{\Sigma_p}(\Z/p); \Z_{(p)}) \cong
\hat H^*(\Sigma_p; \Z_{(p)}).$$

Since the above holds for every simplex $\delta$,
we see that
the spectral sequence \ref{e6} has $E_2$ page

$$E_2^{r,s} = \left\{\matrix{
\Z/p \hfill &r=0 \hbox{ and } s=kn \hfill \cr
H^r(\tilde Q_{k} \times Q_{p-1-k}; \Z/p) \hfill &1 \leq r \leq n-3
\hbox{ and } s=kn \hfill \cr
0 \hfill &\hbox{otherwise } \hfill \cr} \right.$$
because the dimension of $\tilde Q_{k} \times Q_{p-1-k}$
is $2p-5$.

Now apply Lemma \ref{tr32} and Lemma \ref{tr33} to conclude 
that all of the groups
$H^1(\tilde Q_{k} \times Q_{p-1-k};\Z/p)$
are zero.
The lemma now follows for $k \in \{1,\ldots,p-2\}.$

\item $k = 0.$  Then the simplices $\delta$
in spectral sequence \ref{e6}
are all in $Q_{p-1}$.
Since only one graph in $Q_{p-1}$ has
$p$-symmetry, namely the graph $\Theta_{p-1}$,
we have
$$stab_{N_{Aut(F_n)}(B_k)}(\delta) = \left\{\matrix{
N_{\Sigma_p}(\Z/p) \times \Sigma_p \hfill & \hbox{if } \delta
\hbox{ has underlying graph }
\Theta_{p-1} \hfill \cr
N_{\Sigma_p}(\Z/p) \hfill &\hbox{otherwise } \hfill \cr} \right.$$
  
Arguments similar to those in Lemma \ref{t5} show that the $E_2$ page has
the form

$$E_2^{r,s} = \left\{\matrix{
\Z/p^2 \oplus \Z/p \hfill & r=0 \hbox{ and } s=0 \hfill \cr
H^r(Q_{p-1}; \Z/p) \hfill &1 \leq r \leq n-2, s = kn \hfill \cr
(j-1)\Z/p \hfill &r=0 \hbox{ and } |s|=nj-1 \hfill \cr
(j+1)\Z/p \hfill &r=0 \hbox{ and } |s|=nj \not = 0  \hfill \cr
0 \hfill &\hbox{otherwise } \hfill \cr} \right.$$

Because Fact \ref{tr35} and Lemma \ref{tr32} give that
$H^2(Q_{p-1}; \Z/p) = 0$,
the
differentials $d_2: E_2^{0,-nj+1} \to E_2^{2,-nj}$ are zero
and we see that the spectral sequence converges at the $E_2$ page.
Lastly, by Lemma \ref{tr32}, $H^1(Q_{p-1}; \Z/p)=0$
and we are done with the case $k=0$.

\item $k = p-1.$
Then the simplices $\delta$
in spectral sequence \ref{e6}
are all in $\tilde Q_{p-1}$.
Now only two graphs in $\tilde Q_{p-1}$ have
$p$-symmetry.
One is the graph $\Theta_{p-1}^1$
where both $*$ and $\circ$
are the left hand vertex of the
$\Theta$-graph.  The other is
the graph $\Theta_{p-1}^2$
where
$*$ is the vertex on the left side of the $\Theta$-graph
and $\circ$ is the vertex on the right side.
Each of these graphs gives a vertex of $\tilde Q_{p-1}$
with $p$-symmetry.

We have
$$stab_{N_{Aut(F_n)}(B_k)}(\delta) = \left\{\matrix{
N_{\Sigma_p}(\Z/p) \times \Sigma_p \hfill
&\hbox{if } \delta \hbox{ has underlying graph }
\Theta_{p-1}^1 \hfill \cr
\hfill 
&\hfill \hbox{or } \Theta_{p-1}^2 \cr
N_{\Sigma_p}(\Z/p) \hfill &\hbox{otherwise } \hfill \cr} \right.$$

Arguments similar to those in Lemma \ref{t5} show that the $E_2$ page has
the form

$$E_2^{r,s} = \left\{\matrix{
\Z/p^2 \oplus \Z/p \hfill & r=0 \hbox{ and } s=0 \hfill \cr
H^r(\tilde Q_{p-1}; \Z/p) \hfill &1 \leq r \leq n-1, s = kn \hfill \cr
(2j-2)\Z/p \hfill & r=0 \hbox{ and } |s|=nj-1 \hfill \cr
(2j+1)\Z/p \hfill & r=0 \hbox{ and } |s|=nj \not = 0  \hfill \cr
0 \hfill &\hbox{otherwise } \hfill \cr} \right.$$
Let $j \geq 2$. Now the differential
$$d_{n} : E_{n}^{0,nj-1}
\to E_{n}^{n,n(j-1)}=0$$
is necessarily trivial and thus
$\hat H^{nj-1}(N_{Aut(F_n)}(B_{p-1}); \Z_{(p)})$ has a filtration
with successive terms 
$$E_2^{0,nj-1}=(2j-2)\Z/p$$
and
$$E_2^{2p-3,n(j-1)}=H^{2p-3}(\tilde Q_{p-1}; \Z/p).$$
Since
$$\hat H^{nj-1}(N_{Aut(F_n)}(B_{p-1}); \Z_{(p)})=
H^{nj-1}(N_{Aut(F_n)}(B_{p-1}); \Z_{(p)})$$
(because $nj-1$ is above the vcd of $N_{Aut(F_n)}(B_{p-1})$),
we can use the K\"unneth formula for the latter cohomology group to
specify the form that the above filtration takes and obtain that
$$\hat H^{nj-1}(N_{Aut(F_n)}(B_{p-1}); \Z_{(p)})
=(2j-2)\Z/p \oplus H^{2p-3}(\tilde Q_{p-1}; \Z/p).$$

The other tricky cohomology group to calculate is
$$\hat H^{-nj+1}(N_{Aut(F_n)}(B_{p-1}); \Z_{(p)})$$ (again,
for $j \geq 2$.)
This can be computed by noting that
$$d_2: E_2^{0,-nj+1} \to E_2^{2,-nj}$$
is zero by Fact \ref{tr31}
and that
$$E_2^{1,-nj}=H^1(\tilde Q_{p-1}; \Z/p)=0$$
by Lemma \ref{tr33}.
The result follows for the case $k=p-1$. \END
\end{enumerate}

\begin{lemma} \label{t7}
$$\hat H^t(N_{Aut(F_n)}(C); \Z_{(p)}) \cong \left\{\matrix{
\Z/p^2 \oplus \Z/p \hfill &t = 0 \hfill \cr
({3k(p-1) \over 2} + 1) \Z/p \hfill &|t| = kn \not = 0 \hfill \cr
({3k(p-1) \over 2} - 1) \Z/p \hfill &|t| = kn-1 \hfill \cr
0 \hfill &\hbox{otherwise } \hfill \cr
} \right.$$
\end{lemma}

\PF The normalizer $N_{Aut(F_n)}(C)$ acts on the contractible space
$X_n^C$ with finite stabilizers and finite quotient.
Hence we can use the equivariant cohomology
spectral sequence
\ref{e2} to calculate the cohomology of the normalizer.
This gives us:
\begin{equation} \label{e7}
E_1^{r,s} = \prod_{[\delta] \in \Delta^r}
\hat H^s(stab_{N_{Aut(F_n)}(C)}(\delta); \Z_{(p)}) \Rightarrow
\hat H^{r+s}(N_{Aut(F_n)}(C); \Z_{(p)})
\end{equation}
where $[\delta]$ ranges over the set $\Delta^r$ of orbits
of $r$-simplices $\delta$ in $X_n^C$.

\begin{claim} \label{tr37}
The quotient space 
$\bigcup_r\Delta^r$ has $3$ vertices and $2$ edges in it.  The vertices
correspond to marked graphs with underlying graphs $\Phi_{n}$,
$\Omega_{n}$, and $\Psi_{n}$.  The two edges come from the
forest collapses of $\Phi_{n}$ to $\Omega_{n}$ or
$\Psi_{n}$.  Pictorially, we have
$$\matrix{ & & \Phi_{n} & & \cr
           & \swarrow & & \searrow & \cr
           \Omega_{n} & & & & \Psi_{n} \cr}$$
\end{claim}

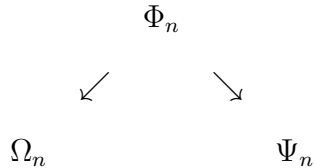
\PF From Theorem \ref{tr4}, if 
$\Gamma_1$ and $\Gamma_2$ are reduced as
$C$-graphs then they are Nielsen equivalent,
up to an equivariant isomorphism
(a basepoint preserving isomorphism.)
That is, any two vertices (marked graphs)
of $X_n^C$ whose underlying graphs are reduced,
can be can be connected by a sequence of Nielsen
transformations.  But the graphs that you get from
$\Psi_n$ by doing Nielsen moves are all isomorphic
to either $\Psi_n$ or $\Omega_n$.

It follows that if $\eta$ is a vertex of $X_n^C$ corresponding to
a reduced marked graph, then the underlying graph of
$\eta$ is either $\Psi_n$ or $\Omega_n$.

It remains to consider which graphs can be blowups
(see Definition \ref{tr22}) of $\Psi_n$ or
$\Omega_n$.  Such a blowup would have a nontrivial
$\Z/p$ action on at least $2p$ of its edges.  From this,
it is not hard to see (using similar methods
to those in Proposition \ref{t1}) that the only possibility
for the underlying graph of such a blowup is
$\Phi_n$. \END
           
Direct examination reveals that the
vertex in $\Delta^0$ corresponding to $\Phi_{n}$ has
automorphism group $\Sigma_p \times \Z/2$.
In the notation used to define $\Phi_{n}$
(refer to the text just above Figure \ref{fig1}),
the $\Sigma_p$ in $\Sigma_p \times \Z/2$
acts on the collections of edges $\{a_i\}$,
$\{b_i\}$, and $\{c_i\}$, respectively,
by permuting their indices.  On the other hand,
the $\Z/2$ in $\Sigma_p \times \Z/2$
fixes the edges $a_i$ and switches the edges
$b_i$ with the edges $c_i$.
The group $C$ is included in $\Sigma_p \times \Z/2$
as the cyclic group generated by the
permutation $(1 2 \ldots p)$ in $\Sigma_p$.
Hence the subgroup of normalizing graph automorphisms in
$\Sigma_p \times \Z/2$ is
$$N_{\Sigma_p}(\Z/p) \times \Z/2.$$
The stabilizer of the vertex in $\Delta^0$
corresponding to $\Phi_n$ has cohomology
$$\hat H^*(N_{\Sigma_p}(\Z/p) \times \Z/2; \Z_{(p)})
=\hat H^*(\Sigma_p; \Z_{(p)}).$$

Similarly, the group of graph automorphisms of
$\Omega_{n}$ is $$\Sigma_p \times \Sigma_p.$$
The group $C$ is included in this as
the subgroup generated by
$$(1 2 \ldots p) \times (1 2 \ldots p).$$
The stabilizer of the vertex in $\Delta^0$ which
corresponds to $\Omega_n$ is the
the normalizer of $C$
in $\Sigma_p \times \Sigma_p$.
Using the definition of a normalizer, we see that
this normalizer is
$$(\Z/p \times \Z/p) \rtimes \Z/(p-1).$$
The generator of $\Z/(p-1)$ acts diagonally on
$(\Z/p \times \Z/p)$
by conjugating the generator of either $\Z/p$
to its $s$-th power for some generator $s$ of
$\F_p^\times$.

The cohomology of $(\Z/p \times \Z/p) \rtimes \Z/(p-1)$
can be calculated in a straightforward way from the cohomology of
$\Z/p \times \Z/p$ using Swan's theorem \ref{tr30}
which gives us that
$$H^*((\Z/p \times \Z/p) \rtimes \Z/(p-1); \Z_{(p)}) =
H^*(\Z/p \times \Z/p; \Z_{(p)})^{\Z/(p-1)}.$$
So we see that
$\hat H^t((\Z/p \times \Z/p) \rtimes \Z/(p-1); \Z_{(p)}) =$
$$\left\{\matrix{
\Z/p^2 \hfill &t = 0 \hfill \cr
(k(p-1) + 1) \Z/p \hfill &|t| = kn \not = 0 \hfill \cr
(k(p-1) - 1) \Z/p \hfill &|t| = kn-1 \hfill \cr
0 \hfill &\hbox{otherwise } \hfill \cr} \right.$$

Lastly, the
group of graph automorphisms of
$\Psi_{n}$ is
$$(\Sigma_p \times \Sigma_p) \rtimes \Z/2.$$
The group $C$ is included in this as
the subgroup generated by
$$(1 2 \ldots p) \times (1 2 \ldots p).$$
The stabilizer of the vertex in $\Delta^0$ which
corresponds to $\Psi_n$ is the
the normalizer of $C$
in $(\Sigma_p \times \Sigma_p) \rtimes \Z/2.$
Using the definition of a normalizer, we see that
this normalizer is
$$(\Z/p \times \Z/p) \rtimes (\Z/2 \times \Z/(p-1)).$$
The generator of $\Z/(p-1)$ acts diagonally on
$(\Z/p \times \Z/p)$
by conjugating the generator of either $\Z/p$
to its $s$-th power for some generator $s$ of
$\F_p^\times$.  The $\Z/2$ acts by exchanging one
$\Z/p$ for the other in $\Z/p \times \Z/p$.
The cohomology of $(\Z/p \times \Z/p) \rtimes (\Z/2 \times \Z/(p-1))$
can be calculated using Swan's theorem \ref{tr30}
which indicates that
$$H^*((\Z/p \times \Z/p) \rtimes (\Z/2 \times \Z/(p-1)); \Z_{(p)}) =
H^*(\Z/p \times \Z/p; \Z_{(p)})^{(\Z/2 \times \Z/(p-1))}.$$
Hence
$\hat H^t((\Z/p \times \Z/p) \rtimes (\Z/2 \times \Z/(p-1));
\Z_{(p)}) =$
$$\left\{\matrix{
\Z/p^2 \hfill &t = 0 \hfill \cr
({k(p-1) \over 2} + 1) \Z/p \hfill &|t| = kn \not = 0 \hfill \cr
({k(p-1) \over 2}) \Z/p \hfill &|t| = kn-1 \hfill \cr
0 \hfill &\hbox{otherwise } \hfill \cr} \right.$$

The two edges in $\Delta^1$ have stabilizers isomorphic to
$N_{\Sigma_p}(\Z/p)$ or $N_{\Sigma_p}(\Z/p) \times \Z/2$.
We omit the argument here, but the stabilizers of
the edges can be found by examining which graph
automorphisms in
$$stab(\Phi_n) = N_{\Sigma_p}(\Z/p) \times \Z/2$$
preserve the relevant forest collapses.
In either case, if we take Farrell
cohomology with $\Z_{(p)}$-coefficients, then both edges have
stabilizers whose cohomology is the same as that of the symmetric
group $\Sigma_p$.

Combining all of this into the spectral sequence \ref{e7} and then
applying the differential on the $E_1$ page, we see that
$$E_2^{r,s} =
\left\{\matrix{
\Z/p^2 \oplus \Z/p \hfill &r = 0, s = 0\hfill \cr
({3k(p-1) \over 2} - 1) \Z/p \hfill &r=0, |s| = kn-1 \hfill \cr
({3k(p-1) \over 2} + 1) \Z/p \hfill &r=0, |s| = kn \not = 0 \hfill \cr
0 \hfill &\hbox{otherwise } \hfill \cr} \right.$$
Thus $E_2=E_\infty$ and $\hat H^*(N_{Aut(F_n)}(C); \Z_{(p)})$
is as stated. \END

\begin{lemma} \label{t8}
$$\hat H^t(N_{Aut(F_n)}(D); \Z_{(p)}) \cong \hat H^t(\Sigma_p
\times \Sigma_p;\Z_{(p)}).$$
$$\hat H^t(N_{Aut(F_n)}(E); \Z_{(p)}) \cong
\hat H^t((\Sigma_p \times \Sigma_p) \rtimes \Z/2;\Z_{(p)})$$
$$ \cong \left\{\matrix{
\Z/p^2 \hfill &t = 0\hfill \cr
([k/2] + 1) \Z/p \hfill &|t| = kn \not = 0 \hfill \cr
([k/2]) \Z/p \hfill &|t| = kn-1 \hfill \cr
0 \hfill &\hbox{otherwise } \hfill \cr
} \right.$$
\end{lemma}

\PF From Lemma \ref{tr11},
$$N_{Aut(F_n)}(D) = N_{\Sigma_p}(\Z/p) \times N_{\Sigma_p}(\Z/p).$$
and
$$N_{Aut(F_n)}(E) =
(N_{\Sigma_p}(\Z/p) \times N_{\Sigma_p}(\Z/p)) \rtimes \Z/2.$$

The detailed description of the cohomology of
$N_{Aut(F_n)}(E)$ in the statement of this lemma
is then obtained by using Swan's theorem \ref{tr30}. \END

\newpage

\chapter{Cohomology of $Aut(F_n)$} \label{c5}

In this chapter, we will use the lemmas of the previous chapter to
calculate the Farrell cohomology of $Aut(F_n)$ with coefficients
in $\Z_{(p)}$, for $n=2(p-1)$.

\begin{thm} \label{t9} Let $p$ be an odd prime, and $n=2(p-1)$.
Then $\hat H^t(Aut(F_{n}); \Z_{(p)})$

\medskip 

$$\matrix{
\cong \left\{\matrix{
\Z/p^2 \oplus p(\Z/p) \hfill &t = 0 \hfill \cr
(p + [{3k \over 2}] - 1)\Z/p \hfill &|t| = kn \not = 0 \hfill \cr
\Z/p \hfill &t = 1 \hfill \cr
0 \hfill &t = kn + 1 > 1 \hfill \cr
([{3k \over 2}] - 1)\Z/p \hfill &t = -kn+1 < 0 \hfill \cr
H^{n-1}(\tilde Q_{p-1}; \Z/p) \oplus
([{3k \over 2}] - 1)\Z/p \hfill &t = kn-1 > 0 \hfill \cr
H^{n-1}(\tilde Q_{p-1}; \Z/p) \hfill &t = -kn-1 < 0 \hfill \cr
H^r(Q_{p}^\omega;\Z/p) \oplus
\sum_{i=0}^{p-1} H^r(\tilde Q_{i} \times Q_{p-1-i};\Z/p)
\hfill &t = kn + r, 2 \leq r \leq n-2 \hfill \cr
} \right. \hfill \cr}$$
\end{thm}

$$\matrix{ \hfill \cr }$$

\bigskip

\PF The cohomology $\hat H^*(Aut(F_{n}); \Z_{(p)})$ will be
calculated using the normalizer spectral sequence \ref{e3}, which has
$E_1$ page
$$E_1^{r,s} = \prod_{(A_0 \subset \cdots \subset A_r) \in |\mathcal{B}|_r}
\hat H^s( \bigcap_{i=0}^r N_{Aut(F_n)}(A_i); \Z_{(p)})
\Rightarrow \hat H^{r+s}(Aut(F_n); \Z_{(p)})$$
where $\mathcal{B}$ denotes the poset of conjugacy classes of
nontrivial elementary
abelian $p$-subgroups of $Aut(F_n)$, and $|\mathcal{B}|_r$ is the
set of $r$-simplices in $|\mathcal{B}|$.
We computed $|\mathcal{B}|$ in Proposition \ref{tr5}.  It is 
$1$-dimensional, so the above spectral sequence is zero except in
the columns $r=0$ and $r=1$.

Recall that the realization $|\mathcal{B}|$ of the poset $\mathcal{B}$ has $p$
path components.  One component just consists of a point corresponding
to the subgroup $A$.  In addition, $p-2$ other components are also
just points corresponding to the subgroups $B_k$ for
$k \in \{1,\ldots,p-2\}$.  Finally, the last component is a
$1$-dimensional simplicial complex corresponding to the subgroups
listed in diagram \ref{e4}, which we duplicate here:
$$\matrix{ & & B_0 & & \cr
           & & \downarrow & & \cr
           & & D & & \cr
           & \nearrow & & \nwarrow & \cr
           B_{p-1} & & & & C \cr
           & \searrow & & \swarrow & \cr
           & & E & & \cr}$$

We have already calculated (in the lemmas of
the previous chapter) the contributions of all of
the vertices in $|\mathcal{B}|$ to the $E_1$ page in \ref{e3}.

The contribution of a $1$-simplex in $|\mathcal{B}|$ can be obtained
by taking the cohomology of the intersections of the normalizers of
the vertices of the $1$-simplex.  Note that each of these
intersections is a finite group (since each is a subgroup of either
the finite group $N_{Aut(F_n)}(D)$ or $N_{Aut(F_n)}(E)$ of normalizing graph automorphisms.)
In this way, we can calculate the (now just Tate) cohomological
contributions of the $1$-simplices in \ref{e4} to be:
$$\matrix{ \hat H^*(N_{Aut(F_n)}(D) \cap N_{Aut(F_n)}(B_0); \Z_{(p)}) \hfill
&= \hat H^*(N_{\Sigma_p}(\Z/p)
 \times N_{\Sigma_p}(\Z/p);\Z_{(p)}) \hfill \cr
 \hfill &= \hat H^*(\Sigma_p \times \Sigma_p;\Z_{(p)}). \hfill}$$
$$\matrix{ \hat H^*(N_{Aut(F_n)}(D) \cap N_{Aut(F_n)}(B_{p-1}); \Z_{(p)}) \hfill
&= \hat H^*(N_{\Sigma_p}(\Z/p)
\times N_{\Sigma_p}(\Z/p);\Z_{(p)}) \hfill \cr
 \hfill &= \hat H^*(\Sigma_p \times \Sigma_p;\Z_{(p)}). \hfill}$$
$$\hat H^t(N_{Aut(F_n)}(D) \cap N_{Aut(F_n)}(C); \Z_{(p)}) =
\hat H^t((\Z/p \times \Z/p) \rtimes \Z/(p-1); \Z_{(p)}) \hfill$$
$$ =
\left\{\matrix{
\Z/p^2 \hfill &t = 0 \hfill \cr
(k(p-1) + 1) \Z/p \hfill &|t| = kn \not = 0 \hfill \cr
(k(p-1) - 1) \Z/p \hfill &|t| = kn-1 \hfill \cr
0 \hfill &\hbox{otherwise } \hfill \cr} \right.$$
$$\matrix{\hat H^*(N_{Aut(F_n)}(E) \cap N_{Aut(F_n)}(B_{p-1}); \Z_{(p)}) \hfill
&= \hat H^*(N_{\Sigma_p}(\Z/p)
\times N_{\Sigma_p}(\Z/p);\Z_{(p)}) \hfill \cr
 \hfill &= \hat H^*(\Sigma_p \times \Sigma_p;\Z_{(p)}).\hfill}$$
$$\hat H^t(N_{Aut(F_n)}(E) \cap N_{Aut(F_n)}(C); \Z_{(p)}) =
\hat H^t((\Z/p \times \Z/p) \rtimes (\Z/2 \times \Z/(p-1));
\Z_{(p)}) \hfill$$
$$ =
\left\{\matrix{
\Z/p^2 \hfill &t = 0 \hfill \cr
({k(p-1) \over 2} + 1) \Z/p \hfill &|t| = kn \not = 0 \hfill \cr
({k(p-1) \over 2}) \Z/p \hfill &|t| = kn-1 \hfill \cr
0 \hfill &\hbox{otherwise } \hfill \cr} \right.$$

(Note that $N_{Aut(F_n)}(E) \cap N_{Aut(F_n)}(B_{p-1})$ is not
$(N_{\Sigma_p}(\Z/p) \times N_{\Sigma_p}(\Z/p)) \rtimes \Z/2$
as you might expect,
because the map that ``flips'' the two sides of $\Psi_n$ is not in
the normalizer $N_{Aut(F_n)}(B_{p-1})$.)

We are now ready to compute the $E_2$ page of the spectral sequence
\ref{e3}.  The contributions coming from the isolated points of
$|\mathcal{B}|$ (i.e., from $A$,$B_1$, $\ldots$, $B_{p-2}$) survive
unaltered from the $E_1$ page.  The contributions from the $1$-dimensional
component
of $|\mathcal{B}|$ pictured in \ref{e4} will be what we concentrate
on from now on.  

First off, it is easy to compute the values for the $E_2$ page in rows $s$
where $s=nj+k$ with $2 \leq k \leq 2(p-2).$  This is because the
$E_1$ page is only nonzero in the column $r=0$ for these rows.  Hence
the entries in $E_1^{0,s}$ necessarily survive to the $E_2$ page and
from there survive to the $E_\infty$ page.  This gives us that 
$\hat H^t(Aut(F_{n}); \Z_{(p)})$ is as the proposition
claims for $t=nj+k$ with $2 \leq k \leq 2(p-2).$

For the rest of our calculations, we will use the fact that the
boundary map on the $E_1$ page is just the restriction map.  From a
comment by Brown in \cite{[B]} on page 286, we know that we can
compute these restriction maps (from normalizers of $p$-subgroups to
finite subgroups of those normalizers) just by looking at the
$E_2$ pages of the various spectral sequences used to compute the
cohomologies of the normalizers (in the lemmas of the previous
chapter.)  This will help us to compute, for any row $s$, the value
$E_2^{0,s}$.  As an example of this, consider the copy
of $\hat H^*(N_{Aut(F_n)}(B_0); \Z_{(p)})$ contained in the column
$E_2^{0,s}$ of our spectral sequence \ref{e3}.   Recall
from Lemma \ref{t6} that the cohomology of
$N_{Aut(F_n)}(B_0)$ is 
$$\hat H^t(N_{Aut(F_n)}(B_0); \Z_{(p)}) = 
\left\{\matrix{
\Z/p^2 \hfill &t = 0 \hfill \cr
(j+1) \Z/p \hfill &|t| = nj \not = 0 \hfill \cr
(j-1) \Z/p \hfill &|t| = nj-1 \hfill \cr
0 \hfill &|t| \equiv 1  \hbox{ } (\hbox{mod } n) \hfill \cr
0 \hfill &t \equiv 2  \hbox{ } (\hbox{mod } n) \hfill \cr
H^r(Q_{p-1};\Z/p) \hfill
&t \equiv r  \hbox{ } (\hbox{mod } n), \hfill \cr
 \hfill &3 \leq r \leq n-2 \hfill \cr
}\right.$$
and this was calculated by looking at a spectral
sequence whose $\mathcal{E}_2$ page
(which we now denote with a script
$E$ as $\mathcal{E}$ to distinguish it
from the spectral sequence \ref{e3} above)
is
$$\mathcal{E}_2^{r,s} = \left\{\matrix{
\Z/p^2 \hfill & r=0 \hbox{ and } s=0 \hfill \cr
H^r(Q_{p-1}; \Z/p) \hfill &1 \leq r \leq 2(p-2), s = kn \hfill \cr
(j-1)\Z/p \hfill & r=0 \hbox{ and } |s|=nj-1  \hfill \cr
(j+1)\Z/p \hfill & r=0 \hbox{ and } |s|=nj \not = 0  \hfill \cr
0 \hfill &\hbox{otherwise } \hfill \cr} \right.$$
From Proposition $(4.6)$ in \cite{[B]}, the restriction
map from $\hat H^*(N_{Aut(F_n)}(B_0); \Z_{(p)})$
to the ring $\frak{H}^*(N_{Aut(F_n)}(B_0); \Z_{(p)})$ 
deriving from the cohomology of the finite subgroups
of $N_{Aut(F_n)}(B_0)$ (see \cite{[B]} for a definition of this ring)
is just the canonical surjection coming from the
vertical edge homomorphism
$$\hat H^t(N_{Aut(F_n)}(B_0); \Z_{(p)}) \to
\mathcal{E}_2^{0,t} =
\hat H^t(N_{\Sigma_p}(\Z/p) \times N_{\Sigma_p}(\Z/p);\Z_{(p)})$$
where
$$\hat H^t(N_{\Sigma_p}(\Z/p) \times N_{\Sigma_p}(\Z/p);\Z_{(p)})=
\left\{\matrix{
\Z/p^2 \hfill & t = 0 \hfill \cr
(j+1)\Z/p \hfill & |t| = nj \not = 0  \hfill \cr
(j-1)\Z/p \hfill & |t|= nj-1  \hfill \cr
0 \hfill &\hbox{otherwise } \hfill \cr} \right.$$
Recall that in order to know the column $E_2^{0,s}$ of
\ref{e3}, we want to calculate the
restriction map
$$\matrix{
\hfill \hat H^*(N_{Aut(F_n)}(B_0); \Z_{(p)}) & \to &
\hat H^*(N_{Aut(F_n)}(B_0) \cap N_{Aut(F_n)}(D)) \hfill \cr
\hfill & = &
\hat H^*(N_{\Sigma_p}(\Z/p) \times N_{\Sigma_p}(\Z/p);\Z_{(p)}) \hfill \cr}$$
So from the above work this restriction map is
just the vertical edge homomorphism calculated above.
We
can similarly
look at the $\mathcal{E}_2$-pages of the other spectral sequences in
chapter \ref{c4} to calculate the
restriction maps from the cohomologies of
$N_{Aut(F_n)}(B_{p-1})$ and $N_{Aut(F_n)}(C)$ to their finite
subgroups
$$N_{Aut(F_n)}(B_{p-1}) \cap N_{Aut(F_n)}(D), \hbox{ }
N_{Aut(F_n)}(B_{p-1}) \cap N_{Aut(F_n)}(E),$$
$$N_{Aut(F_n)}(C) \cap N_{Aut(F_n)}(D), \hbox{ and }
N_{Aut(F_n)}(C) \cap N_{Aut(F_n)}(E).$$
On the other hand, the restrictions from
$N_{Aut(F_n)}(D)$ and $N_{Aut(F_n)}(E)$ to their subgroups are easy to
compute as $N_{Aut(F_n)}(D)$ and $N_{Aut(F_n)}(E)$ are just well
known finite groups.

After calculating all of
the terms in the column $E_2^{0,s}$ of \ref{e3}
as  above, the value $E_2^{1,s}$ 
is found from 
an Euler characteristic argument using $E_1^{0,s}$ and
$E_1^{1,s}$.  That is, if $s \not = 0$ then $E_1^{0,s}$ and
$E_1^{1,s}$ are the only nonzero terms on the row $s$ and they
are both $\F_p$-vector spaces.  As the boundary map goes from
$E_1^{0,s} \to E_1^{1,s}$, this yields
$dim_{\F_p}(E_1^{0,s}) - dim_{\F_p}(E_1^{1,s})
= dim_{\F_p}(E_2^{0,s}) - dim_{\F_p}(E_2^{1,s})$
by the standard Euler characteristic argument.

Assume $s=-kn+1<0$ (the case $s=kn-1>0$
follows similarly, with the only exception being the additional
summand of $H^{n-1}(\tilde Q_{p-1}; \Z/p)$ that needs
to be dealt with.)  Then
$E_2^{0,s} = (k+[k/2]-1)\Z/p$.  Now $k-1$ of these $\Z/p$'s come from
$0$-cocycles that are summations of cocycles from:
\begin{itemize}
\item The portion of $\hat H^s(N_{Aut(F_n)}(B_{p-1});\Z_{(p)})$ coming from the
graph $\Omega_n$. (In the spectral sequence we used to calculate
$\hat H^s(N_{Aut(F_n)}(B_{p-1});\Z_{(p)})$ in part $3.$ of
Lemma \ref{t6}, this was the contribution given by the
stabilizer of the graph $\Theta_{p-1}^2$.)
\item $\hat H^s(N_{Aut(F_n)}(D);\Z_{(p)}).$
\item $\hat H^s(N_{Aut(F_n)}(B_{0});\Z_{(p)}).$
\item The portion of $\hat H^s(N_{Aut(F_n)}(C);\Z_{(p)})$ coming from the
graph $\Omega_n$ (in the spectral sequence we used to calculate
$\hat H^s(N_{Aut(F_n)}(C));\Z_{(p)})$.)
\end{itemize}
The other $[k/2]$ of the $\Z/p$'s come from $0$-cocycles that are
summations of cocycles from:
\begin{itemize}
\item The portion of $\hat H^s(N_{Aut(F_n)}(B_{p-1});\Z_{(p)})$ coming from the
graph $\Psi_n$. 
(In the spectral sequence we used to calculate
$\hat H^s(N_{Aut(F_n)}(B_{p-1});\Z_{(p)})$ in part $3.$ of
Lemma \ref{t6}, this was the contribution given by the
stabilizer of the graph $\Theta_{p-1}^1$.)
\item $\hat H^s(N_{Aut(F_n)}(E);\Z_{(p)}).$
\item The portion of $\hat H^s(N_{Aut(F_n)}(C);\Z_{(p)})$ coming from the
graph $\Psi_n$ (in the spectral sequence we used to calculate
$\hat H^s(N_{Aut(F_n)}(C);\Z_{(p)})$.)
\end{itemize}
This allows us to find that the contribution of the component of
$|\mathcal{B}|$
corresponding to diagram \ref{e4} in the row $s=-kn+1<0$
is
$E_2^{0,s} = (k+[k/2]-1)\Z/p.$
Since
$dim_{\F_p}(E_1^{0,s})-dim_{\F_p}(E_1^{1,s})=$
$$(4(k-1)+{3k(p-1) \over 2}+[k/2]-1) -
(3(k-1)+{3k(p-1) \over 2})=k+[k/2]-1,$$
we find that $E_2^{1,s}=0$.

Assume $s=kn>0$ or $s=-kn<0$.  Then
$E_2^{0,s} = (k+[k/2]+p-1)\Z/p$.  First, observe that $p-1$ of these
$\Z/p$'s come from the points in $|\mathcal{B}|$ corresponding to
$A$,$B_1$, $\ldots$, $B_{p-2}$.
Also, $k-1$ of these $\Z/p$'s come from summations of
cocycles in
$H^s(N;\Z_{(p)})$ for $N$ equal to
$N_{Aut(F_n)}(B_0)$, $N_{Aut(F_n)}(B_{p-1})$, $N_{Aut(F_n)}(C)$, and
$N_{Aut(F_n)}(E)$ where in all cases the cohomology
came from the graph $\Omega_n$ (in the spectral sequences
used to calculate the cohomology groups of the various
normalizers.)
In addition, $[k/2]$ of the $\Z/p$'s come from $0$-cocycles that are
summations of cocycles from $H^s(N;\Z_{(p)})$ for $N$ equal to
$N_{Aut(F_n)}(B_{p-1})$, $N_{Aut(F_n)}(C)$, and $N_{Aut(F_n)}(E)$ where in all cases the cohomology
came from the graph $\Psi_n$.
Finally, one cocycle comes from summing up cocycles from
$H^s(N;\Z_{(p)})$ for $N$ equal to $N_{Aut(F_n)}(B_0)$, $N_{Aut(F_n)}(B_{p-1})$,
$N_{Aut(F_n)}(C)$, $N_{Aut(F_n)}(D)$, and $N_{Aut(F_n)}(E)$.
This allows us to find that the contribution of the component of $B$
corresponding to diagram \ref{e4} in the row $s=kn>0$ or $s=-kn<0$
is
$E_2^{0,s} = (k+[k/2]+p-1)\Z/p.$
Note that because
$dim_{\F_p}(E_1^{0,s})-dim_{\F_p}(E_1^{1,s})=$
$$(4k+{(3k+2)(p-1) \over 2}+[k/2]+5) -
(3k+{3k(p-1) \over 2}+5)=k+[k/2]+p-1,$$
we find that $E_2^{1,s}=0$.

Next, $E_2^{0,0}$ is readily computed to be $\Z/p^2 \oplus p(\Z/p)$,
where $p-1$ of the $\Z/p$'s come from the isolated points of
$|\mathcal{B}|$ and $\Z/p^2 \oplus \Z/p$ comes from vertices in
the $1$-dimensional component of $|\mathcal{B}|$.
Finally, $E_2^{0,1} = \Z/p$ where
the $\Z/p$ corresponds to edges of $|\mathcal{B}|$ which give
cohomology classes in $E_1^{0,1}$ that are not mapped
onto by cohomology classes in $E_1^{0,0}$.\END

\chapter{Notes about the usual cohomology of $Aut(F_{2(p-1)}$)} \label{c6}

We now state what our calculations using Farrell
cohomology 
tell us about the usual cohomology of
$Aut(F_{n})$, $n=2(p-1)$.

\begin{prop} \label{t10}
$$\matrix{
\hfill H^t(Aut(F_{n}); \Z_{(p)}) &=& &\left\{\matrix{
\hat H^t(Aut(F_{n}); \Z_{(p)}) \hfill &t > 4p-6 \hfill \cr
H^t(Q_{n}; \Z_{(p)}) \hfill &0 \leq t \leq 2p-3 \hfill \cr}
\right. \hfill \cr }$$
\end{prop}

\PF Since the virtual cohomological dimension of
$Aut(F_{n})$ is $4p-6$, it follows directly (see \cite{[B]}) that
$H^t(Aut(F_{n}); \Z_{(p)}) = \hat H^t(Aut(F_{n});
\Z_{(p)})$
for $t > 4p-6$ and that the sequence
$$H^{m}(\Gamma); \Z_{(p)}) \to H^{m}(Aut(F_{n});
\Z_{(p)})
\to \hat H^{m}(Aut(F_{n}); \Z_{(p)}) \to 0$$
is exact for $m = 4p-6$, where $\Gamma$ is any torsion free subgroup
of $Aut(F_{n})$ of finite index
and
where the first map is the transfer map.

To establish the final part of the proposition,
we can use the spectral sequence \ref{e1} to calculate
$H^t(Aut(F_{n}); \Z_{(p)})$ for $0 \leq t \leq 2p-3$.
We want to show that for $0 < s \leq 2p-3$,
the rows $\tilde E_2^{t,s}$ on the $E_2$
page of this spectral sequence, are all zero.  If this is true,
then the cohomology of $Aut(F_{n})$ for $0 \leq s \leq 2p-3$
comes just from the terms $H^s(Q_{n}; \Z_{(p)})$ on the
horizontal axis of the $E_2$ page, which proves the proposition.

On the $E_1$-page of the spectral sequence \ref{e1}, the only
simplices that contribute to terms above the horizontal axis
are ones with ``$p$-symmetry''.  In other words, let $X_s$ be
the $p$-singular locus of $X_{n}$, or the set of all simplices
in $X_{n}$
such that $p$ divides $|stab(\delta)|$.  Then $Aut(F_{n})$
also acts on $X_s$ and there is a spectral sequence
analogous to \ref{e1} which calculates the equivariant
cohomology $H_{Aut(F_{n})}^*(X_s;\Z_{(p)})$.  Above the
horizontal axes, the $E_1$-pages of the
two spectral sequences (i.e., the
ones corresponding to $X_{n}$ and $X_s$) are the
same.  In particular, the boundary maps $E_1^{r,s} \to E_1^{r+1,s}$
are the same so that the entries $E_2^{r,s}$ of the two
spectral sequences are identical for $s > 0$.  So we have reduced
the problem to showing that
for $0 < s \leq 2p-3$, the rows $\tilde E_2^{*,s}$ on the
$E_2$-page of this spectral sequence for
calculating $H_{Aut(F_{n})}^*(X_s;\Z_{(p)})$, are all zero.
But note that $X_s$ separates into $p$ disjoint,
$Aut(F_{n})$-invariant (not necessarily connected)
subcomplexes corresponding to the $p$ path components
of $|\mathcal{B}|$.  (These path components were mentioned in
the proof of Theorem \ref{t9}.)
So we see that
$$H_{Aut(F_{n})}^*(X_s; \Z_{(p)})
= \prod_Y \hat H_{Aut(F_{n})}^*(Y; \Z_{(p)})$$
where $Y$ ranges over the $p$ distinct, disjoint,
$Aut(F_{n})$-invariant subcomplexes.  Hence the $E_1$
and $E_2$
pages can be calculated separately for each $Y$.  
Recall that we are only concerned with calculating the $E_2$ page in
the horizontal rows between $1$ and $2p-3$.  The only simplices of
$X_{n}^s$ that contribute terms to these rows are the ones with
dihedral symmetry.  Now only one of the subcomplexes $Y$ has dihedral
symmetry in it and is relevant;  
namely, the subcomplex $Y_A$ which has
marked graphs with underlying graph the
rose $R_{n}$ in it. (This subcomplex $Y_A$ corresponds to the
fixed point 
subcomplex $X_{n}^A$.  More precisely, $Y_A$ is
$Aut(F_{n}) \cdot X_{n}^A$, several disjoint copies of
$X_{n}^A$ grouped together.)  Let $\Xi_p$ be the
$1$-skeleton of the cone over a
$p$-gon.  Note that $\Xi_p$ has dihedral symmetry; moreover, marked
graphs with underlying graph of the form $\Xi_p \vee \Gamma_{p-2}$
are in $Y_A$ and have dihedral symmetry, where $\Gamma_{p-2}$ is any
graph with fundamental group isomorphic to $F_{p-2}$.

Now consider a specific path component of $Y_A$, say $X_{n}^A$.
There is an action of $A \cong \Z/p$ on each simplex of
$X_{n}^A$.  In $\cite{[K-V]}$, Krstic and Vogtmann define an
equivariant deformation retract $L_A$ of $X_{n}^A$ by only
keeping
the ``essential graphs'' in $X_{n}^A$.  Among the graphs that are
in $X_{n}^A - L_A$ are all of
the ones with dihedral symmetry, along with various other
inessential graphs.
The ones with
dihedral symmetry are of the form
$\Xi_p \vee \Gamma_{p-2}$ mentioned above (where the
basepoint is in $\Gamma_{p-2}$ and the wedge doesn't necessarily
take place at the basepoint.)
We can use the same idea that Krstic and Vogtmann used
and define $L_A'$ to be the subcomplex of
$X_{n}^A$ obtained by collapsing out graphs with
dihedral symmetry.  By the same poset lemma used in \cite{[K-V]},
the obvious collapse from $X_{n}^A$
to $L_A'$ is an equivariant deformation retraction.
Observe that
$X_{n}^A \supseteq L_A' \supseteq L_A$.

The contribution of the equivariant cohomology of $Y_A$
to $\tilde E_1^{r,s}$ is the same as that for the $E_1$-page of the
spectral sequence used to calculate the equivariant cohomology
of $N_{Aut(F_n)}(A)$ acting on $X_{n}^A$.
(This follows immediately, as $Y_A$ is
$Aut(F_{n}) \cdot X_{n}^A$ in $X_s$ so that we see that the
$E_1$-pages of the equivariant spectral sequences used to calculate
$\hat H_{Aut(F_{n})}^*(Y_A; \Z_{(p)})$ and
$\hat H_{N_{Aut(F_n)}(A)}^*(X_{n}^A; \Z_{(p)})$
are identical.)

To avoid
confusion, denote the $E_1$-page of
the spectral sequence used to calculate
$\hat H_{N_{Aut(F_n)}(A)}^*(X_{n}^A; \Z_{(p)})$
by  $\mathcal{E}_1^{r,s}$.
We want to show that $\mathcal{E}_1^{*,s}=0$
for $0 < s \leq 2p-3$.  For $0 < s \leq 2p-3$ and
$s \not = 4k$, the row $\mathcal{E}_1^{*,s}$ is all zero.
On the other hand, for $0 < s=4k < 2p-3$, direct examination
(since only simplices $\sigma$ in $X_{n}^A - L_A'$ contribute to the
row $\mathcal{E}_1^{*,s}$ and $H^s(stab(\sigma); \Z_{(p)}) = \Z/p$
for such simplices) reveals that the
row $\mathcal{E}_1^{*,s}$ is just the relative cochain
complex $C^*(X_{n}/N_{Aut(F_n)}(A),L_A'/N_{Aut(F_n)}(A);\Z/p)$.
Because the homotopy
in the deformation retraction from $X_{n}^A$ to $L_A'$
is $N_{Aut(F_n)}(A)$-equivariant,
we see that the
relative cohomology groups
$$H^*(X_{n}/N_{Aut(F_n)}(A),L_A'/N_{Aut(F_n)}(A);\Z/p)$$
are all zero.  Hence all of the rows $\mathcal{E}_2^{r,s}$
are zero for $0 < s \leq 2p-3$, which completes
the proof. \END

\chapter{The Farrell cohomology of $Aut(F_l)$ for $l < 2(p-1)$} \label{c7}

For $l < p-1$ an easy spectral sequence
argument yields that $\hat H^*(Aut(F_l); \Z_{(p)}) =0$.
We have already remarked that Glover and Mislin's work
in \cite{[G-M]} directly implies
that $\hat H^*(Aut(F_{p-1}); \Z_{(p)}) = \hat H^*(\Sigma_p; \Z_{(p)})$,
$\hat H^*(Aut(F_{p}); \Z_{(p)}) = 3\hat H^*(\Sigma_p; \Z_{(p)})$,
and we have noted that Yu Qing Chen's work in \cite{[C]}
shows that 
$$\hat H^*(Aut(F_{p+1}); \Z_{(p)}) = 4\hat H^*(\Sigma_p; \Z_{(p)})$$
and that 
$\hat H^t(Aut(F_{p+1}); \Z_{(p)}) = 5\hat H^t(\Sigma_p; \Z_{(p)})
\oplus \hat H^{t-4}(\Sigma_p; \Z_{(p)})$.

In this chapter, we show how to 
calculate $\hat H^*(Aut(F_{l}); \Z_{(p)})$
for $p \leq l < n$, where $n=2(p-1)$.
We do this by modifying in a direct manner
the arguments we have made
in the previous chapters to calculate
$\hat H^*(Aut(F_{2(p-1)}); \Z_{(p)})$. 
(It should be noted that
using the normalizer spectral sequence perhaps amounts
to overkill in these situations, however, as
the realization $|\mathcal{B}|$ of conjugacy classes
of elementary abelian $p$-subgroups is just
a collection of $l-p+3$ points.)

As before, one of the elements of $\mathcal{B}$
corresponds to the subgroup $A \cong \Z/p$ of
$Aut(F_l)$ obtained by rotating the first
$p$ petals of the rose $R_l$.
The other elements of $\mathcal{B}$
correspond to subgroups $B_k \cong \Z/p$
for $k \in \{0, \ldots, l-p+1\}$.
The subgroup $B_k$ is obtained,
as before, by rotating the $p$-edges
of the $\Theta$-graph in the middle of
$R_k \vee \Theta_{p-1} \vee R_{l-p+1-k}$.

As in Lemma \ref{t5}, we see that
$$N_{Aut(F_l)}(A) \cong
N_{\Sigma_p}(\Z/p) \times
((F_{l-p} \times F_{l-p}) \rtimes (\Z/2 \times Aut(F_{l-p}))$$
where the semi-direct product corresponds to the obvious
action of $\Z/2 \times Aut(F_{l-p})$
on $F_{l-p} \times F_{l-p}$.
Let $\langle \omega \rangle  \cong \Z/2$ be the subgroup of
$Aut(F_{l-p+2})$ corresponding to the
action given by
switching the first two petals of the rose $R_{l-p+2}$.
Accordingly,
$$N_{Aut(F_l)}(A) \cong C_{Aut(F_l)}(A) \cong
(F_{l-p} \times F_{l-p}) \rtimes (\Z/2 \times Aut(F_{l-p}).$$
As in Definition \ref{tr17}, let
$X_{l-p+2}^{\langle \omega \rangle}$ be the fixed point set of $\omega$
in $X_{l-p+2}$.  Then let $X_{l-p+2}^\omega$
be the deformation retract of $X_{l-p+2}^{\langle \omega \rangle}$
obtained by collapsing out inessential edges of
marked graphs.  Finally, define $Q_{l-p+2}^\omega$
to be the quotient of $X_{l-p+2}^\omega$
by $N_{Aut(F_{l-p+2})}(\omega)$.

Putting all of this together, and using the same
methods as those used in Lemma \ref{t5},
Lemma \ref{t6} and Theorem \ref{t9}, we see that

$$\hat H^t(N_{Aut(F_l)}(A); \Z_{(p)})
= \left\{\matrix{
\Z/p \hfill &t \equiv 0  \hbox{ } (\hbox{mod } n) \hfill \cr
H^r(Q_{l-p+2}^\omega; \Z/p) \hfill &t \equiv r  \hbox{ } (\hbox{mod } n), \hfill \cr
 &2 \leq r \leq 2(l-p) \hfill \cr
0 \hfill &\hbox{otherwise } \hfill \cr} \right.$$

and that for $k \in \{0,\ldots,l-p+1\}$,

$\hat H^t(N_{Aut(F_l)}(B_k); \Z_{(p)}) = $
$$\matrix{
\hfill \cr
\left\{\matrix{
\Z/p \hfill &t \equiv 0  \hbox{ } (\hbox{mod } n)  \hfill \cr
H^r(\tilde Q_{k} \times Q_{p-1-k};\Z/p) \hfill
&t \equiv r  \hbox{ } (\hbox{mod } n), \hfill \cr
\hfill &2 \leq r \leq 2(l-p)+1 \hfill \cr
0 \hfill &\hbox{otherwise } \hfill \cr} \right.
\hfill \cr
\hfill & \hfill \cr}$$

so we have that for $p$ odd and $l \in \{p, \ldots, 2p-3\}$,

$\hat H^t(Aut(F_{l}); \Z_{(p)}) = $
$$\matrix{
\hfill \cr
\left\{\matrix{
(l-p+3)\Z/p \hfill &t \equiv 0  \hbox{ } (\hbox{mod } n)  \hfill \cr
\hfill & \hfill \cr
\sum_{k=0}^{l-p+1} H^r(\tilde Q_{k} \times Q_{l-p+1-k};\Z/p)
\hfill &t \equiv r  \hbox{ } (\hbox{mod } n), \hfill \cr
\oplus H^r(Q_{l-p+2}^\omega;\Z/p) \hfill &2 \leq r \leq 2(l-p)+1 \hfill \cr
\hfill & \hfill \cr
0 \hfill &\hbox{otherwise } \hfill \cr} \right. \hfill \cr}$$

\bigskip

The aforementioned results of Glover and Mislin
and Chen give more information than our results
here, however, as they actually explicitly calculate the
cohomology groups
$$\sum_{k=0}^{l-p+1} H^r(\tilde Q_{k} \times Q_{l-p+1-k};\Z/p)
\oplus H^r(Q_{l-p+2}^\omega;\Z/p)$$
that arise in the above formula in their cases.
	
\part{Modified degree theorem} \label{p2}

\chapter{The modified degree theorem} \label{c9}

Before reading Part \ref{p2}, we strongly recommend that the reader
study the entirety of \cite{[H-V]} by Hatcher and Vogtmann.
This part is essentially a modification of their results on the space
$X_n$ to the space $\tilde X_n$ and is not meant to be read independently
of their work.

Our goal is to use a modified version of Hatcher and Vogtmann's
Degree Theorem
and roughly similar computer calculations to those used by
Hatcher and Vogtmann in
\cite{[V]},
to establish Fact \ref{tr31} which stated
that
$$H^2(\tilde Q_{p-1}; \Z/p)=0$$
and was needed in the last part of the proof of
Lemma \ref{t6}. 

In this chapter and the following one, we will use the notation
developed in Proposition \ref{t11} for describing marked graphs in
$\tilde X_n$.  (For this and the following chapter, take $n$ be
be just any positive
integer, rather than $2(p-1)$ as we often did in
Part \ref{p1}.)

Define the {\em degree} of a marked graph
$$\matrix{\hfill (\alpha,f): R_{n} \coprod I \to \Gamma_{n}\hfill \cr}$$
or, equivalently, the degree of the
underlying graph $\Gamma_{n}$ of the marked graph,
to be
$$deg(\alpha,f) = deg(\Gamma_{n}) = \sum_{v \not = *} |v| - 2,$$
where the sum is over all vertices of $\Gamma_{n}$
except the basepoint, and where the valence $|v|$ of a vertex
$v \not = \circ$
is the number of oriented edges starting at $v$.
The valence of the vertex $\circ$ is defined to be
one plus the number of oriented edges starting at $\circ$.
This definition is similar to the one given in \cite{[H-V]},
but the intuition (which we do not attempt to
make precise) is that
we treat the vertex $\circ$ as if it signifies that
a small germ of an extra edge
is coming into that vertex.
Throughout this section,
our modified definitions of degree, split degree,
canonical splitting, etc., will be motivated by this notion of
thinking that the vertex $\circ$ denotes the germ of another edge
entering that vertex.

Let $\tilde X_{n,k}$, $\tilde Q_{n,k}$, and $\tilde \A_{n,k}$
be the subspaces of
$\tilde X_{n}$, $\tilde Q_{n}$, and $\tilde \A_{n}$,
respectively,
where only marked graphs of degree at most $k$ are considered.

Our goal is to show that the arguments which Hatcher and Vogtmann
used in \cite{[H-V]} carry over in an essentially unchanged way
to this new context.  In particular, we want to prove the following
analogue of their degree theorem:

\begin{thm} \label{t13}  A piecewise linear map $f_0:D_k \to \tilde
\A_{n}$ is homotopic to a map $f_1:D_k \to \tilde \A_{n,k}$
by a homotopy $f_t$ during which degree decreases monotonically,
i.e., if $t_1 < t_2$ then $deg(f_{t_1}(s)) \geq deg(f_{t_2}(s))$
for all $s \in D_k$.
\end{thm}

As in \cite{[H-V]}, the immediate corollary is

\begin{cor} \label{t14} The pair $(\tilde \A_{n}, \tilde \A_{n,k})$
is $k$-connected.
\end{cor}

In Part \ref{p3} we will show that $\tilde \A_{n}$ is contractible, so
that Corollary \ref{t14} is equivalent to the statement that
$\tilde \A_{n,k}$ is $(k-1)$-connected.

Hatcher and Vogtmann prove their degree theorem by
using various homotopies to deform the
underlying graphs of marked graphs
$$\alpha: R_n \to \Gamma_n.$$
We will use these same homotopies
to deform the marked graphs
$$\matrix{\hfill (\alpha,f): R_{n} \coprod I \to \Gamma_{n}\hfill \cr}$$
that appear in the context of $\tilde \A_n$.
In a remark
following ``Stage 1: Simplifying the critical point'' in
\cite{[H-V]}, Hatcher and Vogtmann
mention that it is obvious
where their homotopies of the
underlying graph $\Gamma_n$ send the basepoint $*$
and the marking $\alpha$.
Our task is to decide where these homotopies send the
extra point $\circ$ on the graph.  It will then be
clear where the homotopies send the path $f$
from $*$ to $\circ$ in the graph $\Gamma_n$.

As in the case of $\A_n$, graphs in $\tilde \A_n$ come
equipped with a ``height function'' measuring distance
to the basepoint.  A point in the graph is a
{\em critical point} if near that point there is more than one
way to travel ``down'' toward the basepoint.

A few notes are in order about specific parts of the
paper by Hatcher and Vogtmann and how these should be
modified:
\begin{enumerate}
\item Canonical splittings.  A procedure called ``canonical
splitting'' is defined in \cite{[H-V]} which is used to decrease
the degrees of graphs (with lengths assigned to their edges) in
a canonical way.  In our new case, think of $\circ$ as denoting
the germ of an extra edge attached to the graph at that point.
A canonical splitting should split this extra
edge down to the next critical point, assuming $\circ$ is not
already a critical
point itself.  In other words, a canonical splitting should move
$\circ$ down to the next critical point or the basepoint, provided
$\circ$ is not already a critical point.  Otherwise, do the
canonical splitting as described in \cite{[H-V]}.
See Figure \ref{fig5} for examples.

\bigskip

\input{fig5.pic}
\begin{figure}[here]
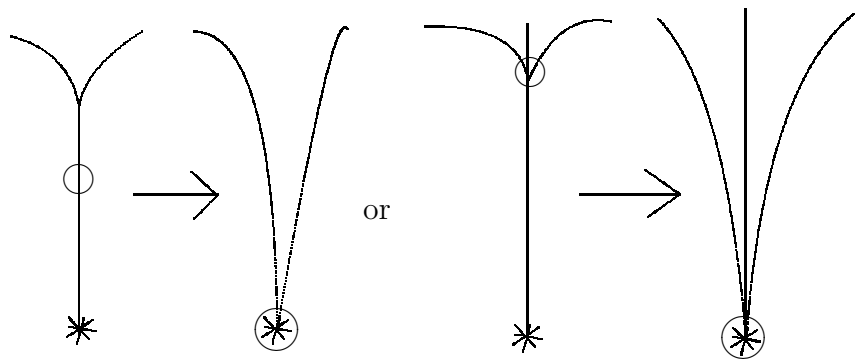

\caption{\label{fig5} Canonical splittings}
\end{figure}

\item Sliding $\epsilon$-cones.  This is another procedure defined
by Hatcher and Vogtmann to reduce degrees, and can be
directly modified by taking $\circ$ to signify that the
germ of an extra edge is attached there.  In other words, $\circ$
signifies the attaching point of an extra branch lying above that
vertex.  Just like the attaching points of other branches
lying above the critical point, it can be perturbed slightly
off of the critical point.  The result of this is that
$\circ$ is perturbed slightly downward off the critical point.

\item Codimension.  Codimension is defined without trouble
as before.  Namely, the {\em codimension} of a point on the
graph is one less than the number of downward directions
from that point.  The {\em codimension} of a graph is the
sum of the codimensions of its critical points.

Figure \ref{fig6} below shows some examples of this when
we are calculating the codimension of $\circ$ in a graph.

\bigskip

\input{fig6.pic}
\begin{figure}[here]
\caption{\label{fig6} Some examples of codimensions concerning $\circ$}
\end{figure}

To illustrate this, consider the second example in the above
figure.   There are two different ways to
start an edgepath going down from the vertex $\circ$.  If both of
these edges leading down from $\circ$, say $a$ and $b$,
went directly to the basepoint, then the
condition that they have equal length would be given by
$2-1=1$ linear equations among the lengths of the edges of
the graph.  In other words, the vertex $\circ$ has codimension
1 and the graph is contained in the hyperplane $l(a)=l(b)$ of
codimension 1.  Thus a stratum where the codimension is equal
to $i$ is contained in a critical plane of codimension $i$.

\item Lemma 4.1.  This lemma is fine as stated because the
same proof with its statements
about codimension remains true
when applied to the new space $\tilde A_n$. 

\item Lemma 4.2.  During the homotopies used in Lemma 4.2, only
the lengths of the edges of the underlying graph $\Gamma_{n}$
are perturbed.  The combinatorial
structure of the graph is not changed at
all.  Hence it is clear where $\circ$ and the path $f$ from $*$ to
$\circ$ are sent during these homotopies.  

\item Complexity.
As before, we think of $\circ$ as meaning that the germ of an
extra edge is attached there.  Previously, Hatcher and Vogtmann
defined a connecting path in the graph as a downward path from one
critical point to another.  We just modify this by saying that
the extra germ attached at $\circ$ counts as the beginning of
a downward path that came from a critical point (although, of course,
there is no critical point at the top of this extra path, which
is why we just consider $\circ$ to be the germ of an attached edge,
rather than a whole extra edge attached.)  With this convention,
complexity $c_s$ is defined as in \cite{[H-V]} as the number of connecting paths
in the graph.  Also, $e_s$ is the number of connecting paths without
critical points in their interiors.  See Figure \ref{fg1}.

\input{fg1.pic}
\begin{figure}[here]
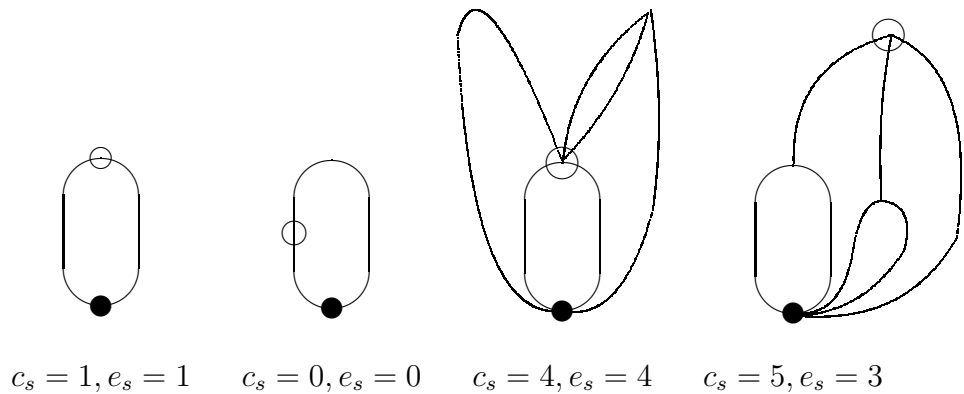

\caption{\label{fg1} Complexity examples with $\circ$}
\end{figure}

With this definition of complexity in mind, the section
``canonical splitting and extension to a neighborhood'' can
be done as before.  For the section
``reducing complexity by sliding in the $\epsilon$-cones'',
consider $\circ$ as giving an attaching point
$\alpha_j$ of a branch $\beta_j$.  Now argue as
directed in \cite{[H-V]}.

\item Reducing the degree.  The final section, ``Reducing the degree'',
presents no further difficulties.  Define the
split degree of a graph as the degree of its
canonical splitting. All of the codimension
arguments check out because Lemma 4.2 still holds in this new context.
\end{enumerate}

Using the above guidelines, the proof of
Theorem \ref{t13} is immediate from the
work of Hatcher and Vogtmann in \cite{[H-V]}.

\chapter{Using the modified degree theorem} \label{c10}

We now describe how the modified degree theorem, Theorem \ref{t13},
can be used to calculate $H_2(\tilde Q_{n}; \Z/p)=0$.
The idea
to use low dimensional calculations and the degree theorem to calculate
rational homology groups of spaces having to do with $Aut(F_n)$
and $Out(F_n)$ comes from
Hatcher and Vogtmann's work in \cite{[V]}.  In addition,
much of the specific methods used to
calculate $H_2(\tilde Q_{n}; \Z/p)=0$
here come from their work, the only real exceptions being
modifying the methods to work with different coefficient
rings like $\Z_{(p)}$ or $\Z/p$
and the idea that the degree theorem
could be modified to calculate homology groups of $\tilde Q_{p-1}$.

As a warm-up exercise to the real problem, we
prove the following proposition using
the original Degree Theorem techniques of
Hatcher and Vogtmann from \cite{[H-V]} and
\cite{[V]}.  Their work directly implies the truth
of the proposition, and it is only because the
results in \cite{[V]} were stated using
rational -- rather than $\Z/p$ --
coefficients that we need to write anything at all
here.

\begin{prop} \label{tr36}
$$H^2(Aut(F_{p-1}); \Z/p)=0.$$
\end{prop}

\PF By dual universal coefficients,
$$H^2(Aut(F_{p-1}); \Z/p) \cong
Hom(H_2(Aut(F_{p-1})),\Z/p) \oplus
Ext(H_1(Aut(F_{p-1}),\Z/p).$$
From Theorem \ref{tr15},
$$Ext(H_1(Aut(F_{p-1}),\Z/p)=0.$$
It remains to show that
$$H_2(Aut(F_{p-1}); \Z/p) = 0.$$
In addition, note from Lemma \ref{tr32} that
$$H_2(Aut(F_{p-1}); \Z/p) = H_2(Q_{p-1}; \Z/p).$$

From Corollary $3.2$ in \cite{[H-V]},
$X_{p-1,3}$ is $2$-connected and so
we can use a spectral sequence to calculate
$H_2(Aut(F_{p-1}); \Z/p)$ by taking a resolution $F \to \Z$ of
$\Z$ as a $\mathbb{Z}(Aut(F_{p-1}))$-module and considering two different
filtrations of the complex $F \otimes C_*(X_{p-1,3}; \Z/p)$.
Using one filtration, one obtains a
spectral sequence that satisfies
$$E_{r,s}^2 \Rightarrow
H_{r+s}(Aut(F_{p-1}) ; \Z/p), \hbox{ for } r+s \leq 2,$$
because $X_{p-1,3}$ is $2$-connected.
Then using the other filtration gives
a spectral sequence
$$E_{r,s}^1 = \prod_{[\delta] \in \Delta^r}
H_s(stab_{Aut(F_{p-1})}(\delta); \Z/p)
\Rightarrow H_{r+s}(Aut(F_{p-1}); \Z/p),
\hbox{ for } r+s \leq 2,$$
where $[\delta]$ ranges over the set $\Delta^r$ of orbits
of $r$-simplices $\delta$ in $X_{p-1,3}$.
As always, for $s > 0$ we
need only be concerned with simplices with $p$-symmetry.  Moreover,
as we are truncating at degree $3$, there are not that
many simplices with $p$-symmetry.
Explicitly, we are looking for graphs with degree no greater
than $3$ and no more than $p-1$ holes.
The only relevant
ones occur when $p=3,5$
and consist of the graphs $\Theta_2$
and $\Theta_4$, respectively.  Each gives a vertex
of $Q_{p-1,3}$ with $p$-symmetry in the case $p=3$
or $p=5$, respectively.
Fortunately, the homology of the stabilizers
of these simplices is $0$ for $0 < t < 3$,  Because
the
differential on the $E^r$-page has bidegree $(-r,r-1)$,
only the horizontal axis of the
$E^1$-page of the spectral sequence is
relevant to our calculation of $H_2(Aut(F_{p-1}); \Z/p)$.
Hence we see that
$$H_2(Q_{p-1,3}; \Z/p) = H_2(Aut(F_{p-1}); \Z/p).$$

Now follow similar reasoning
and consider two different spectral sequences
corresponding to two different filtrations of
$F \otimes \tilde C_*(\mathcal{C}; \Z/p)$ where
$\mathcal{C}$ is the mapping cone of the inclusion
$X_{p-1,3} \to X_{p-1}$.  Because
$(X_{p-1}, X_{p-1,3})$ is $3$-connected,
standard arguments with these spectral sequences imply
that $$H_3(Q_{p-1}, Q_{p-1,3}; \Z/p) = 0.$$
Similarly, $$H_2(Q_{p-1}, Q_{p-1,2}; \Z/p) = 0.$$

Hence the long exact sequence of the triple
$(X_{p-1}, X_{p-1,3}, X_{p-1,2})$
shows that $H_2(Q_{p-1,3}, Q_{p-1,2}; \Z/p) = 0$.
Consequently, the long exact sequence of the pair
$(X_{p-1,3}, X_{p-1,2})$ yields that the sequence
$$H_3(Q_{p-1,3}, Q_{p-1,2}; \Z/p) \to
H_2(Q_{p-1,2}; \Z/p)
\to H_2(Q_{p-1,3}; \Z/p) \to 0$$
is exact.  Now from \S 5 of \cite{[H-V]},
the space $Q_{p-1,2}$ is contractible.
So 
$$H_2(Q_{p-1,2}; \Z/p)=0$$
and the above exact sequence implies that
$$H_2(Q_{p-1,3}; \Z/p)=0.$$
As we have already shown that the latter cohomology
group is isomorphic to 
$$H_2(Aut(F_{p-1}); \Z/p),$$
we are done. \END

In the following proposition, we apply the same methods
that were used in Proposition \ref{tr36}, but this time
using the modified degree theorem:

\begin{prop} \label{t15} 
\begin{enumerate}
\item $H^2(\tilde Q_{p-1}; \Z/p) \cong 
H_2(\tilde Q_{p-1}; \Z/p) \cong H_2(\tilde Q_{p-1,3}; \Z/p)$, and
\item $H_3(\tilde Q_{p-1,3}, \tilde Q_{p-1,2}; \Z/p) \to
H_2(\tilde Q_{p-1,2}; \Z/p)
\to H_2(\tilde Q_{p-1,3}; \Z/p) \to 0$ is exact.
\end{enumerate}
\end{prop}

\PF Let $G = F_{p-1} \rtimes Aut(F_{p-1})$, thought of as a subgroup
of
$$N_{Aut(F_n)}(B_{p-1}) \cong
N_{\Sigma_p}(\Z/p) \times (F_{p-1} \rtimes Aut(F_{p-1})).$$
Recall from chapter \ref{c8} and in particular Proposition \ref{t12}
that $G$ acts on $\tilde X_{p-1}$ with finite quotient
$\tilde Q_{p-1}$ and
finite stabilizers.

By dual universal coefficients,
$$H^2(G; \Z/p) \cong
Hom(H_2(G),\Z/p) \oplus
Ext(H_1(G),\Z/p).$$

Since from Lemma \ref{tr33} $H_1(G; \Z)$ is
a finite group whose order is a power of $2$,
$$Hom(H_2(G),\Z/p)
\cong Hom(H_2(G; \Z/p), \Z/p)
\cong H_2(G; \Z/p)$$
and
$$Ext(H_1(G),\Z/p)=0.$$
Hence $H^2(G; \Z/p) \cong 
H_2(G; \Z/p)$ by our earlier calculation
with dual coefficients.

Using the standard equivariant spectral sequences in cohomology
and homology for $G$ acting on $\tilde X_{p-1}$, we note that
$$H^2(G; \Z/p) \cong H^2(\tilde Q_{p-1}; \Z/p)$$
and
$$H_2(G; \Z/p) \cong H_2(\tilde Q_{p-1}; \Z/p)$$
by the standard arguments we have used throughout this
paper.  (In other words, we need only pay attention to
the simplices of $\tilde X_{p-1}$ with $p$-symmetry if
we are concerned with terms above the horizontal axis
in either spectral sequence.  The only
marked graphs corresponding to vertices in $\tilde X_{p-1}$
which have $p$-symmetry, are the ones whose underlying
graphs are $\Theta_{p-1}$ with either both $*$ and $\circ$
on the same side of the $\Theta$-graph or on opposite sides
of the $\Theta$-graph.
Each contributes one copy of the (co)homology of the
symmetric group $\Sigma_p$ to the vertical axis of
its spectral sequence.  Since the (co)homology of
$\Sigma_p$ is $0$ for $t \in \{1, \ldots, 2p-4\}$, it
is not relevant to our calculation of
$H^2(G; \Z/p)$ or $H_2(G; \Z/p)$.)

Thus we have
$$H^2(\tilde Q_{p-1}; \Z/p) \cong
H_2(G; \Z/p) \cong H_2(\tilde Q_{p-1}; \Z/p),$$
which establishes part of $1.$ in the proposition.

From Theorem \ref{t13}, $\tilde X_{p-1,3}$ is $2$-connected and so
we can use a spectral sequence to calculate
$H_2(G; \Z/p)$ by taking a resolution $F \to \Z$ of
$\Z$ as a $\mathbb{Z}G$-module and considering two different
filtrations of the complex $F \otimes C_*(\tilde X_{p-1,3}; \Z/p)$.
Using one filtration, one obtains a
spectral sequence that satisfies
$$E_{r,s}^2 \Rightarrow H_{r+s}(G ; \Z/p), \hbox{ for } r+s \leq 2,$$
because $\tilde X_{p-1,3}$ is $2$-connected.
Then using the other filtration gives
a spectral sequence
$$E_{r,s}^1 = \prod_{[\delta] \in \Delta^r}
H_s(stab_G(\delta); \Z/p) \Rightarrow H_{r+s}(G; \Z/p),
\hbox{ for } r+s \leq 2,$$
where $[\delta]$ ranges over the set $\Delta^r$ of orbits
of $r$-simplices $\delta$ in $\tilde X_{p-1,3}$.
As always, for $s > 0$ we
need only be concerned with simplices with $p$-symmetry.  Moreover,
as we are truncating at degree $3$, there are not that
many simplices with $p$-symmetry.  In particular, the only relevant
ones occur when
\begin{itemize}
\item $p=3$
where two vertices of $\tilde Q_{p-1,3}$ contribute
$\Sigma_3$-symmetry.  (Explicitly,
since $p-1=3-1=2$ in this case, we are looking for
graphs with degree no greater than 3, with $3$-symmetry,
and with just two holes in them.  The only possibilities
are $\Theta_2$ with both $*$ and $\circ$ the same vertex,
and $\Theta_2$ with $*$ on one side and $\circ$ on the
other.  Graphs like $\Xi_3$ cannot occur because the ranks of
their fundamental groups are too high.)
\item $p=5$
where one vertex of $\tilde Q_{p-1,3}$ contributes
$\Sigma_5$-symmetry.  (Explicitly,
since $p-1=5-1=4$ in this case, we are looking for
graphs with degree no greater than 3, with $5$-symmetry,
and with just four holes in them.  The only possibility
is $\Theta_4$ with both $*$ and $\circ$ the same vertex.)
\end{itemize}
Fortunately, the homology of the stabilizers
of these simplices is $0$ for $0 < t < 3$. 
Since
the
differential on the $E^r$-page has bidegree $(-r,r-1)$,
only the horizontal axis of the
$E^1$-page of the spectral sequence is
relevant to our calculation of $H_2(G; \Z/p)$.
Hence we see that
$$H_2(\tilde Q_{p-1,3}; \Z/p) = H_2(G; \Z/p).$$
Part $1.$ of the proposition follows.

For part $2.$ we also follow similar reasoning
and consider two different spectral sequences
corresponding to two different filtrations of
$F \otimes \tilde C_*(\mathcal{C}; \Z/p)$ where
$\mathcal{C}$ is the mapping cone of the inclusion
$\tilde X_{p-1,3} \to \tilde X_{p-1}$.  Now because
$(\tilde X_{p-1}, \tilde X_{p-1,3})$ is $3$-connected,
standard arguments with these spectral sequences imply
that $$H_3(\tilde Q_{p-1}, \tilde Q_{p-1,3}; \Z/p) = 0.$$
Similarly, $$H_2(\tilde Q_{p-1}, \tilde Q_{p-1,2}; \Z/p) = 0.$$

Hence the long exact sequence of the triple
$(\tilde X_{p-1}, \tilde X_{p-1,3}, \tilde X_{p-1,2})$
shows that $H_2(\tilde Q_{p-1,3}, \tilde Q_{p-1,2}; \Z/p) = 0$.
Consequently, the long exact sequence of the pair
$(\tilde X_{p-1,3}, \tilde X_{p-1,2})$ yields $2.$
of the proposition. \END

The vertices of the space $\tilde Q_{p-1,2}$ are given by
underlying basepointed graphs
$\Gamma_{p-1}$ with one extra distinguished vertex $\circ$
where
$\pi_1(\Gamma_{p-1}) \cong F_{p-1}$ and $deg(\Gamma_{p-1}) \leq 2$.
As a notational device, when writing a graph $\Gamma_{p-1}$
omit any loops at the basepoint.  So if a graph $\Gamma$
with $\pi_1(\Gamma) \cong p-1-k$ is
written down and asserted to represent a vertex of
$\tilde Q_{p-1,2}$, then the vertex it is representing is
the graph $\Gamma_{p-1}$ obtained from $\Gamma$ by adding $k$
loops at the basepoint.  Note that $deg(\Gamma_{p-1})=deg(\Gamma)$
since edges beginning at the basepoint are not relevant to
degree computations.

An example of a $2$-cycle in $\tilde Q_{p-1,2}$ is illustrated
below in Figure \ref{fg2}.
In the figure, a filled dot represents the
basepoint $*$ and a hollow dot represents the
other distinguished point $\circ$.  The cycle forms a $2$-sphere,
and the hash marks on the boundary of the
pictured disk shows how the $2$-sphere is formed.

\bigskip

\newpage  

\input{fg2.pic}
\begin{figure}[here]
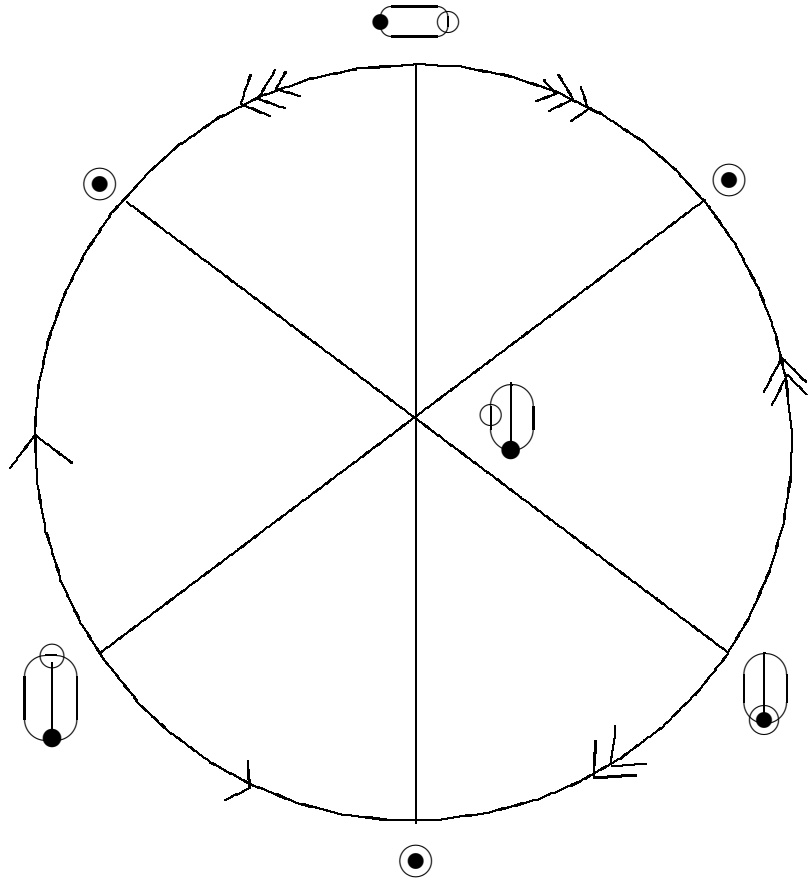

\caption{\label{fg2} The $2$-sphere in $\tilde Q_{p-1,2}$}
\end{figure}

\begin{prop} \label{t16} $H_2(\tilde Q_{p-1}; \Z/p) = 0.$ \end{prop}

\PF We will show that the complex $\tilde Q_{p-1,2}$
deformation retracts to the $2$-sphere given
by Figure \ref{fg2} so that
$H_2(\tilde Q_{p-1,2}; \Z/p) = \Z/p$.
Then an explicit element will be found in
in $$H_3(\tilde Q_{p-1,3}, \tilde Q_{p-1,2}; \Z/p)$$ 
which maps onto the generator of
$$H_2(\tilde Q_{p-1,2}; \Z/p) = \Z/p.$$
So by $2.$ of Proposition \ref{t15},
$H_2(\tilde Q_{p-1,3}; \Z/p)$ will be zero and
then the corollary will follow from $1.$ of Proposition \ref{t15}.

The proof of this adopts much of the same logic as that used by
Hatcher and 
Vogtmann in \cite{[V]}.
In particular, we make use of a
cubical structure on some simplices in
$\tilde Q_{p-1,2}$ and $\tilde Q_{p-1,3}$ and the notion
of ``plusfaces'' and ``minusfaces'' for this cubical structure and
more generally for the simplicial structure.

Let us illustrate this with an example.  Note that
$\tilde Q_{p-1,3}$ is
a $3$-dimensional complex (because we can collapse at most
three $3$-valence vertices or two $3$-valence vertices and
a bivalent $\circ$, to the basepoint.)

\begin{figure}[here]
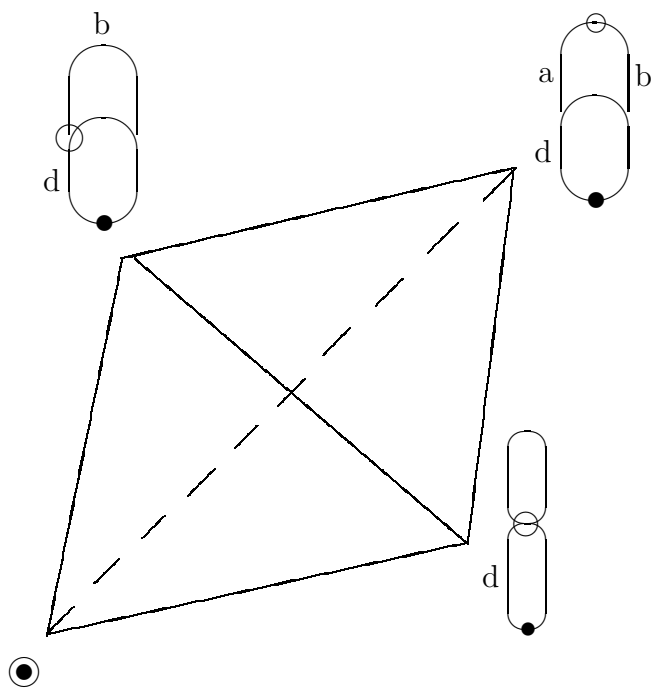

\caption{\label{fg3} One of $6$ $3$-simplices making up the sample cube.}
\end{figure}

\begin{figure}[here]
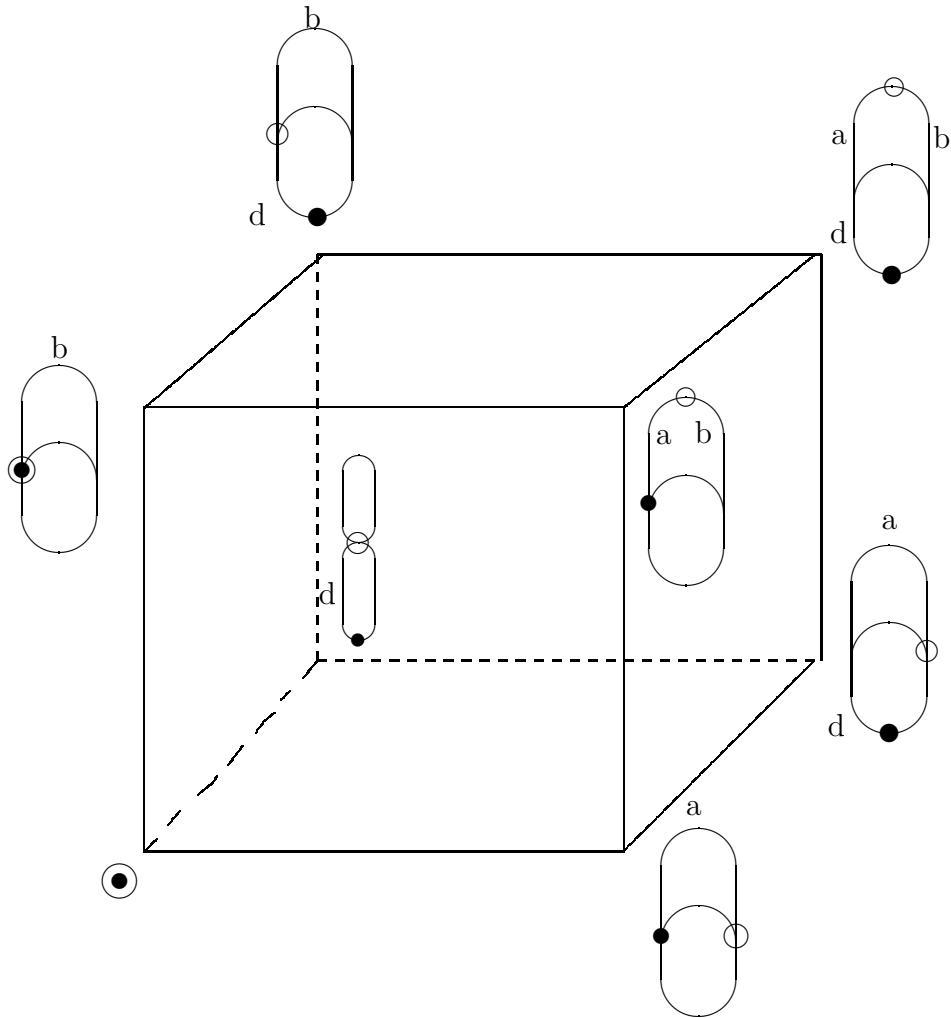

\caption{\label{fg4} A sample cube in $\tilde Q_{p-1,3}$}
\end{figure}

Consider a graph $\Gamma'$ with 4 vertices: $*$, $\circ$.
and two other valence $3$ vertices $x$ and $y$.  The graph $\Gamma'$
has $5$ (unoriented) edges: $a, b, c, d \hbox{ and } e$.
The edge $a$ connects $\circ$ and $x$, $b$ connects $\circ$ and $y$,
$c$ connects $x$ and $y$, $d$ connects $*$ and $x$,
and $e$ connects $*$ and $y$.
Define a forest $F$ in $\Gamma'$ by $F = \{ a,b,d \}$.  There
are $6$ $3$-simplices in $\tilde Q_{p-1,2}$ corresponding
to the forest $F$.  Each corresponds to some collapse of
the edges in $F$ in a particular order.  So, for example, one
particular $3$-simplex comes from collapsing first $a$,
then $b$, and then $d$.  These $6$ $3$-simplices all fit together
in a cube.  Refer to Figures \ref{fg3} and \ref{fg4}
for illustrations of this behavior. 

\newpage

\vspace*{\fill}
\input{fg3.pic}
\vspace*{\fill}

\newpage

\input{fg4.pic}

The $3$ faces of the cube that are adjacent to the vertex $\Gamma'$
are its {\em plusfaces}, and the $3$ faces of the cube that are
adjacent to the vertex $R_{p-1}$ are its {\em minusfaces}.  For the
$3$-simplex pictured above the cube, the $3$ faces of it that are
adjacent to the vertex $\Gamma'$ will be called its {\em plusfaces},
and the $1$ face of it that is not adjacent to the
vertex $\Gamma'$ will be called its {\em minusface}.
(Sometimes, the $3$ one-dimensional edges of the cube adjacent to the
vertex $\Gamma'$ will also be called (lower dimensional) plusfaces
of the cube; and the $3$ edges of the above pictured $3$-simplex
that are adjacent to the vertex $\Gamma'$
will be called plusfaces of the $3$-simplex.)

The subset $F$ of edges of $\Gamma'$ actually gives
the interior of a $3$-dimensional cube as
pictured above because all of the $6$ $3$-simplices that form
the cube are distinct.  This in turn is true because no
graph automorphism of $\Gamma'$ takes the forest $F$ to itself
except the trivial graph automorphism.  Despite this, however,
parts of the boundary of the cube are self-identified.  For example,
the plusface corresponding to $\{a,b\}$ ``folds over'' along
its diagonal so that the square forming the plusface is
glued together to form just one $2$-simplex, a triangle.
All of the other faces (the two other plusfaces and the three
other minusfaces) of the pictured cube are actually squares, and
not self-identified into triangles.

In general, suppose we have a
graph $\Gamma$ (now not necessarily the specific
graph $\Gamma'$ referred to in the above two figures)
which has degree $3$
and a forest
$F = \{ a,b,c \}$ in $\Gamma$.  Then the pair $(\Gamma,F)$ will give
a cube in $\tilde Q_{p-1}$ if no nontrivial graph automorphism of $\Gamma$ sends
$F$ to itself.  Note that if this is the case, then even though
$(\Gamma,F)$ gives a cube, its faces might be
identified or glued to each other in various ways.
This can happen with both the plusfaces and the
minusfaces.
For example, say $\hat \Gamma$ is the graph obtained from
$\Gamma$ by collapsing the edge $a$.  Then if a nontrivial
graph automorphism of $\hat \Gamma$ switches $b$ and $c$, the minusface
(which would have been a square if not for this nontrivial
graph automorphism)
of the cube corresponding to $\hat \Gamma$ is glued to itself along
a diagonal of the square and so is really just a $2$-simplex
or a triangle.

\begin{figure}[here]
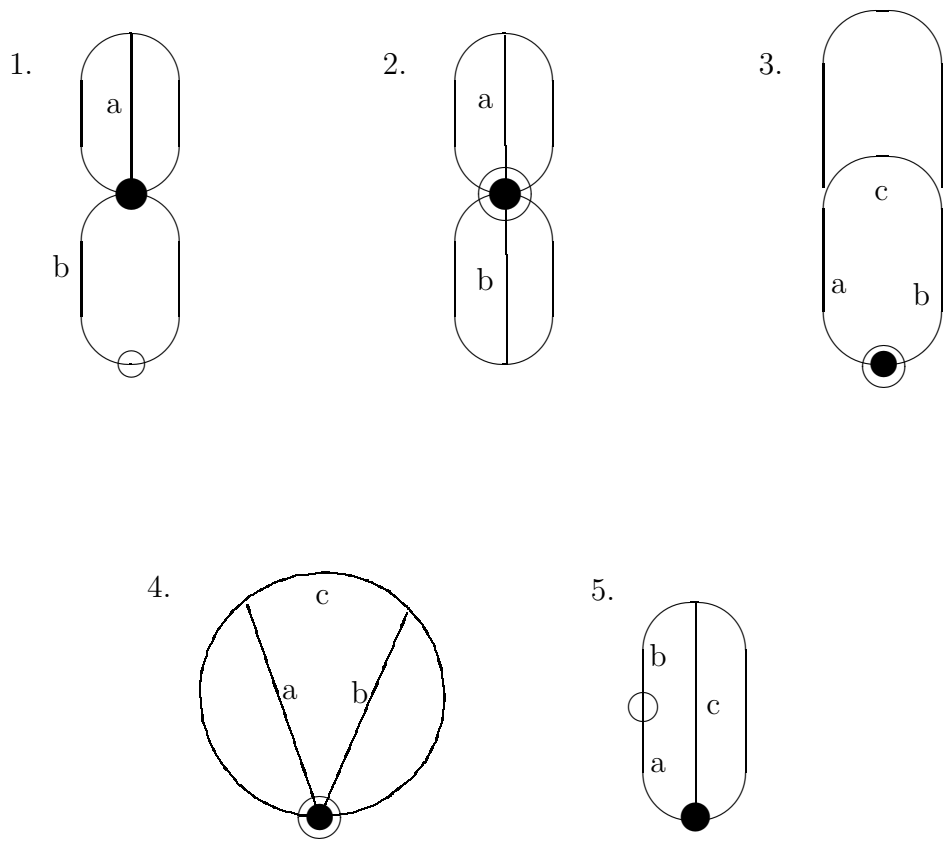

\caption{\label{fg5} Graphs giving $2$-simplices}
\end{figure}

Our first goal is to
show that
$\tilde Q_{p-1,2}$
deformation retracts to the $2$-sphere pictured in
Figure \ref{fg2}.  We do this by listing all of the
graphs that give $2$-simplices in
Figure \ref{fg5} and then handling these graphs in order.

\bigskip

\input{fg5.pic}

\begin{enumerate}
\item This graph only has one maximal subforest $\{a,b\}$ so that the
corresponding square (i.e., two $2$-simplices that join together to
form a square) has free plusfaces.  So we can collapse this
square away.
\item This graph also only has one maximal subforest $\{a,b\}$
and moreover there is an automorphism of the graph that switches $a$
and $b$.  So this graph just contributes one $2$-simplex, the
diagonal of which is automatically a free plusface.
\item The $2$-simplex corresponding to $\{a,b\}$ has a free
diagonal plusface and so can be removed.  In addition, the
square corresponding to $\{a,c\}$ has free plusface $c$.
\item Use exactly the same argument as that for $3.$
\item The three squares that this graph contributes join
together to form the $2$-sphere in Figure \ref{fg2}, as
previously mentioned.
\end{enumerate}

To complete the proof, it remains
to find a specific relative cycle 
in $$H_3(\tilde Q_{p-1,3}, \tilde Q_{p-1,2}; \Z/p)$$ 
which maps onto the generator of
$$H_2(\tilde Q_{p-1,2}; \Z/p) = \Z/p.$$
We use the graph $\Gamma'$, mentioned earlier
in this proof, to construct the relative cycle.
Basically, the cycle is formed by joining together
the cubes corresponding to the
subforests $\{a,b,d\}$, $\{a,c,d\}$,
and $\{a,c,e\}$.  

The plusface corresponding to $\{a,b\}$ of
$\{a,b,d\}$ folds back onto itself along its
diagonal and so it not free.  The plusface
$\{a,d\}$ of $\{a,b,d\}$ connects up with
the plusface $\{a,d\}$ of $\{a,c,d\}$.  The
remaining plusface $\{b,d\}$ of $\{a,b,d\}$
joins up with the plusface $\{a,e\}$ of
$\{a,c,e\}$.

Moreover, the plusface $\{a,c\}$ of
$\{a,c,d\}$ joins up with the plusface
$\{a,c\}$ of $\{a,c,e\}$.  Finally, the
last remaining plusface $\{c,d\}$ of
$\{a,c,d\}$ connects with the last
remaining plusface $\{c,e\}$ of
$\{a,c,e\}$.

The minusface obtained by collapsing $a$
in $\{a,b,d\}$ is the same as that
obtained by collapsing $a$ in
$\{a,c,d\}$.  The one obtained by
collapsing $b$ in $\{a,b,d\}$ is the
same as what we get if we collapse
$a$ in $\{a,c,e\}$.  The remaining
minusface of $\{a,b,d\}$ procured by collapsing
$d$ is one of the three squares pictured in
Figure \ref{fg2} and corresponds to the
subforest $\{a,b\}$ of graph $5.$

Following the same logic, we see that
the minusface of $\{a,c,d\}$ corresponding to
collapsing the edge $c$, is the same
as that of $\{a,c,e\}$ obtained by collapsing
$c$.  The square $\{a,c\}$ of the
graph $5.$ is now seen to be the remaining minusface
acquired from $\{a,c,d\}$ by collapsing $d$.
Similarly, the square $\{b,c\}$ of the graph
$5.$ is the remaining minusface of $\{a,c,e\}$
which we get if we collapse the edge $e$.

The above argument shows that is it at least
possible for the three cubes $\{a,b,d\}$,
$\{a,c,d\}$, and $\{a,c,e\}$ to join together
to form a solid ball with boundary the
sphere pictured in Figure \ref{fg2}.  (Our only remaining concern
might be sign considerations.)
Specific signed computations using the
18 $3$-simplices (6 for each cube) involved
show that they glue together as
we have claimed. \END

\part{Contractibility of fixed point sets} \label{p3}

\chapter{Preliminaries} \label{c11}

Before reading Part \ref{p3}, we strongly recommend that the reader
study the entirety of \cite{[K-V]} by Krstic and Vogtmann.
This part is essentially a modification of their results on
fixed point spaces of outer space to fixed point spaces
of auter space,
and is not meant to be read independently of their work.

We show that for a finite subgroup $G$ of $Aut(F_n)$, the fixed
point set $X_n^G$ is contractible.  This was shown for subgroups
of $Out(F_n)$ by Krstic and Vogtmann in \cite{[K-V]} and our
work here borrows their basic approach and methods (which in
turn are similar to the methods used by Culler and
Vogtmann in \cite{[C-V]} to show that
outer space is contractible).

In particular, Krstic and Vogtmann define a complex $L_G$ of
``essential marked $G$-graphs'' that the fixed point set
$X_n^G$ in outer space deformation retracts to.  Then they
order the reduced marked $G$-graphs in $L_G$ using
a norm $\|\cdot\|_{out}$.  Using this norm
to determine which reduced marked $G$-graph
should be considered next, Krstic and
Vogtmann performed a transfinite induction argument to
show that $L_G$ is contractible, by
building $L_G$ up as the union of stars
of reduced marked $G$-graphs.

We will follow a similar approach, and define norms
$$\|\cdot\|_{aut}
\hbox{ and }
\|\cdot\|_{tot}=\|\cdot\|_{out} \times \|\cdot\|_{aut}$$
to order the reduced marked essential $G$-graphs
in auter space.  For technical reasons,
$\|\cdot\|_{tot}$ will be the appropriate
norm to use when performing the transfinite
induction argument to show the contractibility
of the corresponding $L_G$ in auter space.

The norm $\|\cdot\|_{out}$ was defined by Krstic and
Vogtmann as follows.  Order the set $\mathcal{W}$ of
conjugacy classes of elements of $F_n$ as
$\mathcal{W} = \{w_1, w_2, \ldots \}.$
Totally order $\Z^\mathcal{W}$ by the lexicographic order.
 Let
$\sigma = [s, \Gamma]$ be a marked graph and
define $\|\sigma\|_{out} \in \Z^\mathcal{W}$
by letting $(\|\sigma\|_{out})_i$ be the sum
over all $x \in G$ of the
lengths in $\Gamma$ of the reduced loops (given by the
marking $s$) corresponding to $xw_i$.  Equivalently,
they define an absolute value
$|\cdot|_{out} \in \Z^\mathcal{W}$ on the edges
of $\Gamma$ and set
$$\|\sigma\|_{out} = \frac{1}{2}\sum_{e \in E(\Gamma)} |e|_{out}.$$
The $i$th coordinate of $|e|_{out}$ is simply the
sum for all $x \in G$ of the contributions
of $e$ or $\bar e$ to the loop $xw_i$ in $\Gamma$.
In other words, it is the sum over all $x \in G$ of
the number of times $e$ or $\bar e$ appears in the
cyclically reduced edge path representing $xw_i$.
For $A,B \subseteq E(\Gamma)$ define
$(A.B)_{out} \in \Z^\mathcal{W}$
to be the function whose $i$th coordinate is
the sum over all $x \in G$ of the number of
times $a\bar b$ or $b\bar a$ appears in the
reduced loop in $\Gamma$ corresponding to $xw_i$.
Finally, for $C \subseteq E(\Gamma)$, define
$|C|_{out}$ inductively by the formula
$$\matrix{\hfill |A\coprod B|_{out} = |A|_{out} + |B|_{out} - 2(A.B)_{out} \hfill \cr}$$
for disjoint subsets $A$ and $B$ of $E(\Gamma)$.
Note that with the above definition,
$|A|_{out}=(A.E(\Gamma)-A)_{out}=|E(\Gamma)-A|_{out}$.

The corresponding quantities for $Aut(F_n)$ are defined
in much the same way, the basic difference being that
we think of the lengths of reduced paths rather than reduced
loops.  Order $F_n$ as $F_n = \{ \alpha_1, \alpha_2, \ldots \}$,
and give $\Z^{F_n}$ the lexicographic order.
For a finite subgroup $G$ of $Aut(F_n)$, consider
a basepointed marked $G$-graph $\sigma = [s,\Gamma]$.
Define the norm $\|\sigma\|_{aut} \in \Z^{F_n}$
to be $|G| \cdot L$, where $L : F_n \to \Z$ is
the Lyndon length function of the marked graph.
In other words, the basepointed marked graph $\sigma$ corresponds
to an action of $F_n$ on a rooted $\Z$-tree $T$.
Define $$L(\alpha_i) = \{\hbox{the distance } \alpha_i  \hbox{ moves
the root of } T \}.$$
Equivalently, the $i$th coordinate of $\|\sigma\|_{aut}$
is the sum over all $x \in G$ of the lengths in
$\Gamma$ of the reduced (but {\em not} cyclically reduced)
paths corresponding to $x\alpha_i \in \pi_1(\Gamma,*)$.

As before in the case of $Out(F_n)$, we can
define an absolute value
$|\cdot|_{aut} \in \Z^{F_n}$ on the edges
of $\Gamma$ and set
$$\|\sigma\|_{aut} = \frac{1}{2}\sum_{e \in E(\Gamma)} |e|_{aut}.$$
The $i$th coordinate of $|e|_{aut}$ is simply the
sum of for all $x \in G$ of the contributions
of $e$ or $\bar e$ to the reduced
(but not cyclically reduced) path $x\alpha_i$ in $\pi_1(\Gamma,*)$.
Hence it is the sum over all $x \in G$ of
the number of times $e$ or $\bar e$ appears in the
reduced edge path representing $x\alpha_i$.
For $A,B \subseteq E(\Gamma)$ define
$(A.B)_{aut} \in \Z^\mathcal{W}$
to be the function whose $i$th coordinate is
the sum over all $x \in G$ of the number of
times $a\bar b$ or $b\bar a$ appears in the
reduced path in $\Gamma$ corresponding to $x\alpha_i$.
Finally, for $C \subseteq E(\Gamma)$, define
$|C|_{aut}$ inductively by the formula
$$\matrix{\hfill |A\coprod B|_{aut} = |A|_{aut} + |B|_{aut} - 2(A.B)_{aut} \hfill \cr}$$
for disjoint subsets $A$ and $B$ of $E(\Gamma)$.
In contrast to the case with $Out(F_n)$ the formula
$|A|_{aut}=(A.E(\Gamma)-A)_{aut}$
certainly does not hold any longer.

Our final norm $\|\cdot\|_{tot}$ is just the
product of the previous two.  That is, let
$\sigma = [s,\Gamma]$ be a basepointed marked
$G$-graph for a finite subset $G$ of $Aut(F_n)$
and totally order $\Z^\mathcal{W} \times \Z^{F_n}$
by the lexicographic order.  Define
$\|\sigma\|_{tot} \in \Z^\mathcal{W} \times \Z^{F_n}$
as $\|\sigma\|_{tot} = \|\sigma\|_{out} \times \|\sigma\|_{aut}$,
where to calculate $\|\sigma\|_{out}$ we just
forget that $\Gamma$ has a basepoint.
The functions $|e|_{tot}$, $(A.B)_{tot}$, and
$|A|_{tot}$ are defined similarly.

For a vertex $v$, let $E_v$ be the set of oriented
edges ending at $v$.  
The notion of ideal edges is defined as in \cite{[K-V]},
with the exception that if the ideal edge
$\alpha  \subseteq E_*$ then condition $(i)$ of their definition
should be changed to:
$$(i) \hbox{ } card(\alpha) \geq 2 \hbox{ and } card(E_* - \alpha) \geq 1.$$
That is, ideal edges at the basepoint can contain all except
one of the edges of $E_*$.  The definition of blowing up an
ideal edge is taken exactly as defined in \cite{[K-V]}.  Hence
if we are blowing up an ideal edge $\alpha \subseteq E_*$ then
we are pulling the edges of $\alpha$ away from the basepoint
along a new edge $e(\alpha)$ we just constructed.  If
$card(E_* - \alpha) = 1$, this will result in a graph whose
basepoint has valence $2$.

Let $\alpha$ be an ideal edge of $\sigma = [s, \Gamma]$ and
$\sigma^{G\alpha} = [s^{G\alpha}, \Gamma^{G\alpha}]$
be the result of blowing up the ideal edge $\alpha$.
Then it is easy to show that $|\alpha|_{aut}$ in
$\Gamma$ is equal to $|e(\alpha)|_{aut}$ in $\Gamma^{G\alpha}$
(which was the whole point of defining $|\cdot|_{aut}$
on subsets of edges.)  Hence $|\alpha|_{tot}= |e(\alpha)|_{tot}$
also, as Krstic and Vogtmann show the corresponding formula
for $|\cdot|_{out}$.  From this, the analogues of Proposition
6.4  about Whitehead moves in \cite{[K-V]} are true
for the norms $\|\cdot\|_{aut}$ and
$\|\cdot\|_{tot}$.  That is, for an ideal edge $\alpha$ define
$D(\alpha)$ by
$$D(\alpha) = \{ a \in \alpha : stab(a)=stab(\alpha) \hbox{ and }
\bar a \not \in \bigcup G \alpha \}.$$
Then the Whitehead move $(G\alpha,Ga)$ is the result
of first blowing up $\alpha$ in $\Gamma$ to get $\sigma^{G\alpha}$ and
then collapsing $Ga$ in $\Gamma^{G\alpha}$ to get
$\sigma'$.  Proposition 6.4 of \cite{[K-V]} states that
$$\|\sigma'\|_{out} = \|\sigma\|_{out} + [G:stab(\alpha)]
(|\alpha|_{out} - |a|_{out}).$$
As mentioned before, this remains true if out-norms and
absolute values are replaced by aut- or tot-norms and
absolute values.

The value $[G:stab(\alpha)] (|a|_{out} - |\alpha|_{out})$
is called the {\em out-reductivity} of $(\alpha,a)$
and is denoted $red_{out}(\alpha,a)$.  Similar
notions of {\em aut-reductivity} and {\em tot-reductivity}
are defined as well.  A Whitehead move reduces the norm
iff the corresponding reductivity is greater than zero,
in which case the Whitehead move is called
{\em reductive}.
The {\em $x$-reductivity} of
an ideal edge $\alpha$ is the maximum over all
elements $a \in D(\alpha)$ of $red_x(\alpha,a)$,
where $x$ is out, aut, or tot.  It thus makes
sense to talk of an ideal edge $\alpha$ as being
{\em out-reductive}, etc.  The main reason we
are interested in the norm $\|\cdot\|_{tot}$
rather than just using the norm $\|\cdot\|_{aut}$ is the following
proposition.

\begin{prop} \label{t17} Let $\alpha \subseteq E_v$
be a tot-reductive ideal edge of a reduced
marked $G$-graph $\rho$.  Suppose $\alpha$ is invertible
(that is, $E_v - \alpha \not \subseteq G\alpha$
and $E_v - \alpha$ is an ideal edge.)  Then
$\alpha^{-1} = E_v - \alpha$ is tot-reductive.
\end{prop}

\PF Assume $v=*$, else the proof is trivial.  Say $(\alpha,a)$ is the
reductive ideal edge.  Since $stab(*)=G$ and $\alpha$ is
invertible, the analogue of Lemma 5.1 of \cite{[K-V]}
gives us that $stab(\alpha)=stab(a)=G$.  Say $\rho=[s,\Gamma].$
As before, let $\rho^{G\alpha}=[s^{G\alpha},\Gamma^{G\alpha}]$
be the result of blowing up the ideal edge $(\alpha,a)$.
Then let $\rho'=[s',\Gamma']$ be the result of collapsing
$a$ in $\Gamma^{G\alpha}$.  We know that
$\|\rho'\|_{tot} < \|\rho\|_{tot}$ as $(\alpha,a)$
is tot-reductive.

Assuming the claim below, it will be easy to complete the
proof as follows:  Since $\|\rho'\|_{out} \not = \|\rho\|_{out}$
and $\|\rho'\|_{tot} < \|\rho\|_{tot}$,
we must have $\|\rho'\|_{out} < \|\rho\|_{out}$.
Let $\rho''$ be the result of doing the Whitehead move
$(\alpha^{-1},a^{-1})$ to $\rho$.  Because
$red_{out}(\alpha,a)=red_{out}(\alpha^{-1},a^{-1})$
(see the comments in \cite{[K-V]} following the proof
of \S 6.4), it follows that
$\|\rho''\|_{out} = \|\rho'\|_{out}$.
So $\|\rho''\|_{out} < \|\rho\|_{out}$ and
hence $\|\rho''\|_{tot} < \|\rho\|_{tot}$.
Thus $\alpha^{-1}$ is reductive. \END

\begin{claim} \label{t18} $\|\rho'\|_{out} \not = \|\rho\|_{out}$.
\end{claim}

\PF Since $stab(a)=G$ and $\rho = [s,\Gamma]$ is reduced, the edge $a$
must both begin and end at $*$.  (The only
edges which end at $*$ but begin somewhere else
are the ones on which $G$ acts nontrivially, else
they are a $G$-invariant forest which we could collapse.)
Since $a \in D(\alpha)$, $a^{-1} \not \in \alpha$.
There is at least one edge in $\alpha - \{a\}$.
Enumerate the edges of $\alpha - \{a\}$ as
$b_0, \ldots, b_r$.  Similarly, there is at
least one edge in $\alpha^{-1} - \{a^{-1}\}$.
Enumerate these as $c_0, \ldots, c_s \in \alpha^{-1} - \{a^{-1}\}$.
We have three cases, which are not disjoint but
are exhaustive.
\begin{enumerate}
\item Some $b_i$ is a loop at $*$ and $b_i^{-1} \not \in \alpha$.
Let $w_k \in \mathcal{W}$ be an element that maps to the
loop $b_i$.  Then $(|a|_{out})_k = 0$ and
$(|\alpha|_{out})_k \geq 1$ since the loop $b_i$ is
sent to $b_ie(\alpha)$.

\item Some $b_i$ starts at another vertex $v_1 \not = *$.
Since $G$ acts nontrivially on $b_i$ and because $b_i$
must be elliptic (as it is clearly not bent hyperbolic),
there must be another $b_j \not = b_i$
also going from $*$ to $v_1$.
(For the definitions
of elliptic and bent hyperbolic see \S 4A in the paper
by Krstic and Vogtmann.)  There are two subcases:
\begin{itemize}
\item There is an edge $c_l$ in $\alpha^{-1} - \{a^{-1}\}$ that
begins and ends at $*$.  We can assume
$c_l^{-1} \in \alpha^{-1} - \{a^{-1}\}$ also,
else we are in case 1.  Choose a $w_k \in \mathcal{W}$ that
maps to the loop $b_ib_j^{-1}c_l$.
Now $(|a|_{out})_k = 0$ and
$(|\alpha|_{out})_k \geq 1$ since $w_k$ is
sent to $e(\alpha)^{-1}b_ib_j^{-1}e(\alpha)c_l$.
\item There is an edge $c_l$ in $\alpha^{-1} - \{a^{-1}\}$ that
begins at $v_2 \not = *$ and ends at $*$.
Because $G$ acts nontrivially on $c_l$ and $c_l$ is
elliptic, there is another
edge $c_m \not = c_l$ also going from $v_2$ to $*$.
Choose a $w_k \in \mathcal{W}$ that
maps to the loop $b_ib_j^{-1}c_lc_m^{-1}$.
Then $(|a|_{out})_k = 0$ but
$(|\alpha|_{out})_k \geq 1$ as $b_ib_j^{-1}c_lc_m^{-1}$ is
sent to $e(\alpha)^{-1}b_ib_j^{-1}e(\alpha)c_lc_m^{-1}$.
\end{itemize}

\item Some $b_i$ is a loop at $*$ and $b_i^{-1} \not \in \alpha$
also.  As in case $2.$ above, there are two subcases.
\begin{itemize}
\item Same as in case $2.$ above.  Choose a $w_k \in \mathcal{W}$
that maps to $b_ic_l$.
Then $(|a|_{out})_k = 0$ but
$(|\alpha|_{out})_k \geq 1$ as $b_ic_l$ is
sent to $e(\alpha)^{-1}b_ie(\alpha)c_l$.
\item Same as in case $2.$ above.  Choose a $w_k \in \mathcal{W}$
that maps to $b_ic_lc_m^{-1}$.
Hence $(|a|_{out})_k = 0$ and yet
$(|\alpha|_{out})_k \geq 1$ because $b_ic_lc_m^{-1}$ is
mapped to $e(\alpha)^{-1}b_ie(\alpha)c_lc_m^{-1}$.
\end{itemize}
\end{enumerate}

In each case we have $|a|_{out} \not = |\alpha|_{out}$;
therefore, $red_{out}(\alpha,a) \not = 0$
and $\|\rho'\|_{out} = \|\rho\|_{out}$. \END

Proposition $6.1$ of \cite{[K-V]}, states that
$$\matrix{\hfill |A\coprod B|_{out} = |A|_{out} + |B|_{out} - 2(A.B)_{out} \hfill \cr}$$
for disjoint subsets $A$ and $B$ of $E(\Gamma)$.
This also holds for $aut$-norms because it
is our definition of the absolute values $|\cdot|_{aut}$
for sets of edges and can be inductively shown
to be well-defined.
It is important that this
property holds for $aut$-norms because
it is used by
many of the later propositions in Krstic and
Vogtmann (e.g., Proposition $6.2$ of \cite{[K-V]}
which will correspond to our Proposition \ref{t19}.)

Proposition $6.2$ of \cite{[K-V]} states that:

\begin{prop}[Krstic-Vogtmann] \label{tr38}
Let $K$ be a subgroup of $G$, let $A$ be a $K$-invariant
subset of $E(\Gamma)$, and let $e$ be an edge of $\Gamma$
with $stab(e)$ contained in $K$.  Then
$$((Ke).A)_{out} = [K:stab(e)](e.A)_{out}.$$
\end{prop}

We now show Proposition $6.2$
of \cite{[K-V]} also holds
for the aut-norm,
which will be useful in some
combinatorial lemmas later in this chapter.
Once we show that the analogue of Proposition \ref{tr38}
is true for the $aut$-norm,
it
will be true for both the out- and aut-norms
on a component-by-component basis.  In other words,
the equality stated in the proposition
is true for each component of
$\Z^\mathcal{W}$ or $\Z^{F_n}$ and does not use the
total (lexicographic) order on those sets.  Hence
it is automatically true for the tot-norm, as the
tot-norm is just the product of the out-norm
and the aut-norm.  We will be able to use the
same approach (that of just showing something to be
true for the aut-norm)
in some lemmas later on in this chapter.

\begin{prop} \label{t19}
Let $K$ be a subgroup of $G$,
$A$ be a $K$-invariant subset of $E(\Gamma)$,
and $e$ be an edge of $\Gamma$ with $stab(e)$
contained in $K$.  Then
$$((Ke).A)_{aut}=[K:stab(e)](e.A)_{aut}.$$
\end{prop}

\PF To simplify the notation in the proof below,
we write (just for this proof)
$\|\cdot\|$ for $\|\cdot\|_{aut}$,
$|\cdot|$ for $|\cdot|_{aut}$,
reductive for aut-reductive, etc.
 
Examine $((Ke).A)_i$.  It is the number of
times one of the strings $(ke)a^{-1}$ or
$a(ke)^{-1}$ appears in one of the
$x\alpha_i$, for all $k \in K$, $a \in A$, and
$x \in G$.

Now $stab(e) \subseteq K$ and we can write
$$\matrix{\hfill K = stab(e) \coprod k_2stab(e) \coprod \ldots
\coprod k_{[K:stab(e)]}stab(e) \hfill \cr}$$
using coset representatives $k_i$.
Further note that the number of times one
of the strings $ea^{-1}$ or $ae^{-1}$
appears in one of the strings $x\alpha_i$
for $a \in A$, $x \in G$ is exactly the same
as the number of times one of the strings
$k_iea^{-1}$ or $a(k_ie)^{-1}$ appears in
one of the $x\alpha_i$
for $a \in A$, $x \in G$.  This is because
each $k_i$ is in $G$ and $A$ is
$K$-invariant so if $ea^{-1}$ is in $x\alpha_i$
then $(k_ie)(k_ia)^{-1}$ is in
$(k_ix)\alpha_i$.
So $((Ke).A)_i = [K:stab(e)] (e.A)_i$. \END

Because of Proposition \ref{t17}, tot-reductivity
will be the most useful of the three types of reductivity
(out, aut, and tot) for us.
From now on when we say that $\rho$ is
{\em reductive}, this is just shorthand for saying $\rho$ is
tot-reductive.

\begin{prop} \label{t20} The set of basepointed marked
$G$-graphs is well-ordered by the tot-norm.
\end{prop}

\PF Note that Proposition $6.3$ of \cite{[K-V]} establishes the
corresponding result for the out-norm.  (It probably also
holds for the aut-norm, but as we do not need this here we
make no attempt to prove it.)

Let $\mathcal{A}$ be a nonempty collection of basepointed
marked $G$-graphs.  We must find a least element of $\mathcal{A}$.
Let $[\mathcal{A}]$ be the set of equivalence classes of
marked $G$-graphs in $\mathcal{A}$ obtained by forgetting
the basepoint $*$.  Since from
Proposition $6.3$ of \cite{[K-V]} the out-norm
well orders marked $G$-graphs, $[\mathcal{A}]$ has a
least element $U \subseteq \mathcal{A}$.

Say $\sigma = [s,\Gamma]$ is the (non-basepointed) marked
$G$-graph representing this $U$.  The marked graph $\sigma$
corresponds to an action of $F_n$ on the tree
$\tilde \Gamma = \Lambda$.
From \cite{[C-V]} this action is free, {\em minimal} (there
are no invariant proper subtrees), and not abelian
(an action is {\em abelian} iff every element of the
commutator $[F_n,F_n]$ has length $0$.)  It also follows that
the action is without inversions.  (Suppose $1 \not = a \in F_n$
acts as an inversion.  Then $a^2$ has a fixed point.
Because the action of
$F_n$ on $\Lambda$ is free, this means $a^2=1$.  In the
free group $F_n$, however, you can not have $a \not = 1$
with $a^2=1$.)

The action has an associated length function $l$ on $F_n$.
Since the action is non-abelian, the length function is
non-abelian.  (See Alperin and Bass in
\cite{[A-B]} for the definition of an
abelian length function.  Further see Theorem $7.6(a)$ of
\cite{[A-B]} for a theorem which gives the above implication.)
By Theorem $7.4$ of \cite{[A-B]}, there exist $\alpha_n,
\alpha_m \in F_n$ such that the characteristic
subtrees $A_{\alpha_n}$ and $A_{\alpha_m}$ are disjoint.
(Recall that $F_n = \{\alpha_1, \alpha_2, \ldots\}$ as
an ordered set.)  Without loss of generality,
we may assume
$n < m$.

Since the action of $F_n$ is free, both $\alpha_n$ and $\alpha_m$
are hyperbolic and $A_{\alpha_n}, A_{\alpha_m}$
are linear subtrees.  Recall that we wish to find the least
element of $U$ in the tot-norm.  Following the
proof of Proposition $6.3$ in \cite{[K-V]}, we set
$U_0 = U$ and define $U_i$ inductively for $i \geq 1$.
Let $\gamma_i = min\{(\|\delta\|_{aut})_i : \delta \in U_{i-1}\}.$
Next define $U_i$ to be the subset of $U_{i-1}$ consisting of
$\delta$ with $(\|\delta\|_{aut})_i = \gamma_i$.
Our goal is to show that $U_m$ contains only
finitely many elements, because then we could
just compare these finitely many elements to 
get the least one.

In order to make explicit a minor technical point,
a couple of definitions are appropriate here.
Say we have an action of the free group $F_n$ on a
basepointed tree $T$.  Then the tree $T'$ obtained by
{\em forgetting the basepoint of $T$} is defined as
follows.  First, if the basepoint $*$ of $T$
had valence greater than or equal to $3$, then $T'$ is
obtained by replacing the vertex $*$ of $T$ with
a normal (i.e., not distinguished
and not a basepoint) vertex.
On the other hand, if $*$ of $T$ had valence $2$,
then we form $T'$ by removing from $T$ all of the vertices
in the $F_n$-orbit of $*$.  Hence, in both cases,
$T'$ is a tree which does not have a basepoint $*$
and all of whose vertices have valence at least $3$.

Next, we define an inverse of the above operation.
That is, say that we have an action of the free
group $F_n$ on a (non-basepointed) tree $T'$.
Choose a point $v$ on $T'$, where $v$ can either be a
vertex of $T'$ or the midpoint of an edge of $T'$.
Then the tree $T$ obtained by
{\em placing a basepoint on $T'$ at $v$}
is defined as follows.  First, if $v$ is a vertex of $T'$,
then $T$ is just the tree $T'$ with the vertex $v$
replaced by the basepoint $*$.
On the other hand, if $v$ is the midpoint of an
edge $e$ of $T'$, then define $T$ to be the result of
adding a basepoint $*$ to $T'$ at $v$ and
subdividing all of the edges in the $F_n$ orbit of $e$.
Hence, in both cases, we have an induced action of
$F_n$ on a simplicial, basepointed tree $T$.


Each element of $U$ corresponds to an action of $F_n$
on a basepointed tree.  In each case, if we forget the
basepoint then the tree is just $\Lambda$.  For each of these
actions corresponding to elements of $U$, we can obtain a Lyndon
length function on $F_n$ by seeing how far each element
moves the basepoint.  This map from $U$ to
Lyndon length functions is injective (see \cite{[H-V]},
\cite{[A-B]}.)  Note that in each case, the action of $F_n$
on the underlying non-basepointed tree $\Lambda$ is the
same.  We are only varying where we place the basepoint of $\Lambda$
and seeing how far elements of $F_n$ move this
basepoint.

The elements of $U_1$ are those where the basepoint is located
closest to the linear subtree $A_{\alpha_1} \subset \Lambda$.
Hence $U_1$ could be infinite.  (For example, if
one element of $U_1$ has its basepoint on $A_{\alpha_1}$,
then $U_1$ could be all of the vertices on $A_{\alpha_1}$.)
By the $m$th stage, however, we have already
considered both $\alpha_n$ and $\alpha_m$.  Because
$A_{\alpha_n}$ and $A_{\alpha_m}$ are disjoint, there
are only finitely many elements in $U_m$,
as we show in the next two paragraphs.

Let $B$ be the bridge joining $A_{\alpha_n}$
and $A_{\alpha_m}$.  Elements of $U_m$ have basepoints
that are both closest first to $A_{\alpha_n}$
and then to $A_{\alpha_m}$.  So they must be both
close to $A_{\alpha_n}$ and close to the endpoint
of the bridge $B$ on $A_{\alpha_m}$.
Since only finitely many basepoints can be
closest to the endpoint of this bridge, we
are done and $U_m$ is finite.

More formally, it suffices to show that there are only
finitely many points at fixed distances $d_1$ and
$d_2$ from $A_{\alpha_n}$ and $A_{\alpha_m}$, respectively.
If $d(x,A_{\alpha_n})=d_1$ and $d(x,A_{\alpha_m})=d_2$, then
choose paths $p_1$ and $p_2$ of lengths $d_1$ and $d_2$
from $x$ to $q_1 \in A_{\alpha_n}$ and
$q_2 \in A_{\alpha_m}$, respectively.  The union of these
two paths $p_1$ and $p_2$ contains the bridge $B$.
Consequently, $d(x,B) \leq d_1 + d_2$.  Since the tree
is locally finite and $B$ is finite, $x$ is one of a finite number
of vertices. \END

A few definitions are in order at this point.  Basically,
we are trying to find the appropriate parallels of
definitions in \cite{[K-V]}.  Fix a reduced marked
$G$-graph $\rho = [s,\Gamma]$.  Let $(\mu,m)$ be
a maximally reductive ideal pair of $\rho$.  That is,
$\mu$ is the maximally reductive ideal edge in $\rho$
and $m \in D(\mu)$ is an edge in $\mu$ which
allows the Whitehead move $(\mu,m)$ to realize this
maximum.

Let $\alpha \subset E_u$ and $\beta \subset E_v$ be ideal
edges of $\rho$.  Then the ideal edge orbits $G\alpha$
and $G\beta$ are {\em compatible} if one of the
following holds:
\begin{enumerate}
\item $G\alpha \subseteq G\beta.$
\item $G\beta \subseteq G\alpha.$
\item $G\alpha \cap G\beta = \void$ and $\alpha \not = \beta^{-1}.$
\item $G\alpha \cap G\beta = \void$ and $u=v=*.$
\end{enumerate}
The ideal edge orbits $G\alpha$
and $G\beta$ are {\em pre-compatible} if one of the
following holds:
\begin{enumerate}
\item They are compatible.
\item $\alpha$ is invertible and $\alpha^{-1} \subseteq \beta$.
\item $\beta$ is invertible and $\beta^{-1} \subseteq \alpha$.
\end{enumerate}
Note that $2.$ and $3.$ above would be equivalent if we did not
need to consider ideal edges of the form $\gamma = E_* - \{c^{-1}\}$
which have $stab(\gamma)=G$ but are not invertible.

An {\em oriented ideal forest} is a collection of
pairwise compatible ideal edge orbits.  These can be blown up
to obtain marked graphs in the star in $L_G$ of $\rho$.  The
correspondence is not unique, however, as two different
oriented ideal forests can be blown up to yield the same
marked graph.  This problem is solved by defining ideal
forests.  There is a poset isomorphism between the poset
of ideal forests and the star of $\rho$ in $L_G$.

An {\em ideal forest} is a collection $\Phi = \Phi_1 \coprod \Phi_2$
where $\Phi_1$ are the edges at $*$ and $\Phi_2$ are the edges not at $*$,
such that
\begin{enumerate}
\item The elements of $\Phi_2$ are pairwise pre-compatible
and $\Phi_2$ contains the inverse of each of its invertible
edge orbits; and
\item The elements of $\Phi_1$ are pairwise compatible.
\end{enumerate}

With respect to a particular reduced marked $G$-graph $\rho$
and maximally reductive ideal edge $(\mu,m)$, the following definitions
will be used frequently in the next chapter (which contains the core
proof of the contractibility of $L_G$.)
\begin{itemize}
\item $\mathcal{R} = \{\hbox{reductive ideal edges}\}$.
\item If $\mathcal{C}$ is a set of ideal edges, then let
$\mathcal{C}^\pm$ denote the set obtained by adjoining
to $\mathcal{C}$ the inverses of its invertible
elements that are not at the basepoint.
\item Let $S(\mathcal{C})$ be the subcomplex of
the star $st(\rho)$ spanned by ideal forests of
$\rho$, all of whose edges are in $\mathcal{C}$.
Note:  The empty forest should \underline{not}
be taken to be in $S(\mathcal{C})$.
\item $\mathcal{C}_0 = \{\alpha \in \mathcal{R} :
\alpha \hbox{ is compatible with } \mu\}$.
\item $\mathcal{C}_0^{'} = \mathcal{C}_0 \cup \{\alpha \in \mathcal{R} :
\hbox{ if } \alpha \subset E_v \hbox{ then }
stab(\alpha)=stab(v)\}$ (cf. Lemma $5.1$ of \cite{[K-V]}.)
\item $\mathcal{C}_1 = \mathcal{C}_0^{'} \cup
\{\alpha \in \mathcal{R} : m \in G\alpha \hbox{ and }
N(G\alpha,G\mu)=1\}$.
\end{itemize}
The definition of the crossing number $N(G\alpha,G\mu)$
comes from \S 7 of \cite{[K-V]} where it and other
combinatorial notions are defined.  For the
reader's convenience, we briefly state their
definitions again here.  Say
$\alpha$ and $\beta$ are two ideal edges at some
vertex $v$, with stabilizers $P$ and $Q$, respectively,
of indices $p$ and $q$ in $G$.
Choose double coset representatives
$x_1, \ldots, x_k$ of $P\backslash G/Q$. The intersection
$\delta = \alpha \cap G\beta$ breaks up as
a disjoint union
$$\matrix{\hfill \delta = \gamma_1 \coprod \ldots \coprod \gamma_k \hfill \cr}$$
with each $\gamma_i = \alpha \cap Px_i\beta$.
The $\gamma_i$ are called the {\em intersection
components} of $\alpha$ with $\beta$
and the number $N(G\alpha,G\beta)$ of nonempty
intersection components is
called the {\em crossing number}.
If $N(G\alpha,G\beta)=1$ then $G\alpha$ and $G\beta$
are said to {\em cross simply}.

The following two lemmas are stated for the out-norm
by Krstic and Vogtmann.  We will show them for
the aut-norm.
The proofs will be routine,
although they are not the same as the proofs given in
\cite{[K-V]}. This is because their proofs use the
fact that $|A|_{out} = (A.E(\Gamma)-A)_{out}$, which
is no longer true with the new norms. 
As with Proposition \ref{t19}, the
lemmas are true for both the out- and aut- norms
on a component-by-component basis.  That is,
the inequalities stated in the lemmas
are true for each component of
$\Z^\mathcal{W}$ or $\Z^{F_n}$ and do not use the
total (lexicographic) order on those sets.  Hence
it suffices to show them for the aut-norm, as the
tot-norm is just the products of the out-norm
and the aut-norm.

\begin{lemma} \label{t21} 
Suppose $G\alpha$ and $G\beta$ cross simply,
with $P \leq Q$, then
$$p|\alpha \cap \beta|_{aut} + q|\beta \cup Q\alpha|_{aut}
\leq p|\alpha|_{aut} + q|\beta|_{aut}.$$
\end{lemma}

\PF To simplify the notation in the proof below,
we write (just for this proof)
$\|\cdot\|$ for $\|\cdot\|_{aut}$,
$|\cdot|$ for $|\cdot|_{aut}$,
reductive for aut-reductive, etc.

Let $[Q:P]=n$. Then $p=nq$.  Dividing by $q$,
we see that we want to show that
$$n|\alpha \cap \beta| + |\beta \cup Q\alpha|
\leq n|\alpha| + |\beta|.$$

Let $q_1, \ldots, q_n$ be a set of coset
representatives for $P$ in $Q$.
Let $\delta = \alpha \cap \beta$,
$A=\alpha - \delta$, and
$B=\beta - Q\delta$.

Then the left hand side of what we want to show is
$$\matrix{
\hfill n|\delta| + |B \coprod Q\alpha| &= n|\delta| + |B| + |Q\alpha| - 2 B.Q\alpha \hfill \cr
\hfill &= n|\delta| + |B| + |Q\delta \coprod QA| - 2 B.Q\alpha \hfill \cr
\hfill &= n|\delta| + |B| + |Q\delta| + |QA| - 2 Q\delta.QA - 2 B.Q\alpha. \cr}$$
Similarly, the right hand side of what we want to show is
$$\matrix{\hfill n|\delta \coprod A| + |B \coprod Q\delta| =
n|\delta| + n|A| - 2n \delta.A + |B| + |Q\delta| - 2 B.Q\delta. \hfill \cr}$$
So we have reduced the problem to showing that
$$|QA| - 2 Q\delta.QA - 2 B.Q\alpha \leq
n|A| - 2n \delta.A - 2 B.Q\delta.$$

Note that $-2 B.Q\alpha \leq -2 B.Q\delta$
as $Q\delta \subseteq Q\alpha$.
Also note that
$$\matrix{
\hfill Q\delta.QA &= Q\delta.(q_1A \coprod \ldots \coprod q_nA) \hfill \cr
\hfill &= Q\delta.q_1A + \ldots + Q\delta.q_nA \hfill \cr
\hfill &= Q\delta.A + \ldots + Q\delta.A \hfill \cr
\hfill &= n Q\delta.A \hfill \cr
\hfill &\geq n \delta.A. \hfill \cr}$$
Hence $-2 Q\delta.QA \leq -2n \delta.A$.
In addition, we could use induction to show
$$\matrix{
\hfill |QA| &= |q_1A \coprod \ldots \coprod q_nA| \hfill \cr
\hfill &\leq |q_1A| + \ldots + |q_nA| \hfill \cr
\hfill &= |A| + \ldots + |A| \hfill \cr
\hfill &= n|A|. \hfill \cr}$$
This establishes the lemma. \END

\begin{lemma} \label{t22}
Suppose $G\alpha$ and $G\beta$ {\em cross}
(i.e., $N(G\alpha,G\beta) \not = 0$).
Just as $\delta$ breaks up into
intersection components of $\alpha$ with $\beta$,
let $\delta' = \beta \cap G\alpha$ give
the analogous disjoint components
$$\matrix{\hfill \delta' =
\gamma_1' \coprod \ldots \coprod \gamma_k' \hfill \cr}$$
with $\gamma_i' = \beta \cap Qx_i^{-1}\alpha$.
Then for all $i$,
$$p|\alpha-\gamma_i|_{aut} + q|\beta-\gamma_i'|_{aut} \leq
p|\alpha|_{aut} + q|\beta|_{aut}.$$
\end{lemma}

\PF To simplify the notation in the proof below,
we write (just for this proof)
$\|\cdot\|$ for $\|\cdot\|_{aut}$,
$|\cdot|$ for $|\cdot|_{aut}$,
reductive for aut-reductive, etc.

Let $A = \alpha - \gamma_i$ and $B = \beta - \gamma_i'$.
We must show that
$$p|\gamma_i| + q|\gamma_i'| \geq 2\gamma_i.A + 2 \gamma_i'.B.$$
Note that $G\gamma_i=G\gamma_i'$.
Choose coset representatives $y_1, \ldots, y_p$ for
$P$ in $G$ and $z_1, \ldots, z_q$ for $Q$ in $G$.
Accordingly, we can express the left hand
side of what we want to show as
$$\matrix{
\hfill p|\gamma_i| + q|\gamma_i'| &=
\sum_{n=1}^p |\gamma_i| + \sum_{m=1}^q |\gamma_i'| \hfill \cr
\hfill &= \sum_{n=1}^p |y_n\gamma_i| + \sum_{m=1}^q |z_m\gamma_i'| \hfill \cr
\hfill &\geq |\coprod_{n=1}^p y_n\gamma_i| + |\coprod_{m=1}^q z_m\gamma_i'| \hfill \cr
\hfill &= |G\gamma_i| + |G\gamma_i'| \hfill \cr
\hfill &= 2 |G\gamma_i| = 2 |G\gamma_i'|. \hfill \cr}$$
Moreover, the right hand side of what we want to show is
$$\matrix{
\hfill 2\gamma_i.A + 2\gamma_i'.B &\leq
2G\gamma_i.A + 2G\gamma_i'.B \hfill \cr
\hfill &= 2G\gamma_i.(A\coprod B) \hfill \cr
\hfill &\leq 2G\gamma_i.(E(\Gamma) - G\gamma_i). \hfill \cr}$$
So to prove the lemma it suffices to show
$$|G\gamma_i| \geq G\gamma_i.(E(\Gamma) - G\gamma_i).$$
This can be shown by induction on $|G\gamma_i|$.
In other words, say $|C| \geq C.(E(\Gamma)-C)$ and assume
$y \not \in C$.  Observe that
$$\matrix{
\hfill |C \coprod \{y\}| &= |C| + |y| - 2C.y \hfill \cr
\hfill &\geq C.(E(\Gamma) - C) + |y| - 2C.y \hfill \cr
\hfill &= (C.(E(\Gamma)-C) - C.y) + (|y| -  C.y) \hfill \cr
\hfill &\geq C.(E(\Gamma) - (C \coprod \{y\}))  + 0 \hfill \cr}$$
so that the inductive step is completed. \END

Next we show that the Pushing and
Shrinking Lemmas of Krstic and Vogtmann
hold in the context of aut-norms and
absolute values.  Unlike the proofs of the
previous two lemma, the proofs for the
next two follow exactly the same lines
as the original proofs by Krstic and Vogtmann
for out-norms and absolute values.
Basically, the only way that the new proofs will differ
from the old ones is that occasionally we
will have to verify the new
cardinality conditions for ideal
edges $\alpha_0 \subseteq E_v$, namely that:
\begin{itemize}
\item If $v = *$ then $card(\alpha_0) \geq 2$ and
$card(E_v - \alpha_0) \geq 1$.
\item If $v \not = *$ then $card(\alpha_0) \geq 2$ and
$card(E_v - \alpha_0) \geq 2$.
\end{itemize}
As before, it is easily seen from the proofs of
the lemmas that since they hold for both the out-
and aut-norms and absolute values, they also
hold for the tot-norms and absolute values.

\begin{lemma}[Pushing Lemma] \label{t23} 
Let $(\mu,m)$ be a maximally aut-reductive ideal edge
of a reduced basepointed marked $G$-graph
with $m \in D(\mu)$.  Let $(\alpha,a)$ be
an aut-reductive ideal edge
containing $m$ which simply crosses $\mu$,
and set $P=stab(\alpha)$.
Then either both $\mu - \alpha$ and $\alpha - \mu$
are aut-reductive or both $\alpha \cup P\mu$ and
$\alpha \cap \mu$ are aut-reductive.
\end{lemma}

\PF To simplify the notation in the proof below,
we write (just for this proof)
$\|\cdot\|$ for $\|\cdot\|_{aut}$,
$|\cdot|$ for $|\cdot|_{aut}$,
reductive for aut-reductive, etc.

Note that since $m \in \alpha$, $stab(\mu) \leq P$.
As in \cite{[K-V]}, there are four cases depending
upon where $a^{-1}$ and $m^{-1}$ are located.
Since this follows the proof by Krstic and Vogtmann
so closely, the only real detail will be put into the
first case.

\medbreak

\noindent {\em Case 1.} $a^{-1} \not \in G\mu$.
From Lemma \ref{t21},
$$[G:stab(\mu)]|\alpha \cap \mu| + [G:stab(\alpha)]|\alpha \cup P\mu|
\leq [G:stab(\mu)]|\mu| + [G:stab(\alpha)]|\alpha|.$$
Consequently,
$$[G:stab(\mu)](|m|-|\alpha \cap \mu|) +
[G:stab(\alpha)](|a|-|\alpha \cup P\mu|)$$
is greater than or equal to 
$$[G:stab(\mu)](|m|-|\mu|) + [G:stab(\alpha)](|a|-|\alpha|).$$
In other words,
$$red(\alpha \cap \mu,m) + red(\alpha \cup P\mu)
\geq red(\mu,m) + red(\alpha,a).$$

Since $(\mu,m)$ is maximally reductive and $(\alpha,a)$ is
reductive, both of $(\alpha \cup P\mu,a)$
and $(\alpha \cap \mu, m)$ are reductive. As mentioned
above in the discussion preceeding this lemma, we must
verify the cardinality conditions on these two prospective
ideal edges.

First we deal with $(\alpha \cup P\mu,a)$. The
edge $a$ is either bent hyperbolic or elliptic (see
Corollary $4.5$ of \cite{[K-V]}.)  Assume it is bent hyperbolic.
Then as in \cite{[K-V]} we can choose $x \in G$ such that
$xa^{-1} \in E_v-(G\alpha \cup G\mu).$  If $v \not = *$
and $xa^{-1}$ is the only edge in $E_v-(\alpha \cup P\mu)$
then
$$|\alpha \cup P\mu| = |xa^{-1}| = |a^{-1}| = |a|,$$
where the first equality holds because $v \not = *$,
the second is by the $G$-invariance of $|\cdot|$,
and the third follows from our definition of $|\cdot|$
for edges.

In more detail, the first equality
$|\alpha \cup P\mu| = |xa^{-1}|$ holds since
$$\matrix{\hfill E_v = (\alpha \cup P\mu) \coprod \{xa^{-1}\} \hfill \cr}$$
and $v \not = *$.  For a particular coordinate $i$,
both $|\alpha \cup P\mu|_i$ and $|xa^{-1}|_i$ are
measuring the number of times one of the paths $y\alpha_i$
enters $v$ via $\alpha \cup P\mu$ and
leaves via the reverse of $xa^{-1}$ (i.e., $xa$)
\underline{or} enters $v$ via $xa^{-1}$ and leaves it
via the reverse of something in $\alpha \cup P\mu$.
There would be problems if $v=*$ since the
above paths could then enter $v$
and not have to leave it again.

But $|\alpha \cup P\mu| = |a|$ contradicts the fact that
$(\alpha \cup P\mu,a)$ is reductive because
$$[G:stab(\alpha)] (|a| - |\alpha \cup P\mu|) > 0.$$
So if $v \not = *$
then $xa^{-1}$ is not the only edge in $E_v-(\alpha \cup P\mu)$.

For the next possibility, that $a$ is elliptic,
the proof by Krstic and Vogtmann can be used verbatim.

Second we deal with $(\alpha \cap \mu,m)$.  The
set $\alpha \cap \mu$ must contain more than two
edges because it is reductive:
$$[G:stab(\mu)] (|m| - |\alpha \cap \mu|) > 0.$$
The condition on the cardinality of
$E_v - (\alpha \cap \mu)$ is easily
satisfied because $\alpha$ is an ideal
edge and so satisfies the corresponding
condition with $E_v - \alpha$.

\medskip

\noindent {\em Case 2.} $a^{-1} \in G\mu$ and $m^{-1} \in G\alpha$.

Krstic and Vogtmann show that both $(\alpha-\mu,ym^{-1})$
and $(\mu-\alpha,xa^{-1})$ are reductive.  The
cardinality conditions on the edges for $\alpha-\mu$ are
easily satisfied.  First, $\alpha-\mu$ has more than one edge
because it is a reductive set.  Second,
$E_v - (\alpha-\mu)$ satisfies its cardinality
condition because $E_v - \alpha$ does.  Analogous
arguments also show the cardinality conditions for
$\mu - \alpha$.

\medskip

\noindent {\em Case 3.} $a^{-1} \in G\mu$,
$m^{-1}  \not \in G\alpha$, and $a \in \mu$.

Both $(\alpha \cup \mu,m)$ and $(\alpha \cap \mu,a)$
are reductive.  See case $1.$ for how the
cardinality checks should go.

\medskip

\noindent {\em Case 4.} $a^{-1} \in G\mu$,
$m^{-1}  \not \in G\alpha$, and $a \not \in \mu$.

Both $(\alpha \cup \mu,m)$ and $(\alpha \cap \mu,m)$
are reductive.  See case $1.$ for how the
cardinality checks should go. \END

\begin{lemma}[Shrinking Lemma] \label{t24}
Let $(\mu,m)$ be as before. Let $\alpha$ be
an ideal edge with $N(G\alpha,G\mu) \not = 0$.
Let $\gamma_{i_1}, \ldots, \gamma_{i_k}$ be
the intersection components of $\alpha$ with $\mu$ which
contain no translate of $m$ and let
$\beta = \alpha - \bigcup \gamma_{i_j}$.  Then
$\beta$ or one of the sets $\gamma_{i_j}$ is an
aut-reductive ideal edge.
\end{lemma}

\PF See verbatim the proof by Krstic and Vogtmann.
If $\alpha_0$ is one of the above
sets, we know it is aut-reductive, and
we want to show it is an ideal
edge, then the cardinality checks
are easy.  The set $\alpha_0$ contains more than
one edge because it is aut-reductive.
Moreover, the cardinality checks on
$E_v - \alpha_0$ follow from similar ones
on $E_v - \alpha$, because
$\alpha_0 \subset \alpha$ for each possibility of
$\alpha_0$. \END

The following proposition will also be useful in the
next chapter.  As before, ideal edges are in the
reduced, basepointed marked $G$-graph $\rho$
and $(\mu,m)$ is a maximally reductive ideal edge.

\begin{prop} \label{t25}  There is at most one
reductive ideal edge $(\gamma,c)$ at $*$ with \newline
$stab(\gamma)=G$ but where $\gamma$ is not
invertible.  The Whitehead move $(\gamma,c)$ is just
conjugation by $c$, and $\|\gamma\|_{out}=0$.
If $\gamma$ is not compatible with $\mu$, then
$c=m^{-1}$ and $\mu$ is invertible.
\end{prop}

\PF Since $stab(\gamma)=G$ and yet $\gamma$ is not invertible,
$E_* - \gamma$ must contain just one element.
As $c^{-1} \in E_* - \gamma$, this means that
$E_* - \gamma = \{c^{-1}\}$.  So the Whitehead move
$(\gamma,c)$ consists of first blowing up $\gamma$ and then collapsing
$c$.

\bigskip

\input{fig15.pic}
\begin{figure}[here]
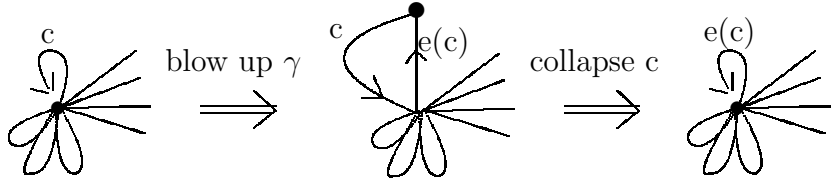

\caption{\label{fig15} The Whitehead move $(\gamma,c)$.}
\end{figure}

The Whitehead move has no effect on the out-norm.  The effect
on the aut-norm can be calculated as follows.
Recall that $F_n = \{\alpha_1, \alpha_2, \ldots \}$.  The
Whitehead move $(\gamma,c)$ conjugates each $\alpha_i$ by
$c$; i.e., each $\alpha_i \mapsto c^{-1}\alpha_ic$ (or
$e(c)^{-1}\alpha_ie(c)$ more accurately, but in the
final graph we can just relabel $e(c)$ as $c$.)

Since $\gamma$ is reductive, it must decrease the length of
some $\alpha_i$.  Say $\alpha_n$ is the first one whose length
is decreased.  Then the $\alpha_1, \ldots, \alpha_{n-1}$ are
of the form
$$(c\leftrightsquigarrow x) \hbox{ or } (y\leftrightsquigarrow c^{-1})$$
where $x \not = c^{-1}$, $y \not = c$ and the
$\leftrightsquigarrow$ are just any words in the edges.
In addition, $\alpha_n$ is
of the form
$$c\leftrightsquigarrow c^{-1}.$$
Thus for any other $\gamma'=E_* - \{ d^{-1} \}$, $d \not = c$,
$\gamma'$ will increase the length of one of the
$\alpha_1, \ldots, \alpha_n$
and it will not be reductive.  So $(\gamma,c)$ is the only
reductive edge at $*$ with $stab(\gamma)=G$ but where
$\gamma$ is not invertible.

Further suppose that $\gamma$ is not compatible with $\mu$.
If $c^{-1} \not \in \mu$, then $\mu \subseteq \gamma$ and they
are compatible.  So $c^{-1} \in \mu$.
As $stab(c^{-1}) = G$, $stab(\mu)=G$.
Since  $\mu$ is reductive and not equal to $\gamma$,
$\mu$ is invertible.  By way of contradiction,
assume $m \in \gamma$.  Then $m \not = c$ else $\gamma$
and $\mu$ are compatible.
We apply the Pushing Lemma to $\gamma$ and $\mu$.
Case $2$ is the relevant case and so
$(\mu-\gamma,a^{-1})$ is reductive.  This
contradicts the fact that $\mu - \gamma = \{a^{-1}\}$
has just one edge in it. So $m \not \in \gamma$
and hence $c=m^{-1}$. \END
 
\chapter{The contractibility lemmas} \label{c12}

\begin{thm} \label{tr13}
The fixed point subcomplex $X_n^G$ is contractible.
\end{thm}

\PF The space $X_n^G$ deformation retracts to $L_G$.
Following the proof of Theorem $8.1$ by Krstic and Vogtmann
in \cite{[K-V]}, we show that the complex $L_G$ is contractible
by setting
$$L_{<\rho} = \bigcup_{\|\rho'\|_{tot} < \|\rho\|_{tot}} st(\rho')$$
and letting
$$S_\rho= st(\rho) \cap L_{<\rho}.$$
As in \cite{[K-V]}, we show that $S_\rho$ is
contractible when it is non-empty, so that a
transfinite induction argument then yields that
for all $\rho$, all of the components of $L_{<\rho}$
are contractible.  Krstic's work in \cite{[K]} shows that
any two reduced graphs in $L_G$ can be connected by Whitehead moves,
so that $L_G$ is connected.  Thus $L_G$ is contractible
if we can perform the above transfinite induction.

As in \cite{[K-V]}, the first step is to
deformation retract $S_{<\rho}$ to to $S(\mathcal{R})$
by the Poset Lemma (stated in \cite{[K-V]}, deriving from
Quillen in \cite{[Q]}.) We can do this
for the case of $Aut(F_n)$ rather than $Out(F_n)$
without any
significant modifications of the arguments in the previous case.
This is because the Factorization Lemma
and Proposition $6.5$ of $\cite{[K-V]}$ let
us identify $S_{<\rho}$ with the poset of ideal forests
which contain a reductive ideal edge (where we must,
of course, use the newly modified definition of an
ideal forest.)  (The Factorization Lemma gives
a certain isomorphism between forest that does {\em not}
preserve basepoints, but the fact that basepoints are
not preserved is not relevant to Proposition $6.5$.)

After contracting $S_{<\rho}$ to $S(\mathcal{R})$, Krstic
and Vogtmann then use a series of lemmas to
deformation retract from $S(\mathcal{R})$ to $S(\mathcal{C}_1)$,
from there to $S(\mathcal{C}_0)$, and finally
to a point.  We more or less follow this, except there
is an additional
intermediate step where we deformation retract from
$S(\mathcal{C}_1)$ to $S(\mathcal{C}_0^{'})$
and from there to $S(\mathcal{C}_0)$.

The rest of this chapter will be devoted to proving the
aforementioned series of lemmas which show that
$S(\mathcal{R})$ deformation retracts to a point. \END

We assume that the maximally reductive ideal edge
$(\mu,m)$ is at the basepoint in all that follows, else the
arguments of Krstic and Vogtmann directly give the
contractibility of $S(\mathcal{R})$.
Moreover, if $(\mu,m) = (\gamma,c)$ where
$\gamma = E_* - \{c^{-1}\}$, then $\mu$ is not
out-reductive at all.  Since $\mu$ is maximally
reductive, Proposition \ref{t25} implies that $\mu$ is the
only reductive ideal edge.  So in this case
$\mathcal{R} = \mathcal{C}_0 = \{\mu\}$ and
$S(\mathcal{R})$ is contractible.  Assume $\mu \not = \gamma$
from now on.

Note the
slight difference in our definition of $\mathcal{C}_1$
from that of \cite{[K-V]},
where here it is phrased to include $\alpha \subset E_v$
which have $stab(\alpha) = stab(v)$, rather than just
invertible $\alpha$.  In other words, from Proposition \ref{t25},
there is at most one reductive ideal edge $(\gamma,c)$
at the basepoint which has $stab(\gamma)=G$ and yet is
not invertible.  This $\gamma$ would be in both $\mathcal{C}_1$
and $\mathcal{C}_0^{'}$.

The next lemma (unlike the ones which follow it) is
essentially the corresponding lemma in \cite{[K-V]} with
minimal modifications.  We repeat their
arguments here for the sake of convenience.

\begin{lemma} \label{t26} The complex $S(\mathcal{R})$ deformation
retracts onto $S(\mathcal{C}_1)$.
\end{lemma}

\PF Let $\mathcal{C}=\mathcal{C}^\pm$ be a subset of $\mathcal{R}$
which contains $\mathcal{C}_1$.  We show that
$S(\mathcal{C})$ deformation retracts to $S(\mathcal{C}_1)$
by induction on the cardinality of $\mathcal{C} - \mathcal{C}_1$.

Choose $\alpha \in \mathcal{C} - \mathcal{C}_1$ which
satisfies both of:
\begin{enumerate}
\item The cardinality $|\alpha \cap G\mu|$ is minimal
(recall that $\mu$ is the maximally reductive ideal edge.)
\item The ideal edge $\alpha$ is minimal with respect to
property $1.$
\end{enumerate}

Using the Shrinking Lemma 7.4 of \cite{[K-V]}
with $\alpha$ and $\mu$, we obtain
a reductive ideal edge $\alpha_0 \subset \alpha$ which is
compatible with $\mu$.
Let $\gamma_i$ be the intersection components of
$\alpha$ with $\mu$ and index them so that
$m \in \gamma_0$.
Now from the Shrinking Lemma, we
can choose $\alpha_0$ so that it
is either one of the intersection
components $\gamma_{i_j}$ of $\alpha$ with $\mu$ which contain
no translate of $m$, or it is $\alpha - \cup \gamma_{i_j}$.
Because $\alpha \in \mathcal{C}
- \mathcal{C}_1$, $stab(\alpha) \not = G$ and $\alpha$ is neither
invertible nor equal to $\gamma = E_* - \{c^{-1}\}$.

\begin{claim} \label{t27} For every $\beta \in \mathcal{C}$, if $G\beta$
is compatible with $G\alpha$, then $G\beta$ is compatible with
$G\alpha_0$.
\end{claim}

\PF The three cases are
\begin{enumerate}
\item $G\alpha \subseteq G\beta$.  In this case, $G\alpha_0
\subseteq G\beta$ as $G\alpha_0 \subseteq G\alpha$.
\item $G\alpha \cap G\beta = \void$.
It follows that
$G\alpha_0 \cap G\beta = \void$ since $G\alpha_0 \subseteq G\alpha$.
\item $G\beta \subseteq G\alpha$.
Without loss of generality $\beta \subseteq \alpha$.
If $\beta \not \in \mathcal{C}_1$, then the minimality
conditions on $\alpha$ imply that $\beta=\alpha$, in which case
$\beta$ is clearly compatible with $\alpha_0$.  So
assume $\beta \in \mathcal{C}_1$.
As $G\beta \subseteq G\alpha$, $stab(\beta) \not = G$
since $stab(\alpha) \not = G$.
So either $\beta \in \mathcal{C}_0$ or $m \in G\beta$ and
$N(G\beta,G\mu)=1$.
If $\beta \in \mathcal{C}_0$ then either
$G\beta \subseteq G\mu$ (in which case $\beta$ is in some
$\gamma_i$ and thus compatible with $\alpha_0$),
$G\mu \subseteq G\beta$ (which can not happen as then
$\alpha$ would be in $\mathcal{C}_0$), or
$G\beta \cap G\mu = \void$ (in which case
$\beta \subseteq \alpha - \cup \gamma_{i_j}$ and thus
compatible with $\alpha_0$.)  Finally,
if $m \in G\beta$ and
$N(G\beta,G\mu)=1$ then $\beta \cap G\mu$ is not in
any of the $\gamma_{i_j}$'s as those are the intersection
components of $\alpha$ with $\mu$ that do not
contain a translate of $m$.  In fact, $\beta \cap G\mu$ is
in $\gamma_0$ and
$\beta \subseteq \alpha - \cup \gamma_{i_j}$.  Thus
$\beta$ is compatible with every choice of $\alpha_0$.\END
\end{enumerate}

Define a poset map
$f : S(\mathcal{C}) \to S(\mathcal{C})$ by sending
an ideal forest $\Phi$ to $\Phi \cup \{G\alpha_0\}$ if
$\Phi$ contains $\alpha$ and
to itself otherwise.  By the Poset Lemma, the
image of $f$ is a deformation retract of $S(\mathcal{C})$,
because $\Phi \subseteq f(\Phi)$
for all $\Phi$.
Define another poset map
$g : f(S(\mathcal{C})) \to f(S(\mathcal{C}))$ by sending
an ideal forest $\Psi$ to $\Psi - \{G\alpha\}$ if
$\Psi$ contains $\alpha$ and
to itself otherwise.  By the Poset Lemma, the
image of $g$ is a deformation retract of $f(S(\mathcal{C}))$,
as $g(\Psi) \subseteq \Psi$ for all $\Psi$.
Hence $S(\mathcal{C})$ deformation
retracts to $S(\mathcal{C}-\{G\alpha\})$, completing
the induction step. \END

\begin{lemma} \label{t28}   The complex $S(\mathcal{C}_1)$ deformation
retracts onto $S(\mathcal{C}_0^{'})$.
\end{lemma}

\PF Let $\mathcal{C}$ be a subset of $\mathcal{C}_1$
which contains $\mathcal{C}_0^{'}$.  We show that
$S(\mathcal{C})$ deformation retracts to $S(\mathcal{C}_0^{'})$
by induction on the cardinality of $\mathcal{C} - \mathcal{C}_0^{'}$.

Choose $\alpha \in \mathcal{C} - \mathcal{C}_0^{'}$
such that $m \in \alpha$ and
\begin{enumerate}
\item The cardinality $|\alpha \cap G\mu|$ is minimal.
\item The ideal edge $\alpha$ is minimal with respect to
property $1.$
\end{enumerate}

We apply the Pushing Lemma 7.3 of \cite{[K-V]} (just look
at the actual results that each of the cases in the proof
gives you) to get a reductive edge $\alpha_0$
with $\alpha_0=\alpha \cap \mu$ or $\alpha_0 = \alpha - \mu$.
Note that $\alpha_0 \in \mathcal{C}_0$.

\begin{claim} \label{t29} For every $\beta \in \mathcal{C}$, if $G\beta$
is compatible with $G\alpha$, then $G\beta$ is compatible with
$G\alpha_0$.
\end{claim}

\PF The three cases are
\begin{enumerate}
\item $G\alpha \subseteq G\beta$.  In this case, $G\alpha_0
\subseteq G\beta$ as $G\alpha_0 \subseteq G\alpha$.
\item $G\alpha \cap G\beta = \void$.
It follows that
$G\alpha_0 \cap G\beta = \void$ since $G\alpha_0 \subseteq G\alpha$.
\item $G\beta \subseteq G\alpha$.
Without loss of generality $\beta \subseteq \alpha$.
If $\beta \not \in \mathcal{C}_0^{'}$, then the minimality
conditions on $\alpha$ imply that $\beta=\alpha$, in which case
$\beta$ is clearly compatible with $\alpha_0$.  So
assume $\beta \in \mathcal{C}_0^{'}$.
Since $\alpha \in \mathcal{C} - \mathcal{C}_0^{'}$,
$stab(\alpha) \not = G$.
As $G\beta \subseteq G\alpha$, $stab(\beta) \not = G$
also, which means that $\beta \in \mathcal{C}_0$.
The three ways in which $\beta$ could be compatible with
$\mu$ are:
\begin{itemize}
\item $G\beta \cap G\mu = \void$.
Then $G\beta$ is
disjoint from $G(\alpha \cap \mu)$ and
contained in $G(\alpha - \mu)$.
\item $G\beta \subseteq G\mu$.
Then $G\beta \subseteq G(\alpha \cap \mu)$ and
$G\beta$ is disjoint from $G(\alpha - \mu)$.
\item $G\mu \subseteq G\beta$.
Then $G\mu \subseteq G\alpha$ and so $G\mu$ and $G\alpha$
are compatible, a contradiction.\END
\end{itemize}
\end{enumerate}

Define a poset map
$f : S(\mathcal{C}) \to S(\mathcal{C})$ by sending
an ideal forest $\Phi$ to $\Phi \cup \{G\alpha_0\}$ if
$\Phi$ contains $\alpha$ and
to itself otherwise.  By the Poset Lemma, the
image of $f$ is a deformation retract of $S(\mathcal{C})$,
because $\Phi \subseteq f(\Phi)$
for all $\Phi$.
Define another poset map
$g : f(S(\mathcal{C})) \to f(S(\mathcal{C}))$ by sending
an ideal forest $\Psi$ to $\Psi - \{G\alpha\}$ if
$\Psi$ contains $\alpha$ and
to itself otherwise.  By the Poset Lemma, the
image of $g$ is a deformation retract of $f(S(\mathcal{C}))$,
as $g(\Psi) \subseteq \Psi$ for all $\Psi$.
Hence $S(\mathcal{C})$ deformation
retracts to $S(\mathcal{C}-\{G\alpha\})$, completing
the induction step. \END

Now we are left with the task of showing that $S(\mathcal{C}_0^{'})$
deformation retracts to $S(\mathcal{C}_0)$
and from there to a point.  From
Proposition \ref{t25}, we see that this can be handled in
three separate cases:
\begin{itemize}
\item The ideal edge $\mu$ is invertible and the
reductive ideal edge $\gamma = E_* - \{c^{-1}\}$ is not
compatible with $\mu$.  In this case, the proposition
gives us that $c=m^{-1}$.
\item The ideal edge $\mu$ is invertible and the
reductive ideal edge $\gamma = E_* - \{c^{-1}\}$ is
compatible with $\mu$.
\item The ideal edge $\mu$ is not invertible.
\end{itemize}

\begin{lemma} \label{t30} Suppose $\mu$ is invertible
and $\gamma= E_* - \{m\}$ is reductive.  Then
$S(\mathcal{C}_0^{'})$ is contractible.
\end{lemma}

\PF We first contract $S(\mathcal{C}_0^{'})$ to
$S(\mathcal{C}_0 \cup \{\gamma\})$.  Let
$\mathcal{C}$ be a subset of $\mathcal{C}_0^{'}$
which contains $\mathcal{C}_0 \cup \{\gamma\}$.
Also assume that if $\alpha \in \mathcal{C}$ is not
pre-compatible with $\mu$, then $\alpha^{-1} \in \mathcal{C}$
also.  We will use induction
on $|\mathcal{C} - (\mathcal{C}_0 \cup \{\gamma\})|$
to show that $S(\mathcal{C})$ deformation retracts to 
$S(\mathcal{C}_0 \cup \{\gamma\})$.

Choose $\alpha \in \mathcal{C} - (\mathcal{C}_0 \cup \{\gamma\})$
such that $m \in \alpha$ and
\begin{enumerate}
\item The cardinality $|\alpha \cap \mu|$ is maximal.
\item The edge $\alpha$ is maximal with respect to property $1.$
\end{enumerate}

There are two main cases, and two subcases in the second case.

\medskip

\noindent {\em Case 1.} $\alpha^{-1}$ is compatible with $\mu$.

Then $\mu \not \subseteq \alpha^{-1}$ because
$m \in \mu$ and $m \in \alpha$. Also,
$\mu \cap \alpha^{-1} \not = \void$ else $\mu \subseteq \alpha$
and $\alpha$ is compatible with $\mu$.  So $\alpha^{-1} \subseteq \mu$.
Let $\alpha_0 = \alpha^{-1}$ and note that $\alpha_0 \in \mathcal{C}$.

\begin{claim} \label{t31} For every $\beta$ in $\mathcal{C}$,
if $G\beta$ is compatible with $\alpha$, then $G\beta$
is compatible with $\alpha_0$.
\end{claim}

\PF If $G\beta \subseteq \alpha$ then $G\beta \cap \alpha_0 = \void$.
On the other hand, if $G\beta \cap \alpha = \void$, then
$G\beta \subseteq \alpha_0$.

Finally, assume $\alpha \subseteq G\beta$.  Note that
$\beta \not = \gamma$ as $m \in \alpha \subseteq G\beta$.
So $\beta$ is invertible and $G\beta=\beta$,
as $\alpha \subseteq G\beta$ and $\alpha$ is invertible.
Since $\alpha \subseteq \beta$ and $\alpha^{-1} \subseteq \mu$,
$\beta^{-1} \subseteq \mu$.
As $m \in \beta \cap \mu$ and $\beta^{-1} \subseteq \mu$,
$\beta$ is not compatible with $\mu$.
So $\beta \in \mathcal{C} - (\mathcal{C}_0 \cup \{\gamma\})$.
But $|\alpha \cap \mu| \leq |\beta \cap \mu|$ and
so the maximality conditions on $\alpha$ give
$\beta=\alpha$.  Hence $\beta$ is compatible with $\alpha_0$. \END

By the claim we can define 
$f : S(\mathcal{C}) \to S(\mathcal{C})$ by sending
an ideal forest $\Phi$ to $\Phi \cup \{G\alpha_0\}$ if
$\Phi$ contains $\alpha$ and
to itself otherwise.
We can also define
$g : f(S(\mathcal{C})) \to f(S(\mathcal{C}))$ by sending
an ideal forest $\Psi$ to $\Psi - \{G\alpha\}$ if
$\Psi$ contains $\alpha$ and
to itself otherwise.
Then $g(f(S(\mathcal{C})))=S(\mathcal{C} - \{\alpha\})$
and $S(\mathcal{C})$ deformation retracts to this by
two applications of the Poset Lemma.

\medskip

\noindent {\em Case 2.} $\alpha^{-1}$ is not compatible with $\mu$.

Since $\alpha$ and $\mu$ cross simply (this is automatic because
$\alpha$ is invertible), the Pushing Lemma applies.
Thus one of the sets $\alpha_0 = \mu - \alpha$ or
$\alpha_0 = \alpha \cup \mu$ is a reductive ideal edge.
As $\mu - \alpha \subseteq \mu$ and $\mu \subseteq \alpha \cup \mu$,
$\alpha_0 \in \mathcal{C}_0$ in either case.

\medskip

\noindent {\em Subcase 1.} $\alpha_0 = \mu - \alpha$.

\begin{claim} \label{t32} For every $\beta \in \mathcal{C}$,
if $G\beta$ is pre-compatible with $\alpha$, then
$G\beta$ is compatible with $\alpha_0$.
\end{claim}

\PF We list what happens in each of the possible ways that
$G\beta$ can be pre-compatible with $\alpha$.
\begin{itemize}
\item $\beta = \gamma$.  Now $m \in \alpha$ and thus
$m \not \in \mu - \alpha$.  So $\alpha_0 \subseteq \gamma$.
\item $G\beta \subseteq \alpha$.  Then $G\beta \cap \alpha_0 = \void$.
\item $\beta$ is invertible and $\beta^{-1} \subseteq \alpha$.
Then $\alpha_0 \subseteq \beta$.
\item $\alpha \subseteq G\beta$.  Then $\beta$ is invertible
and $\alpha \subseteq \beta$.  Also $|\alpha \cap \mu|
\leq |\beta \cap \mu|$.  If $\beta \not \in \mathcal{C}_0 \cup \{\gamma\}$
then the maximality conditions on $\alpha$ give us that $\beta=\alpha$.
As $\alpha \cap \alpha_0 = \void$, they are compatible
in this situation.  Now $\beta \not = \gamma$ as
$m \in \alpha \subseteq \beta$ (or, alternatively, just because
we already handled that case above.)  So the
possibility left is that $\beta \in \mathcal{C}_0$.
If $\beta \subseteq \mu$, then we get a contradiction
because then $\alpha \subseteq \beta \subseteq \mu$ and
$\alpha \in \mathcal{C}_0$.
Similarly, since $m \in \beta \cap \mu$, we also
can not have $\beta \cap \mu = \void$.
So $\mu \subseteq \beta$ and thus
$\beta \supseteq \mu - \alpha = \alpha_0$.
\item $\alpha \cap G\beta = \void$.

If $\beta$ is
invertible, then $\alpha \subseteq \beta^{-1}$
and we can apply the methods of the previous case.
Namely, if $\beta^{-1} \in \mathcal{C}_0$ then
$\mu \subseteq \beta^{-1}$ is the only possible case.
So $\beta$ is disjoint from $\mu - \alpha = \alpha_0$.
Also if $\beta \subseteq \mu$ then
$\beta \subseteq \mu - \alpha = \alpha_0$
since $\alpha \cap G\beta = \void$.
Alternatively, if we assume
$\beta$ is not pre-compatible with $\mu$, then
$\beta^{-1} \in \mathcal{C}$ as $\beta \in \mathcal{C}$.
So we have $\beta^{-1} \in \mathcal{C} - (\mathcal{C}_0 \cup \{\gamma\})$
and $\alpha \subseteq \beta^{-1}$.  Hence the
maximality conditions on $\alpha$ yield that
$\beta^{-1}=\alpha$.
Since $\alpha_0 = \mu - \alpha \subseteq \alpha^{-1} = \beta$,
$\alpha_0$ is compatible with $\beta$ here also.

Next assume $\beta$ is not invertible.
Consequently, $\beta \in \mathcal{C}_0$ and we have the
standard three possibilities for how $\beta$ and $\mu$
are compatible.  If $G\beta \subseteq \mu$, then $G\beta \subseteq \alpha_0$.
On the other hand, if $G\beta \cap \mu = \void$ then
$G\beta \cap \alpha_0 = \void$.
Finally, since $m \in G\beta \cap \alpha$, we can not
have $G\beta \supseteq \mu$. \END
\end{itemize}

So by the claim we can define 
$f : S(\mathcal{C}) \to S(\mathcal{C})$ by sending
an ideal forest $\Phi$ to $\Phi \cup \{G\alpha_0\}$ if
$\Phi$ contains $\alpha$ or $\alpha^{-1}$ and
to itself otherwise.
We can also define
$g : f(S(\mathcal{C})) \to f(S(\mathcal{C}))$ by sending
an ideal forest $\Psi$ to $\Psi - \{\alpha,\alpha^{-1}\}$ if
$\Psi$ contains $\alpha$ or $\alpha^{-1}$ and
to itself otherwise.
Then $g(f(S(\mathcal{C})))=S(\mathcal{C} - \{\alpha,\alpha^{-1}\})$
and $S(\mathcal{C})$ deformation retracts to this by
two applications of the Poset Lemma.

\medskip

\noindent{Subcase 2.} $\alpha_0 = \alpha \cup \mu$.
The ideal edge $\alpha_0$ does not equal $\gamma$ since it
is out-reductive by the proof of the Pushing Lemma (because
$\alpha$ and $\mu$ are not equal to $\gamma$ and thus
out-reductive.)  (Alternatively, it is clear that
$\alpha_0 \not = \gamma$ because $m \in \alpha_0$.)
Accordingly, $\alpha_0$ is invertible and
both $\alpha_0$ and $\alpha_0^{-1}$ are in
$\mathcal{C}_0$.

\begin{claim} \label{t33}  For every $\beta \in \mathcal{C}$,
if $G\beta$ is compatible with $\alpha^{-1}$ then
$G\beta$ is compatible with $\alpha_0^{-1}$.
\end{claim}

\PF We list what happens in each of the possible ways that
$G\beta$ can be compatible with $\alpha^{-1}$.
\begin{itemize}
\item $\beta = \gamma$.  Note that $m \in \alpha_0$
and consequently $m \not \in \alpha_0^{-1}$.  For this
reason, $\alpha_0^{-1} \subseteq \beta$ and they are compatible.
Assume $\beta \not = \gamma$ from now on.
\item $G\beta \cap \alpha^{-1} = \void$.
Hence $G\beta \subseteq \alpha \subseteq \alpha \cup \mu = \alpha_0$
and $G\beta \cap \alpha_0^{-1} = \void$.
\item $\alpha^{-1} \subseteq G\beta$.  In this case,
$\alpha_0^{-1} \subseteq \alpha^{-1} \subseteq G\beta$.
\item $G\beta \subseteq \alpha^{-1}$.

Suppose $\beta$ is invertible.  Then $\alpha \subseteq \beta^{-1}$.
Now if $\beta \subseteq \mu$ then $\beta \cap \alpha_0^{-1} = \void$
and they are compatible.  Also if $\beta^{-1} \in \mathcal{C}_0$,
then neither $\beta^{-1} \subseteq \mu$ nor
$\beta^{-1} \cap \mu = \void$ can occur because the first
contradicts the fact that $\alpha$ is not compatible with $\mu$
and the second contradicts the fact that $m \in \beta^{-1} \cap \mu$.
Accordingly, $\mu \subseteq \beta^{-1}$ and
$\beta^{-1} \supseteq \alpha \cup \mu = \alpha_0$.
Thus $\alpha_0^{-1} \supseteq \beta$ and $\beta$ is
compatible with $\alpha_0^{-1}$.
So the last possibility is that $\beta$ is not pre-compatible
with $\mu$.  Then our assumptions on $\mathcal{C}$ imply
that $\beta^{-1} \in \mathcal{C}$ also.  Since
$\alpha \subseteq \beta^{-1}$, the maximality conditions
on $\alpha$ yield that $\alpha = \beta^{-1}$.  It follows
that $\alpha_0^{-1} \subseteq \beta$ and they are compatible.

If we suppose $\beta$ is not invertible, then $\beta \in \mathcal{C}_0$
and there are three possibilities.  If $G\beta \subseteq \mu$,
then $G\beta \cap \alpha_0^{-1} = \void$.
If $G\beta \cap \mu = \void$ then $G\beta \subseteq (\alpha \cup \mu)^{-1}
= \alpha_0^{-1}$ since $G\beta \subseteq \alpha^{-1}$ also.
Finally, because $m \not \in \alpha^{-1}$ and $G\beta \subseteq \alpha^{-1}$,
the last case $G\beta \supseteq \mu$ will not happen.\END
\end{itemize}

Therefore we can define 
$f : S(\mathcal{C}) \to S(\mathcal{C})$ by sending
an ideal forest $\Phi$ to $\Phi \cup \{\alpha_0^{-1}\}$ if
$\Phi$ contains $\alpha^{-1}$ and
to itself otherwise.
We can also define
$g : f(S(\mathcal{C})) \to f(S(\mathcal{C}))$ by sending
an ideal forest $\Psi$ to $\Psi - \{\alpha^{-1}\}$ if
$\Psi$ contains $\alpha^{-1}$ and
to itself otherwise.
Then $g(f(S(\mathcal{C})))=S(\mathcal{C} - \{\alpha^{-1}\})$
and $S(\mathcal{C})$ deformation retracts to this by
two applications of the Poset Lemma.

We must further contract $S(\mathcal{C} - \{\alpha^{-1}\})$
to $S(\mathcal{C} - \{\alpha,\alpha^{-1}\})$.  This is
done via the following claim:

\begin{claim} \label{t34} For every $\beta \in \mathcal{C} - \{\alpha^{-1}\}$,
if $G\beta$ is compatible with $\alpha$ then $G\beta$ is
compatible with $\alpha_0$.
\end{claim}

\PF Note that $\gamma$ is not compatible with $\alpha$.
As before, we list the standard cases:
\begin{itemize}
\item $G\beta \subseteq \alpha$.
Then $G\beta \subseteq \alpha \cup \mu = \alpha_0$.
\item $\alpha \subseteq G\beta$.  Then $\beta$ is
invertible and $\alpha \subseteq \beta$.
Note that $|\alpha \cap \mu| \leq |\beta \cap \mu|$.
If $\beta \not \in \mathcal{C}_0 \cup \{\gamma\}$ then
the maximality conditions on $\alpha$ give us that
$\beta = \alpha$.  As $\alpha \subseteq \alpha \cup \mu = \alpha_0$,
they are compatible in this situation.
Next we examine what happens if $\beta \in \mathcal{C}_0$.
Neither the case $\beta \subseteq \mu$ nor
$\beta \cap \mu = \void$ can occur as the first contradicts the
fact that $\alpha$ is not compatible with $\mu$ and the second
that $m \in \beta \cap \mu$.
Accordingly, $\mu \subseteq \beta$ and
$\beta \supseteq \alpha \cup \mu = \alpha_0$.
\item $\alpha \cap G\beta = \void$.

If $\beta$ is invertible, then $\alpha \subseteq \beta^{-1}$
and we can apply similar methods to those used in the last case.
Namely, if $\beta^{-1} \in \mathcal{C}_0$ then
$\mu \subseteq \beta^{-1}$ is the only possible case.
So $\beta^{-1} \supseteq \alpha \cup \mu$ and
$\alpha_0 \cap \beta = \void$.
Also if $\beta \subseteq \mu$ then
$\beta \subseteq \alpha \cup \mu = \alpha_0$.
Hence the last case left is where $\beta$ is not
pre-compatible with $\mu$.  Then $\beta^{-1} \in \mathcal{C}$
because $\beta \in \mathcal{C}$.
Since $\alpha \subseteq \beta^{-1}$, the maximality conditions on
$\alpha$ yield that $\alpha = \beta^{-1}$.
As $\beta$ lies in $\mathcal{C} - \{\alpha^{-1}\}$, this
gives a contradiction.  So the last case will
not actually occur.

Now assume $\beta$ is not invertible,
so that $\beta \in \mathcal{C}_0$.
If $G\beta \subseteq \mu$, then $G\beta \subseteq \alpha \cup \mu =
\alpha_0$.
Alternatively, if $G\beta \cap \mu = \void$
then $G\beta \cap (\alpha \cup \mu) = \void$ as
$\alpha \cap G\beta = \void$ also.
Finally, the case $G\beta \supseteq \mu$ does
not take place, else $m \in G\beta \cap \alpha$. \END
\end{itemize}

Using this last claim we can define 
$f : S(\mathcal{C} - \{\alpha_0^{-1}\})
\to S(\mathcal{C} - \{\alpha_0^{-1}\})$ by sending
an ideal forest $\Phi$ to $\Phi \cup \{\alpha_0\}$ if
$\Phi$ contains $\alpha$ and
to itself otherwise.
We can also define
$g : f(S(\mathcal{C} - \{\alpha_0^{-1}\}))
\to f(S(\mathcal{C} - \{\alpha_0^{-1}\}))$ by sending
an ideal forest $\Psi$ to $\Psi - \{\alpha\}$ if
$\Psi$ contains $\alpha$ and
to itself otherwise.
Then $g(f(S(\mathcal{C})))=S(\mathcal{C} - \{\alpha,\alpha^{-1}\})$
and $S(\mathcal{C})$ deformation retracts to this by
two applications of the Poset Lemma.

\medskip

This concludes our argument that
$S(\mathcal{C}_0^{'})$ contracts to
$S(\mathcal{C}_0 \cup \{\gamma\})$.  To eliminate
$\gamma$, note that $\gamma$ is compatible with
$\mu^{-1} \in \mathcal{C}_0$ and proceed to the
following claim:

\begin{claim} \label{t35}  For every $\beta \in \mathcal{C}_0$,
if $G\beta$ is compatible with $\gamma$ then $G\beta$
is compatible with $\mu^{-1}$.
\end{claim}

\PF Observe that if $G\beta$ is not invertible then
the fact that it is compatible with $\mu$ means it is
compatible with $\mu^{-1}$. So we can assume $\beta$
is invertible.  Since $\beta$ is compatible
with $\gamma$, $m \not \in \beta$.
Because $m \in \mu$, we can not have
$\mu \subseteq \beta$.  If $\mu \cap \beta = \void$
then $\beta \subseteq \mu^{-1}$ and we are fine.
Also if $\beta \subseteq \mu$ then
$\beta \cap \mu^{-1} = \void$ and $\beta$ and
$\mu^{-1}$ are compatible here also. \END

Now define poset maps $f$ and $g$ analogously to
all of the times we did this above to add
$\mu^{-1}$ to ideal forests containing $\gamma$,
and then to delete $\gamma$.  This allows us to
deformation retract
$S(\mathcal{C}_0 \cup \{\gamma\})$
to $S(\mathcal{C})$.

The final step of contracting $S(\mathcal{C}_0)$
to a point is done with poset
maps $f$ and $g$ as above which first add $\mu$ to all
ideal forests and then delete from these ideal forests
everything except $\mu$. \END

\begin{lemma} \label{t36}  Suppose $\mu$ is invertible
and the reductive $\gamma= E_* - \{c^{-1}\}$ is
compatible with $\mu$.  Then
$S(\mathcal{C}_0^{'})$ is contractible.
\end{lemma}

\PF The proof of the more complicated case in
Lemma \ref{t30} carries over to this one, with the
exception that the penultimate step of
deformation retracting from
$S(\mathcal{C}_0 \cup \{\gamma\})$
to $S(\mathcal{C})$ is unnecessary, because
$\gamma$ is already compatible with $\mu$.
In addition, various other minor changes need to
be made to the claims listed under Lemma \ref{t30}, because
$\gamma$ is now compatible with $\mu$, but the changes
are all easy and direct. \END

\begin{lemma} \label{t37} Suppose $\mu$ is not invertible.
Then
$S(\mathcal{C}_0^{'})$ is contractible.
\end{lemma}

\PF As before, let $\gamma = E_* - \{c^{-1}\}$ be the
reductive edge that Proposition \ref{t25} gives us (if it exists).
We know that $\gamma$ is compatible with $\mu$ because
$stab(\mu) \not = G = stab(c^{-1})$.

The proof of this lemma is basically the proof in
Subcase 2, Case 2 of Lemma \ref{t30} above; however, since
$\mu$ not being invertible changes a few things,
it seems useful to state the entire proof here.

We first contract $S(\mathcal{C}_0^{'})$ to
$S(\mathcal{C}_0)$.  Let
$\mathcal{C}$ be a subset of $\mathcal{C}_0^{'}$
which contains $\mathcal{C}_0$.
Also assume that if $\alpha \in \mathcal{C}$
and $\alpha$ is invertible, then $\alpha^{-1} \in \mathcal{C}$
also.  We will use induction
on $|\mathcal{C} - \mathcal{C}_0|$
to show that $S(\mathcal{C})$ deformation retracts to 
$S(\mathcal{C}_0)$.

Choose $\alpha \in \mathcal{C} - (\mathcal{C}_0)$
such that $m \in \alpha$ and
\begin{enumerate}
\item The cardinality $|\alpha \cap G\mu|$ is maximal.
\item The edge $\alpha$ is maximal with respect to property $1.$
\end{enumerate}

Since $\alpha$ and $\mu$ cross simply (this is automatic because
$\alpha$ is invertible), the Pushing Lemma applies.
Say $\alpha = (\alpha,a)$.  Neither $a$ nor $a^{-1}$
is in $\mu$ since $stab(a)=G$ and $\mu$ is not invertible
(and not equal to $\gamma$.)
So case 1 of the Pushing Lemma shows that
$\alpha_0 = \alpha \cup G\mu$ is a reductive ideal edge.
As $\mu \subseteq \alpha \cup G\mu$,
$\alpha_0 \in \mathcal{C}_0$.
The ideal edge $\alpha_0$ is not equal to $\gamma$
because it is out-reductive by the proof of the Pushing
Lemma (as both $\alpha$ and $\mu$ are out-reductive.)
So $\alpha_0$ is invertible and
$\alpha_0^{-1} \in \mathcal{C}_0$ also.

\begin{claim} \label{t38} For every $\beta \in \mathcal{C}$,
if $G\beta$ is compatible with $\alpha^{-1}$ then
$G\beta$ is compatible with $\alpha_0^{-1}$.
\end{claim}

\PF We list what happens in each of the possible ways that
$G\beta$ can be compatible with $\alpha^{-1}$.
\begin{itemize}
\item $\beta = \gamma$.  Then $\alpha^{-1} \subseteq \gamma$
so $c^{-1} \not \in \alpha^{-1}$.
Thus $c^{-1} \in \alpha$.
So $c^{-1} \in \alpha \cup G\mu = \alpha_0$
and $c^{-1} \not \in \alpha_0^{-1}$.
Hence $\gamma$ is compatible with $\alpha_0^{-1}$.
Assume $\beta \not = \gamma$ from now on.
\item $G\beta \cap \alpha^{-1} = \void$.
Hence $G\beta \subseteq \alpha \subseteq \alpha \cup G\mu = \alpha_0$
and $G\beta \cap \alpha_0^{-1} = \void$.
\item $\alpha^{-1} \subseteq G\beta$.  In this case,
$\alpha_0^{-1} \subseteq \alpha^{-1} \subseteq G\beta$.
\item $G\beta \subseteq \alpha^{-1}$.

Suppose $\beta$ is invertible.  Then $\alpha \subseteq \beta^{-1}$.
If $\beta^{-1} \in \mathcal{C}_0$,
then neither $\beta^{-1} \subseteq G\mu$ nor
$\beta^{-1} \cap G\mu = \void$ can occur because the first
contradicts the fact that $\alpha$ is not compatible with $G\mu$
(or alternatively the fact that $\mu$ is not
invertible)
and the second contradicts the fact that $m \in \beta^{-1} \cap G\mu$.
Accordingly, $G\mu \subseteq \beta^{-1}$ and
$\beta^{-1} \supseteq \alpha \cup G\mu = \alpha_0$.
Thus $\alpha_0^{-1} \supseteq \beta$ and $\beta$ is
compatible with $\alpha_0^{-1}$.
So the last possibility is that
$\beta^{-1} \in \mathcal{C} - \mathcal{C}_0$.
Since
$\alpha \subseteq \beta^{-1}$, the maximality conditions
on $\alpha$ yield that $\alpha = \beta^{-1}$.  It follows
that $\alpha_0^{-1} \subseteq \beta$ and they are compatible.

If we suppose $\beta$ is not invertible, then $\beta \in \mathcal{C}_0$
and there are three possibilities.  If $G\beta \subseteq G\mu$,
then $G\beta \cap \alpha_0^{-1} = \void$.
If $G\beta \cap G\mu = \void$ then
$G\beta \subseteq (\alpha \cup G\mu)^{-1}
= \alpha_0^{-1}$ since $G\beta \subseteq \alpha^{-1}$ also.
Finally, because $m \not \in \alpha^{-1}$ and
$G\beta \subseteq \alpha^{-1}$,
the last case $G\beta \supseteq G\mu$ will not happen.\END
\end{itemize}

Therefore we can define 
$f : S(\mathcal{C}) \to S(\mathcal{C})$ by sending
an ideal forest $\Phi$ to $\Phi \cup \{\alpha_0^{-1}\}$ if
$\Phi$ contains $\alpha^{-1}$ and
to itself otherwise.
We can also define
$g : f(S(\mathcal{C})) \to f(S(\mathcal{C}))$ by sending
an ideal forest $\Psi$ to $\Psi - \{\alpha^{-1}\}$ if
$\Psi$ contains $\alpha^{-1}$ and
to itself otherwise.
Then $g(f(S(\mathcal{C})))=S(\mathcal{C} - \{\alpha^{-1}\})$
and $S(\mathcal{C})$ deformation retracts to this by
two applications of the Poset Lemma.

We must further contract $S(\mathcal{C} - \{\alpha^{-1}\})$
to $S(\mathcal{C} - \{\alpha,\alpha^{-1}\})$.  This is
done via the following claim:

\begin{claim} \label{t39} For every $\beta \in \mathcal{C} - \{\alpha^{-1}\}$,
if $G\beta$ is compatible with $\alpha$ then $G\beta$ is
compatible with $\alpha_0$.
\end{claim}

\PF As before, we list the standard cases:
\begin{itemize}
\item $\beta=\gamma$. Then $c^{-1} \not \in \alpha$.
As $c^{-1} \not \in G\mu$ also, this means
$c^{-1} \not \in \alpha \cup G\mu = \alpha_0$.
So $\gamma$ is compatible with $\alpha_0$.
Assume $\beta \not = \gamma$ from now on.
\item $G\beta \subseteq \alpha$.
Then $G\beta \subseteq \alpha \cup G\mu = \alpha_0$.
\item $\alpha \subseteq G\beta$.  Then $\beta$ is
invertible and $\alpha \subseteq \beta$.
Note that $|\alpha \cap G\mu| \leq |\beta \cap G\mu|$.
If $\beta \not \in \mathcal{C}_0$ then
the maximality conditions on $\alpha$ give us that
$\beta = \alpha$.  As $\alpha \subseteq \alpha \cup G\mu = \alpha_0$,
they are compatible in this situation.
Next we examine what happens if $\beta \in \mathcal{C}_0$.
Neither the case $\beta \subseteq G\mu$ nor
$\beta \cap G\mu = \void$ can occur as the first contradicts the
fact that $\mu$ is not invertible and the second
that $m \in \beta \cap G\mu$.
Accordingly, $G\mu \subseteq \beta$ and
$\beta \supseteq \alpha \cup G\mu = \alpha_0$.
\item $\alpha \cap G\beta = \void$.

If $\beta$ is invertible, then $\alpha \subseteq \beta^{-1}$
and we can apply similar methods to those used in the last case.
Namely, if $\beta^{-1} \in \mathcal{C}_0$ then
$G\mu \subseteq \beta^{-1}$ is the only possible case.
So $\beta^{-1} \supseteq \alpha \cup G\mu$ and
$\alpha_0 \cap \beta = \void$.
On the other hand,
suppose $\beta \in \mathcal{C} - \mathcal{C}_0$.
Since $\alpha \subseteq \beta^{-1}$, the maximality conditions on
$\alpha$ yield that $\alpha = \beta^{-1}$.
As $\beta$ lies in $\mathcal{C} - \{\alpha^{-1}\}$, this
gives a contradiction.  So the last case will
not actually occur.

Now assume $\beta$ is not invertible,
so that $\beta \in \mathcal{C}_0$.
If $G\beta \subseteq G\mu$, then $G\beta \subseteq \alpha \cup G\mu =
\alpha_0$.
Alternatively, if $G\beta \cap G\mu = \void$
then $G\beta \cap (\alpha \cup G\mu) = \void$ as
$\alpha \cap G\beta = \void$ also.
Finally, the case $G\beta \supseteq G\mu$ does
not take place, else $m \in G\beta \cap \alpha$. \END
\end{itemize}

Using this last claim we can define 
$f : S(\mathcal{C} - \{\alpha_0^{-1}\})
\to S(\mathcal{C} - \{\alpha_0^{-1}\})$ by sending
an ideal forest $\Phi$ to $\Phi \cup \{\alpha_0\}$ if
$\Phi$ contains $\alpha$ and
to itself otherwise.
We can also define
$g : f(S(\mathcal{C} - \{\alpha_0^{-1}\}))
\to f(S(\mathcal{C} - \{\alpha_0^{-1}\}))$ by sending
an ideal forest $\Psi$ to $\Psi - \{\alpha\}$ if
$\Psi$ contains $\alpha$ and
to itself otherwise.
Then $g(f(S(\mathcal{C})))=S(\mathcal{C} - \{\alpha,\alpha^{-1}\})$
and $S(\mathcal{C})$ deformation retracts to this by
two applications of the Poset Lemma.

\medskip

This concludes our argument that
$S(\mathcal{C}_0^{'})$ contracts to
$S(\mathcal{C}_0)$. 
The final step of contracting $S(\mathcal{C}_0)$
to a point is done with poset
maps $f$ and $g$ as above which first add $G\mu$ to all
ideal forests and then delete from these ideal forests
everything except $G\mu$. \END

The sequence of lemmas above concludes our proof of
Theorem \ref{tr13}.

Back in the 
beginning of this paper, in
Definition \ref{tr17},
we asserted that the dimension of $X_p^\omega$
was $2p-4$.  We are now in a position to justify this.

\begin{remark} \label{tr29}
The dimension of $X_p^\omega$ is $2p-4$.
\end{remark}

\PF Just as asserted in \S 9 of \cite{[K-V]}, every maximal
simplex in $X_p^\omega$ has the same dimension.
We can obtain a maximal simplex of $X_p^\omega$
by first taking a $(2p-6)$-dimensional maximal simplex in $X_{p-2}$
which has no inessential edges under the
action of the trivial group $1$, say the simplex
$$\eta : R_{p-2} \to \Gamma_{p-2}^0 \to \Gamma_{p-2}^1 \to \cdots \to \Gamma_{p-2}^{2p-6},$$
and then modifying the simplex as follows.
Choose some edge $e$ of $\Gamma_{p-2}^0$
which survives (uncollapsed) to the graph $\Gamma_{p-2}^{2p-6}$.
and pick two points $x_1 \not = x_2$ in the interior of
that edge. 
Now attach two other edges $f_1$
and $f_2$ to the graph $\Gamma_{p-2}^0$ by letting
them both start at $x_1$ and both end at $x_2$.
Call the resulting graph $\Gamma_p^0$
and define an $\omega$-action on it by having $\omega$
switch the edges $f_1$ and $f_2$ and leave all other
edges fixed.
By adding the vertices $x_1$ and $x_2$ to the edge $e$,
we divided $e$ into three smaller edges, which we will
call $e_1$, $e_2$, and $e_3$.

In a similar manner, for all $i \in \{1, \ldots, 2p-6\}$,
define an $\omega$-graph $\Gamma_p^i$
by attaching
edges $f_1$ and $f_2$ to $\Gamma_{p-2}^i$.

Define $\Gamma_p^{2p-5}$ to be the result of collapsing the
edge $e_1$ in $\Gamma_p^{2p-6}$ and
define $\Gamma_p^{2p-4}$ to be the result of collapsing the
edge $e_2$ in $\Gamma_p^{2p-5}$.

Then it is easy to see that the chain of forest collapses
$$\Gamma_{p}^0 \to \Gamma_{p}^1 \to \cdots \to \Gamma_{p}^{2p-4}$$
gives us that the dimension of a maximal simplex in
$X_p^\omega$ is $2p-4$.  This is because the only places
that the graph $\Gamma_p^0$ can be blown up are at
the vertices $x_1$ or $x_2$, and any such blowup would result in
an inessential edge. \END

\part{Cohomology of moduli spaces of graphs} \label{p4}

\chapter{Symmetry groups of graphs} \label{c13}

We now examine the cohomology of the quotient $Q_n$ of the
spine $X_n$ of auter space.  
The are natural inclusions $Q_m \injarrow Q_{m+1}$, and it is
known \cite{[H-V]} that the induced map
$H^i(Q_{m+1}; \Q) \to H^i(Q_m; \Q)$
is an isomorphism for $m > 3i/2$.
Our goal is to show that, in contrast,
$H^5(Q_m; \Z)$
never stabilizes as $m \to \infty$.  This is done by showing
that as $m$ increases, torsion from increasingly higher primes is
introduced in $H^5(Q_m; \Z)$.  
To this end, we do specific calculations in the spectral sequence
\ref{e1} applied to the action of $Aut(F_n)$ on $X_n$ for
$n=2p-1$.
The $E_2^{r,0}$-term of this spectral sequence is $H^r(Q_n; \Z_{(p)})$,
and the sequence converges to $H^r(Aut(F_n); \Z_{(p)})$.
Results from Hatcher and Vogtmann \cite{[H-V]} on the
cohomology of $Aut(F_n)$ are then used to obtain the result.

In this chapter, we do the ground work necessary to compute the 
$E_1$-page of the spectral sequence: we find all simplices of $X_n$ with 
$p$-symmetry and compute the cohomology of the stabilizers
of these simplices with coefficients in $\Z_{(p)}$.
In Chapter \ref{c14} we will compute the $E_2$-page of the
spectral sequence, and use this calculation to obtain the result.


Unless otherwise stated, $p \geq 5$ will be prime and $n=2p-1$.
The assumption that $p \geq 5$ is for convenience more than
any other reason, as the main results will only consider
arbitrarily large primes $p$ and so we should not
devote extra time to the (fairly easy to resolve) complications
introduced by considering the prime $p=3$. 
These complications arise from the fact that
the dihedral
group $D_6$ is the same as the symmetric group $S_3$,
so that we cannot
distinguish between dihedral and symmetric symmetry in that case.

\begin{figure}[here]
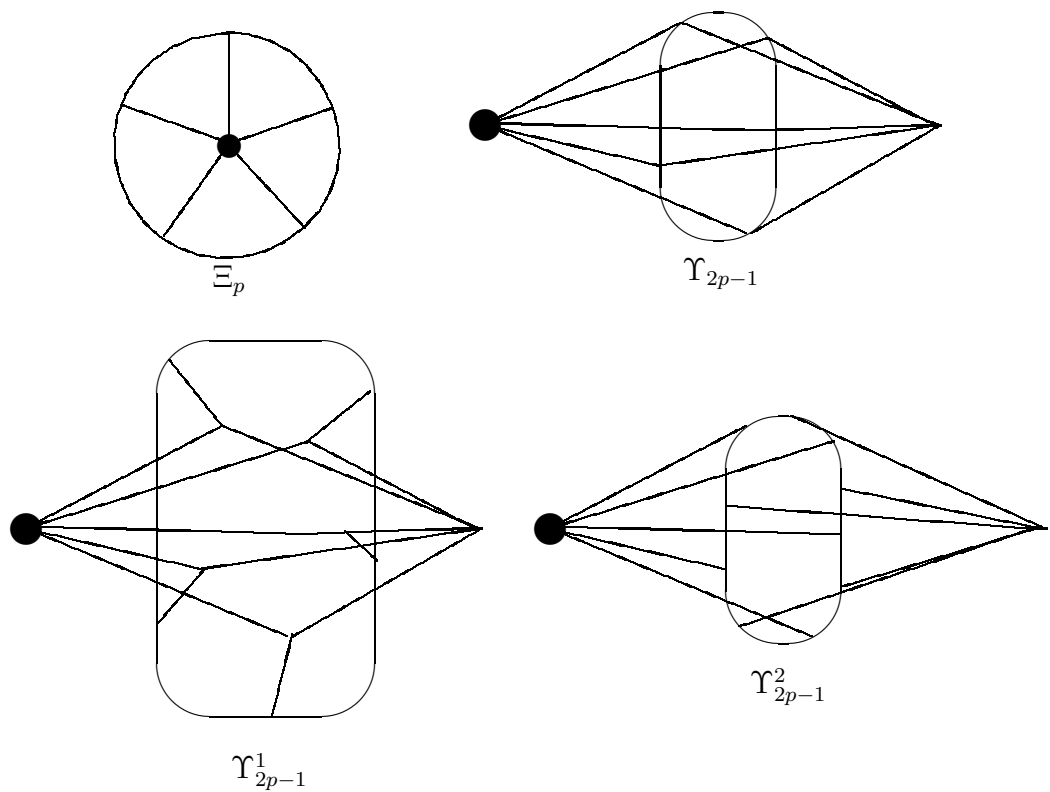

\caption{\label{fig2} Some graphs with $D_{2p}$-symmetry}
\end{figure}

We now define some graphs
that we will need for this chapter.
(Refer to Figure \ref{fig2} for illustrations of most of these graphs.)
Let $\Theta_{p-1}$ and $R_p$ be as before.  That is,
$R_p$ is a rose and
$\Theta_{p-1}$ is the graph with two vertices and $p$ edges, each of
which goes from one vertex to the other.  Let $\Xi_p$ be the
$1$-skeleton of the cone over
a $p$-gon, so that $\Xi_p$ has $p+1$ vertices and $2p$ edges, 
one vertex has valence $p$ and the other $p$ vertices all have valence
$3$.  Let $\Upsilon_{2p-1}$ be the
$1$-skeleton of the 
suspension of a $p$-gon.  Hence
$\Upsilon_{2p-1}$ has $p+2$ vertices and $3p$ edges;  two of the vertices have
valence $p$ and the other $p$ have valence $4$.
Let $\Upsilon_{2p-1}^1$ and $\Upsilon_{2p-1}^2$ be the two possible graphs
that can be obtained from $\Upsilon_{2p-1}$ by
equivariantly blowing
up (see Definition \ref{tr22})
the $p$ valence $4$ vertices into $2p$ valence $3$
vertices.  That is, $\Upsilon_{2p-1}^1$ can be obtained by first taking a
$p$-gon and then attaching $p$ free edges to the $p$ vertices of
the $p$-gon.  Say each of these new edges $e_i$ begins at the vertex
$x_i$ and ends at the vertex $y_i$, and suppose that the vertices
$x_i$ are the ones that are attached to the $p$-gon.
Now form the $1$-skeleton of the
double cone or suspension over the $p$ vertices
$y_i$.  This gives the graph $\Upsilon_{2p-1}^1$.  The graph
$\Upsilon_{2p-1}^2$ can be thought of as follows:  First take a $p$-gon
and
cone off over the $p$ vertices of the $p$-gon.  Now also cone off
over the $p$ midpoints of the $p$ edges of the $p$-gon.
Note that there
is an obvious $\Z_p$-action on each of $\Theta_{p-1}$, $\Xi_p$,
$\Upsilon_{2p-1}$, $\Upsilon_{2p-1}^1$, and $\Upsilon_{2p-1}^2$.

\bigskip

\input{fig2.pic}

\begin{figure}[here]
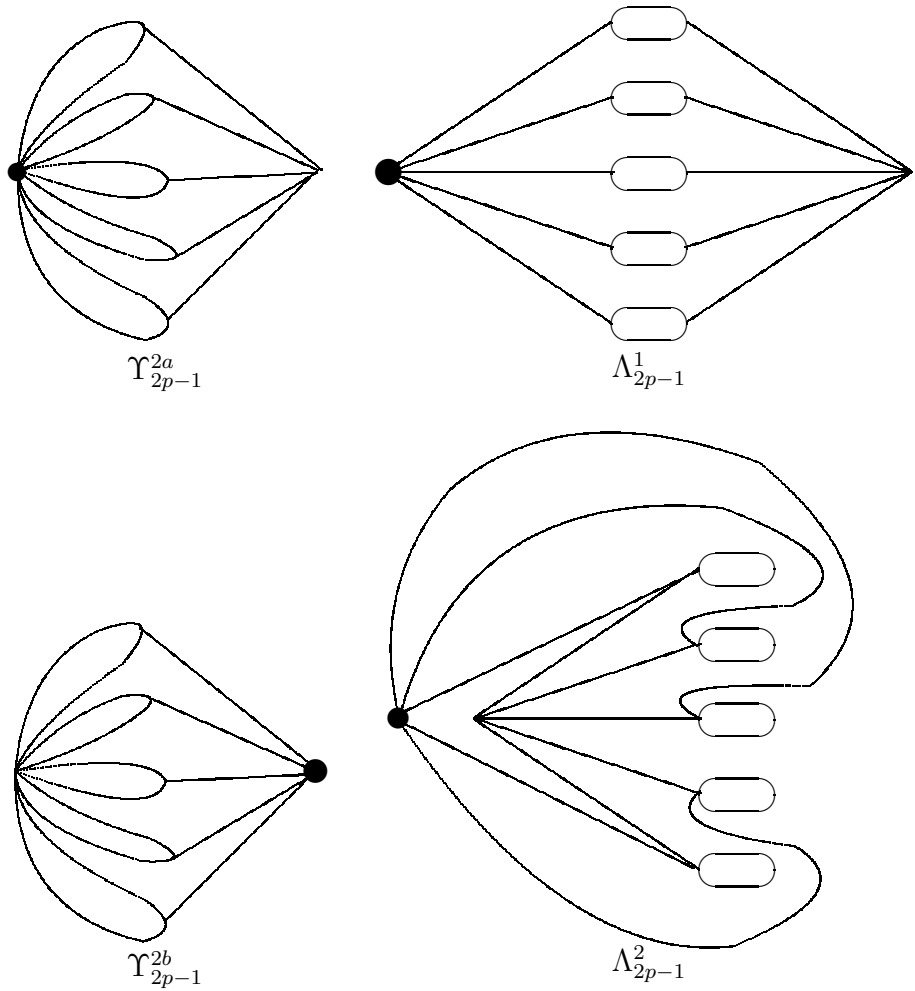

\caption{\label{fil3} Some graphs with (cohomologically) $\Sigma_p$-symmetry}
\end{figure}

Choose basepoints for the graphs $\Theta_{p-1}$, $\Upsilon_{2p-1}$,
$\Upsilon_{2p-1}^1$, and $\Upsilon_{2p-1}^2$ as
illustrated in Figure \ref{fig2}:  
Let the vertex on the ``leftmost'' side of $\Theta_{p-1}$ be
the basepoint.  Additionally, orient $\Upsilon_{2p-1}$,
$\Upsilon_{2p-1}^1$, and $\Upsilon_{2p-1}^2$
so that one of their
valence $p$ vertices is on the ``left'' and the other is on the
``right'' and choose the leftmost vertex to be the basepoint.
Writing $R_p \vee \Theta_{p-1}$ will mean that the two graphs are
wedged together at the basepoint of $\Theta_{p-1}$, while writing
$\Theta_{p-1} \vee R_p$ will mean that the non-basepointed vertex of
$\Theta_{p-1}$ is wedged to the vertex of $R_p$.  Let $\Upsilon_{2p-1}^{2a}$
be the graph obtained from $\Upsilon_{2p-1}^2$ by collapsing the
leftmost $p$ edges and let $\Upsilon_{2p-1}^{2b}$ be the one obtained
by collapsing the rightmost $p$ edges.
Refer to Figure \ref{fil3} for pictures of these graphs.
Figure \ref{fil3} also depicts two graphs
$\Lambda_{2p-1}^1$ and $\Lambda_{2p-1}^2$ which will be
used in the proof of Lemma \ref{t40}.

\bigskip

\newpage

\input{fil3.pic}

For the next lemma (Lemma \ref{t40})
only, we will consider the above graphs not to have
basepoints specified. 
The basepoint is just assumed to be located at some spot which is
invariant under the $\Z/p$ action. (This assumption is for
convenience rather than anything else,
so that we will not need to
introduce several separate subcases, each corresponding to
a different location for the basepoint.)

We will be looking at the standard equivariant spectral
sequence \ref{e1} applied to calculating the cohomology
groups $H^*(Aut(F_n); \Z_{(p)})$.  In particular, we
will be looking at the $E_1$ page of this
spectral sequence only in rows $0$ through $2p-3$, and often
just in rows $1$ though $2p-3$.  One interesting fact about these
rows is that they allow us to distinguish between simplices 
that have stabilizers whose cohomology is $\Z/p$, $D_{2p}$,
or $\Sigma_p$.
It is well
known that
$$\matrix{\hfill H^*(\Z/p; \Z_{(p)}) &=& \Z_{(p)}[x_2]/(px_2), \hfill \cr
\hfill H^*(D_{2p}; \Z_{(p)}) &=& \Z_{(p)}[x_4]/(px_4), \hfill \cr
\hfill \hbox{and } H^*(\Sigma_p; \Z_{(p)}) &=&
\Z_{(p)}[x_{2(p-1)}]/(px_{2(p-1)}), \hfill \cr}$$
where $x_2$, $x_4$ and $x_{2(p-1)}$ are generators of dimensions
$2$, $4$ and
$2(p-1)$, respectively.  Hence if a simplex of $X_n$ has
stabilizer isomorphic 
to $\Z/p$ or $D_{2p}$ then it will
contribute something to the $E_1$ page of the spectral
sequence in some of the rows $1$ through $2p-3$.  On the other hand,
if its stabilizer is isomorphic to $\Sigma_p$, then it
will not contribute anything to the $E_1$ page of the
spectral sequence in the given rows.

Define an $r$-simplex
in the $p$-singular locus of $X_n$ to have
{\em exactly $\Z/p$ symmetry} if it contributes exactly one copy of
$\Z/p$ to each of the entries $E_1^{r,2k}$, $0 < 2k < 2(p-1)$,
in the $E_1$ page of the spectral sequence.  Define an $r$-simplex
in the $p$-singular locus of $X_n$ to have {\em at most dihedral
symmetry}
if it contributes exactly one copy of
$\Z/p$ to each of the entries $E_1^{r,4k}$, $0 < 4k < 2(p-1)$,
in the $E_1$ page of the spectral sequence.

The next lemma examines which vertices in $X_n$ 
contribute to the spectral sequence in the given rows.
The proof of the
lemma actually explicitly enumerates which graphs have $p$-symmetry,
which will be very useful to us later.

\begin{lemma} \label{t40} Let $p \geq 5$ be prime and set $n=2p-1$.
Let the marked graph $(\xi,G,x_0)$ be a vertex in the
$p$-singular locus of $X_n$. 
Then the cohomology
$H^*(Aut(G,x_0); \Z_{(p)})$
of the stabilizer of this vertex is the
same as the cohomology with $\Z_{(p)}$ coefficients
of one of $D_{2p}$, $D_{2p} \times \Sigma_p$,
$\Sigma_p$, $\Sigma_p \times \Sigma_p$, 
$(\Sigma_p \times \Sigma_p) \rtimes \Z/2$
or $\Sigma_{2p}$.
\end{lemma}

\PF The interested reader should refer to chapter 6 of
Adem and Milgram \cite{[A-M]} for information about
the cohomology of the symmetric groups in general.
There they calculate,
for example, $H^*(\Sigma_{2p}; \Z/p)$
and more generally $H^*(\Sigma_{m}; \Z/p)$.  We will not
go into this here, as we are not so much interested
in the structure of $H^t(\Sigma_{2p}; \Z_{(p)})$
as we are interested in the fact that
it is $0$ for $1 \leq t \leq 2p-3$.

From \cite{[B-T]} and \cite{[M]},
we see that $p^2$ is an upper bound for the order of any $p$-subgroup
of $Aut(F_n)$.  Thus $p^2$ is an
upper bound for the order of a maximal $p$-subgroup $P$ of $Aut(G,x_0)$.
Since all possible choices for $P$ are abelian (i.e., there are only
three possibilities: $\Z/p$, $\Z/p^2$, and $\Z/p \times \Z/p$), we can
apply Swan's theorem \ref{tr30}
to see that
\begin{equation} \label{e8}
H^*(Aut(G,x_0); \Z_{(p)}) = H^*(P; \Z_{(p)})^{N_{Aut(G,x_0)}(P)}.
\end{equation}

We now look at each of the individual cases $P=\Z/p, \Z/p^2, \hbox{
and } \Z/p \times \Z/p$.

\begin{enumerate}

\item   We will first examine the case where
$P=\Z/p^2$.  In this case, we have that
$p^2$ edges $e_1, \ldots, e_{p^2}$ in $G$ are rotated around by $P$.
An examination of all possible ways that these edges could be
connected together, keeping in mind that $G$ is admissible, reveals
that this case is impossible.  For example, the first subcase is that
all of the $e_i$ begin and end at the same vertex.  This is not
possible because the fundamental group of $R_{p^2}$ is too large for
it to be a subgraph of $G$. For the next subcase, suppose each edge
goes from some vertex $y_1$ to some other vertex $y_2$.
Then they form a $\Theta_{p^2-1}$ inside $G$, which is also impossible.
In the next subcase, the edges begin at one common vertex $y_0$
and end at $p^2$ distinct vertices $y_1, \ldots, y_{p^2}.$  Since the
graph $G$ is admissible, it has no free edges and every nonbasepointed
vertex has valence at least $3$.  So all of the vertices
$y_1, \ldots, y_{p^2}$ have to connect up in some manner, and in doing
so they will violate the fact that $\pi_1(G)=F_{2p-1}$.  The final
three subcases, in which the $e_i$ either form a $p^2$-gon, have no
common vertices,
or form loops with $p^2$ distinct endpoints,
are similar.  Hence $P$ will never
be $\Z/p^2$ and this case will not occur.

\item Next we will examine the case where
$P=\Z/p \times \Z/p= (\alpha) \times (\beta)$.  The first cyclic
summand $(\alpha)$ must
rotate $p$ edges $e_1, \ldots, e_p$ of $G$.  Without loss of generality,
we may assume that
$(\beta)$ fixes these edges (by setting the new $\beta$ to be
$\beta - \alpha^j$ for some $j$.)  Now $(\beta)$ must rotate $p$ other edges
$f_1, \ldots, f_p$ of $G$.
(Note that we cannot have each $e_i$ forming a loop with
endpoints $y_i$, where $y_1, \ldots, y_p$ are all distinct.
This is because we could then choose $p$
$\alpha$-equivariant paths
from each of these points to the basepoint.  This would give us at
least $p$ different separating edges in the graph,
and hence this cannot occur.
The same argument also holds for the $f_i$ instead of the
$e_i$.)
By doing the sort of case-by-case analysis
that we did in the previous paragraph, we see that $G$ must be one of the
following graphs (listed in increasing order with respect
to the number of vertices):
\begin{itemize}
\item $R_p \vee \Theta_{p-1}$, whose automorphism group has
the same cohomology as $\Sigma_p \times \Sigma_p.$
\item $\Theta_{2p-1}$, whose
automorphism group is $\Sigma_{2p}$.
\item $\Theta_{p-1} \vee \Theta_{p-1}$, plus one additional edge $e$
attached in some manner to the existing vertices.
The automorphism group here will have the same
cohomology as either
$\Sigma_p \times \Sigma_p$ or 
$(\Sigma_p \times \Sigma_p) \rtimes \Z/2$.
\item $\Theta_{p-1} \coprod \Theta_{p-1}$,
with one additional
edge $e_1$ attached going from an already existing vertex
of one of the $\Theta$-graphs to one on the other $\Theta$-graph,
after which we sequentially attach another
edge $e_2$ to that resulting graph.
The endpoints of $e_2$ can be attached to any of the
already existing vertices, or they can be attached anywhere in
the interior of $e_1$. 
The automorphism group here will have the same
cohomology as either
$\Sigma_p \times \Sigma_p$ or 
$(\Sigma_p \times \Sigma_p) \rtimes \Z/2$.
\item $\Theta_{p-1} \vee \Xi_p$, with
automorphism group $\Sigma_p \times D_{2p}$.
\end{itemize}

\item  For the final case, $P=\Z/p=(\alpha)$.
The idea behind our proof of
part 3. is that we want to show that all of the $p$-symmetries in
the graph $G$ are also at least $D_{2p}$-symmetries.  That is, in addition
to the rotation by $\Z/p$, there is also a
dihedral ``flip''.
We will be
able to get this result because $n=2p-1$ is not ``large enough'' with
respect to $p$ for us to be able to generate graphs $G$ with
$\pi_1(G)=F_n$ that have $\Z/p$-symmetries but not
$D_{2p}$-symmetries.  

We have $P = \Z/p$
acting on a graph whose fundamental group has
rank $n=2p-1$.
As before, there exist at least $p$
edges $e_1, \ldots, e_p$ that $P$ rotates.  If these edges form a
$\Theta_{p-1}$ or an $R_p$, we are done.  This is because now
$P$ cannot move any other edges of $G$, else we are in the case of the
previous paragraph where  $P=\Z/p \times \Z/p$.
As the automorphism groups of both $\Theta_{p-1}$
and $R_p$ contain the symmetric group $\Sigma_p$,
we are done in this subcase.

Now suppose we are in the other extreme subcase, the one where
$e_1, \ldots, e_p$ have no endpoints in common.  Choose a minimal path
$\gamma_1$ from $e_1$ to the basepoint.  Since $P=(\alpha)$ fixes the
basepoint, we have that $\gamma_i := \alpha^{i-1} \gamma_1$ is a
minimal path from $e_i$ to $x_0$ for all $i=1, \ldots, p$.  Since none
of the endpoints of the $e_i$ can be the basepoint, there are at least
$p$ distinct edges $f_1, \ldots, f_p$ in $\gamma_1, \ldots, \gamma_p$,
respectively.  We can also assume that each $f_i$ has at least one
endpoint that is not the basepoint. Because $G$ has no free edges,
no separating edges, 
and all non-basepointed vertices have valence at least three, another
case-by-case analysis reveals that
since $\pi_1(G)$ must have rank less than $2p$, the graph is forced
to be either $\Upsilon_{2p-1}^1$, $\Upsilon_{2p-1}^2$,
or $\Lambda_{2p-1}^1$.  The first two of these graphs have
dihedral symmetry, while the last has
automorphism group with the same cohomology as $\Sigma_p$.

The next case is the one in
which the $e_i$ are all loops with $p$ distinct
endpoints $y_i$.  As in the previous case, we can choose
a $\Z/p$ equivariant path from each $y_i$ to the
basepoint.  The admissibility conditions on the graph
only allow one possibility, namely the graph $\Lambda_{2p-1}^2$.
As this graph has automorphism group with the same
cohomology as $\Sigma_p$, we are finished with this case.

Next consider the case where the $e_i$ form a $p$-gon.  Since the
vertices of $G$ have valence at least $3$, there must be $p$ other
edges $f_1, \ldots, f_p$ in $G$ that each start at one of the $p$
vertices of the $p$-gon.  Since $G$ is admissible and the rank of
$\pi_1(G)$ is $2p-1$, these additional edges cannot also join up to
form a $p$-gon.  (Why?  Both $p$-gons still need to connect up to the
basepoint in some way, and in connecting up to $x_0$ the rank of the
fundamental group of $G$ will be forced too high.)  In addition, the
edges must have some vertices in common, else we reduce to the
previous case; therefore, the $f_i$ are all forced to end at some common
vertex $y_0$.  In other words, we have a $\Xi_p$ embedded in $G$.  If
$P$ doesn't move any other edges in $G$, we are done since $\Xi_p$
has dihedral symmetry.  If some other edges $g_1, \ldots, g_p$
are moved, they must also
be attached to the $p$-gon that the $e_i$ form, or we will have two
independent $\Z/p$-actions and be in the case $P=\Z/p \times \Z/p$.
None of the following cases can happen, else
rank$(\pi_1(G)) \geq 2p$:
\begin{itemize}
\item The other endpoints of the $g_i$ all connect to $y_0$.
\item The other endpoints of the $g_i$ also connect to the $p$-gon
formed by the $e_i$.
\item The other endpoints of the $g_i$ form $p$ other distinct
vertices.
\end{itemize}
Hence these other endpoints all have to connect to some other common
vertex $y_1$, forming another copy of $\Xi_p$ in $G$.  Thus $G$ must
be the graph $\Upsilon_{2p-1}$, which certainly has dihedral symmetry.

For the final case, the edges $e_i$ have one common vertex $y_0$,
and end in $p$ other distinct vertices $y_1, \ldots, y_p$.
In addition, i) there are no $p$-gons in $G$, ii) there
are no collections
of $p$ edges in $G$ that are rotated by $P$ and that have no common
vertices,
and iii) there
are no collections
of $p$ edges in $G$ that are rotated by $P$ and that
each form loops with distinct endpoints.
Since all of $y_1, \ldots, y_p$ have valence three, $P$
must rotate two other collections of edges
$\{f_1, \ldots, f_p\}$ and $\{g_1, \ldots, g_p\}$ that begin at the
vertices $y_1, \ldots, y_p$ and end at the vertices $z_1$ and
$z_2$, respectively.  Note that since $\pi_1(G)=F_{2p-1}$,
$|\{y_0,z_1,z_2\}| \geq 2$.  Also note that $P$ cannot move any other
edges of $G$ except the ones we have listed.  In this case, the
symmetric group $\Sigma_p$ acts on the collections of edges defined above,
and so the cohomology of the group of graph automorphisms of
the graph is the same as that of the symmetric group.
If $|\{y_0,z_1,z_2\}| = 2$ then the only edges in the
graph are the $e_i$, $f_i$, and $g_i$ and the graph is
either $\Upsilon_{2p-1}^{2a}$ or $\Upsilon_{2p-1}^{2b}$.
On the other hand, if $|\{y_0,z_1,z_2\}| = 3$, then
the graph has one additional edge besides the
$e_i$, $f_i$, or $g_i$.
Accordingly, the graph looks like
a $\Phi_{2p-1}$ (See Figure \ref{fig1}) with
one additional edge added.  This additional edge can go
from any of the $\{y_0,z_1,z_2\}$ to any other one,
including possibly the
same one.  In any case, it is definitely true that the
graph has automorphism group with the same
cohomology as $\Sigma_p$.
The lemma follows. \END
\end{enumerate}

For an example of what we were trying to avoid in the proof
of the above lemma,
refer to the three examples given in the Figure \ref{fil1} below.

\input{fil1.pic}
\begin{figure}[here]
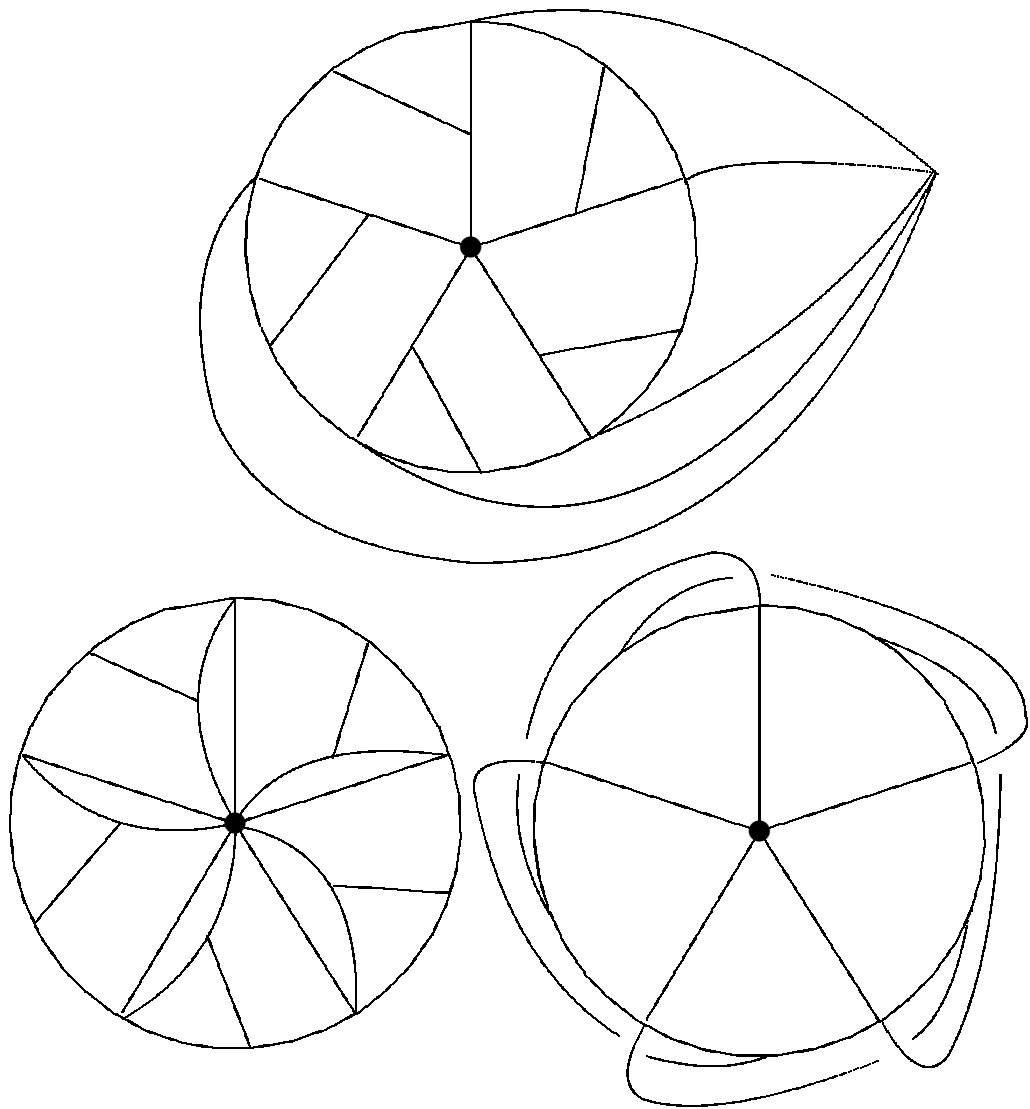

\caption{\label{fil1} Graphs whose symmetry groups are exactly $\Z/p$}
\end{figure}

The graphs pictured above have an obvious $\Z/p$-symmetry
given by rotation about the basepoint, which is
indicated by a solid  dot.  But they have no dihedral
flip, and their basepoint-preserving automorphism groups
are all exactly $\Z/p$,
where $p=5$ in the examples pictured and where obvious
analogues exist for other odd primes.  The ranks of the
fundamental groups of the graphs pictured are
$3p-1$, $3p$, and $2p$, respectively.  The last rank, $2p$,
is the lowest rank possible where one can have a graph with
exactly $\Z/p$ symmetry.  (Note that one can always construct
graphs with exactly $\Z/p$ symmetry, but higher rank
fundamental groups, from any of the above graphs, by adding
concentrically embedded
circles around the basepoint of the graph,
where the circles are attached to the ``spokes'' radiating
away from the basepoint.)

\begin{cor} \label{ttt2}
A vertex in the $p$-singular locus of $X_n$ has at most dihedral symmetry 
if its cohomology is the same as that of $D_{2p}$ or $D_{2p} \times \Sigma_p$,
and vertices in the $p$-singular locus
will never have exactly $\Z/p$ symmetry.
\end{cor}

By analyzing the $\Z/p$-invariant subforests of all of the graphs
explicitly listed in the proof of Lemma \ref{t40}, we can see what
types of stabilizers higher dimensional simplices
(rather than just vertices) have.

We will show that the
simplices with at most dihedral symmetry will fall into
two (exhaustive but not disjoint) categories.  The first
category
consists of those that are listed in
Figure \ref{fig4} below.
The second category consists of simplices whose maximal vertex
(recall that $X_n$ is the realization of a poset)
has the form $\Xi_p \vee \Gamma_{p-1}$ 
where $\Gamma_{p-1}$ is some basepointed
graph with fundamental group of rank $p-1$, the wedge does
not necessarily take place at the basepoint, and where the
forest collapses of the simplex respect the $\Z/p$ action
on $\Xi_p$.

We will also show that the simplices listed in
Figure \ref{fil2} are the only ones with exactly $\Z/p$
symmetry.

In the figures below, a dotted line or a hollow dot indicates
that the given edge or vertex, respectively, does not have
the indicated property.  A solid dot, a solid line, or
a $2$-simplex with an X in it, means that the given
vertex, edge, or $2$-simplex, respectively, does have
the indicated property.

\bigskip

\bigskip

\newpage

\input{fig4.pic}
\begin{figure}[here]
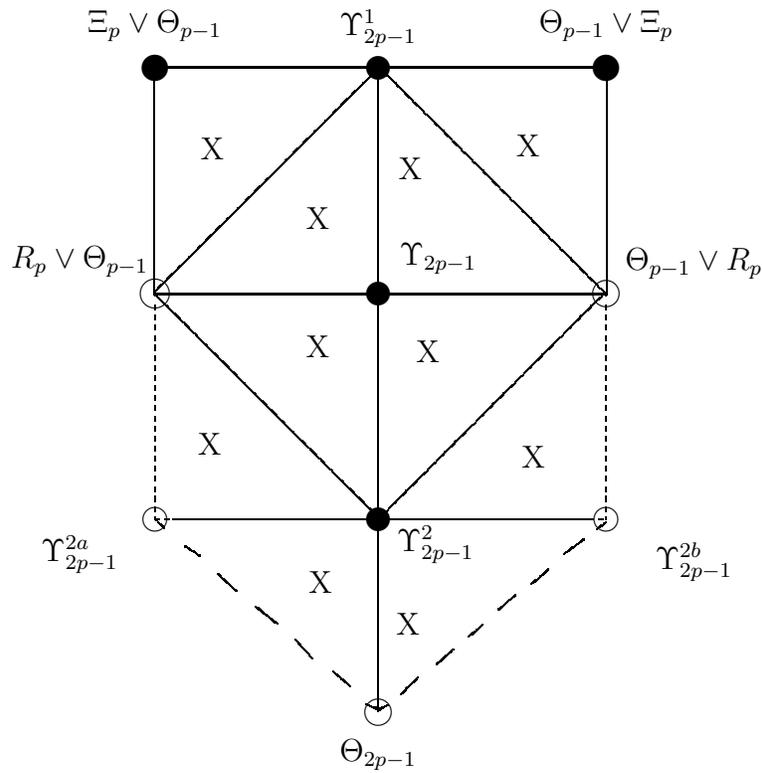

\caption{\label{fig4} Some simplices with at most dihedral symmetry}
\end{figure}

\bigskip

\input{fil2.pic}
\begin{figure}[here]
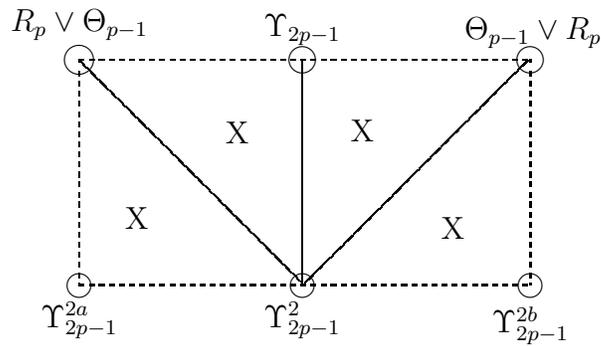

\caption{\label{fil2} Simplices with exactly $\Z/p$ symmetry}
\end{figure}

\bigskip

\begin{cor} \label{t41} Let $p \geq 5$ be prime, $n = 2p-1$,
and consider the $p$-singular locus of the spine $X_n$ of
auter space.
\begin{itemize}
\item
The only simplices
with at most dihedral symmetry are either: (i) listed in
Figure \ref{fig4};  or (ii) have maximal vertex
of the form $\Xi_p \vee \Gamma_{p-1}$.
\item The only simplices with exactly
$\Z/p$ symmetry are those listed in Figure \ref{fil2}.
\end{itemize}
\end{cor}

\PF 
We examine each of the graphs listed in
Lemma \ref{t40} separately.  By enumerating the $\Z/p$ invariant
subforests of each of these graphs, one can list all of the simplices
in the $p$-singular locus of $X_n$.  We can ignore the graphs
in Lemma \ref{t40} that do not have dihedral symmetry, as all of
their symmetry comes from symmetric groups.  When you collapse
invariant subforests of these graphs, you still get
graphs with symmetry coming from the symmetric group.

So we are left with analyzing the graphs from Lemma \ref{t40} with
dihedral symmetry, which were:
\begin{itemize} 
\item $\Theta_{p-1} \vee \Xi_p$. (There are actually two
possibilities here as the enumeration in Lemma \ref{t40}
did not specify basepoints.  The central vertex of $\Xi_p$
could be attached to either the basepoint of $\Theta_{p-1}$
or the other vertex of $\Theta_{p-1}$.)
\item $\Gamma_{p-1} \vee \Xi_p$, where $\Gamma_{p-1}$
is a basepointed graph with fundamental group of rank $p-1$
which has
no $p$-symmetry (or where $\Gamma_{p-1}$ is $\Theta_{p-1}$
but the central vertex of $\Xi_p$ is attached to the
midpoint of an edge of $\Gamma_{p-1}$.)
\item $\Upsilon_{2p-1}$.
\item $\Upsilon_{2p-1}^1$.
\item $\Upsilon_{2p-1}^2$.
\end{itemize}

For the first two types of graphs, you can obtain
simplices with dihedral symmetry by collapsing all of the spokes
of $\Xi_p$ and/or any forest in the other graph
of the wedge sum
(either $\Theta_{p-1}$ or $\Gamma_{p-1}$.)  The resulting
simplex with maximal vertex $\Theta_{p-1} \vee \Xi_p$
or $\Gamma_{p-1} \vee \Xi_p$ will
clearly have at most dihedral symmetry, and 
will also just as clearly not give you a graph with exactly $\Z/p$
symmetry.

In a similar manner, simplices in the $p$-singular locus
of $X_n$ with maximal vertex $\Upsilon_{2p-1}$ or
$\Upsilon_{2p-1}^1$ are (exhaustively) listed in
Figure \ref{fig4}.  Note that $\Upsilon_{2p-1}$
can only be blown up (while still preserving the $\Z/p$
action so that we stay in the $p$-singular locus of $X_n$)
in two ways, to either
$\Upsilon_{2p-1}^1$ or $\Upsilon_{2p-1}^2$. The latter two
graphs cannot be blown up at all.

Finally, the simplices with maximal vertex  $\Upsilon_{2p-1}^2$
are listed in Figure \ref{fig4} or Figure \ref{fil2}.
Note that we can obtain edges and $2$-simplices with
exactly $\Z/p$ symmetry, even though no actual vertex
of $X_n$ has exactly $\Z/p$ symmetry.  This is because
you can choose subforests of $\Upsilon_{2p-1}^2$
that do respect the dihedral ``flip'' of
$\Upsilon_{2p-1}^2$.  In other words, this flip will
not take the subforest to itself again.  Hence the
resulting simplex will just have symmetry group $\Z/p$.
Last of all, note that you can also choose subforests
of $\Upsilon_{2p-1}^2$ which do respect the dihedral
flip, and these give simplices with dihedral symmetry. \END

\chapter{The integral cohomology of the quotient never stabilizes} \label{c14}

The main result of this chapter will be to show that
for each $k \in \Z^+$, either
$$H^{4k}(Aut(F_\infty); \Q) \not = 0$$
or
$$H^{4k+1}(Aut(F_m); \Z) \hbox{ never stabilizes as } m \to \infty.$$
In particular, since $H^{4}(Aut(F_\infty); \Q) = 0$ by
\cite{[V]}, this will prove that $H^{5}(Aut(F_m); \Z)$
never stabilizes.
As in Chapter \ref{c13}, all primes $p$
considered are assumed to be greater than or equal to $5$.

\begin{lemma} \label{t42}
For the rows $0 \leq s < 2(p-1)$, the $E_2$ page of
the spectral sequence \ref{e1}
applied to calculate $H^*(Aut(F_{2p-1}); \Z_{(p)})$
is given by
\bigskip
$$\matrix{
\hfill E_2^{r,s} &=& &\left\{\matrix{
H^r(Q_{2p-1}; \Z_{(p)})  \hfill &s = 0 \hfill \cr
\Z/p \hfill &r=2 \hbox{ and } s = 4k-2 > 0, k \in \Z^+ \hfill \cr
0 \hfill &\hbox{otherwise} \hfill \cr} \right. \hfill \cr
}$$
\end{lemma}
\bigskip

\PF As $E_1^{r,0}$ is the cochain complex $C^r(Q_{2p-1}; \Z_{(p)})$,
it follows that $E_2^{r,0} = H^r(Q_{2p-1}; \Z_{(p)})$ as
claimed above.

None of the simplices in $X_{2p-1}$ contribute anything to the
odd rows between $0$ and $2(p-1)$ of the above spectral sequence,
from Corollary \ref{t41}.  Also from Corollary \ref{t41}, the ones that contribute
to rows of the form $4k - 2$,
$k \in \Z^+$, are all listed in Figure \ref{fil2}.
Let $A$ be the subcomplex of $Q_n$ generated by all of the
simplices pictured in Figure \ref{fil2} and let $B$ be the
subcomplex generated by just the simplices corresponding to
dotted lines or hollow dots in Figure \ref{fil2}.  Then the row
$s=4k-2$ on the $E_1$ page of the spectral sequence is
$C^r(A,B; \Z/p)$.  Examining Figure \ref{fil2} we see that
$$\matrix{
\hfill H^r(A,B; \Z/p) &=& &\left\{\matrix{
\Z/p \hfill &r=2 \hfill \cr
0 \hfill &\hbox{otherwise} \hfill \cr} \right. \hfill \cr
}$$
Consequently the $E_2$ page is as claimed for the
rows $s = 4k-2$.

Our final task is to calculate the $E_2$ page for the
rows $s=4k$.  Simplices in the $p$-singular locus
of the spine with ``at most dihedral symmetry'' contribute
to these rows.  From Corollary \ref{t41}, we have a characterization
of such simplices. 
Define the subcomplex $M$ of the
$p$-singular locus of the spine
$X_{2p-1}$ of auter space to be the
subcomplex generated by
simplices with ``at most dihedral symmetry''.
More precisely, from Corollary \ref{t41}, we know it is
generated by the simplices corresponding
to those in Figure \ref{fig4}
(i.e., corresponding in the sense that we are
taking $M$ to be a subcomplex of the spine
rather than its quotient and Figure \ref{fig4} is
a picture in the quotient)
in addition to simplices whose maximal vertex
has underlying graph of the
form $\Xi_p \vee \Gamma_{p-1}$
(where the
forest collapses in the simplices respect
the $\Z/p$ action on $\Xi_p$.)
Recall that an $r$-simplex 
with at most dihedral symmetry
contributes exactly one $\Z/p$ to 
$E_1^{r,4k}$, while all other
simplices (those without dihedral or
exactly $\Z/p$ symmetry)
contribute nothing to this row.
 
Let $N$ be the full
subcomplex of $M$ generated by simplices in $M$
which do not have at most dihedral symmetry.
Observe that none of the simplices in $N$ have
at most dihedral symmetry.
Also note that the row $E_1^{*,4k}$ is the
relative cochain complex
$$C^*(M/Aut(F_{2p-1}),N/Aut(F_{2p-1}); \Z/p).$$

Let $M'$ be the subcomplex of $M$ generated by $N$ and by
simplices whose maximal vertex is 
$\Upsilon_{2p-1}^2$.  Hence $M'$ is the subcomplex
consisting of $N$ and the bottom two thirds
of Figure $\ref{fig4}$.  There is an 
$Aut(F_{2p-1})$-equivariant deformation retraction
of $M$ onto $M'$, given on the vertices of the
poset by:
\begin{itemize}
\item Contracting the spokes of the
graph $\Xi_p$ in $\Theta_{p-1} \vee \Xi_p$.
\item Contracting the spokes of the graph $\Xi_p$ in
$\Gamma_{p-1} \vee \Xi_p$, where $\Gamma_{p-1}$ has
no $p$-symmetry.
\item Contracting the $p$ outward radiating
edges attached to the $p$-gon in the center of
the graph $\Upsilon_{2p-1}^1$.  In the terminology
used at the beginning of Chapter \ref{c13} while defining
$\Upsilon_{2p-1}^1$, we are contracting the edges $e_i$.
\end{itemize}
That it is a deformation retraction follows from the
Poset Lemma in \cite{[K-V]} attributed to Quillen.

In other words, the deformation retraction of $M$ to
$M'$ comes from collapsing the ``inessential edges
of the $D_{2p}$-action'' (see part \ref{p3} and
\cite{[K-V]}.)  Two words of caution are in order here
with respect to this use of the word ``inessential''.
First, note that we are not collapsing \underline{all} of the 
inessential edges with respect to some $D_{2p}$
action, just those edges that are associated with
some $p$-gon giving dihedral symmetry.  Hence, for example,
we are not collapsing any of the $2p$ inessential
edges of the graph $\Lambda_{2p-1}^1$
with respect to some $D_{2p}$ action on that graph.
Second,
do not be fooled into thinking that there is only one
conjugacy class of subgroups of $Aut(F_{2p-1})$ that
is providing all of these $D_{2p}$-symmetries for $M$.
In other words,
our calculations are not taking place in the fixed
point subcomplex $X_{2p-1}^G$ for just one particular
$D_{2p} \cong G \subset Aut(F_{2p-1})$.  For example,
it can be seen from Theorem \ref{tr4}
that the dihedral groups corresponding to
$\Upsilon_{2p-1}^1$ and  $\Xi_p \vee R_{p-1}$ are
not conjugate to each other.

As the homotopy retracting $M$ to $M'$ is
$Aut(F_{2p-1})$-invariant, it descends to a deformation
retraction of $M/Aut(F_{2p-1})$
to $M'/Aut(F_{2p-1})$.  Hence the relative cohomology
groups
$$H^*(M/Aut(F_{2p-1}),N/Aut(F_{2p-1}); \Z/p)$$
and
$$H^*(M'/Aut(F_{2p-1}),N/Aut(F_{2p-1}); \Z/p)$$
are isomorphic.  Now referring to
Figure $\ref{fig4}$,
we see that
$$H^t(M'/Aut(F_{2p-1}),N/Aut(F_{2p-1}); \Z/p) = 0$$
for all $t$
because we can contract all of the
simplices in $M'/Aut(F_{2p-1})$
uniformly into $N/Aut(F_{2p-1})$. \END

An immediate consequence is

\begin{remark} \label{t43}
Let $k$ be such that $0 < 4k - 2 < 2(p-1)$.
In the spectral sequence \ref{e1}
applied to calculate $H^*(Aut(F_{2p-1}); \Z_{(p)})$,
a nontrivial torsion cohomology class in $E_1^{2,4k-2}$
survives at least until it would be mapped to the $r$-axis;
that is, it survives at least up to the
$E_{4k-1}$-page, after which the
differential
$$d_{4k-1} : E_{4k-1}^{2,4k-2} \cong \Z/p
\to E_{4k-1}^{4k+1,0} \cong H^{4k+1}(Q_{2p-1}; \Z_{(p)})$$
might kill the class by mapping it to
something nonzero.
\end{remark}

\begin{thm} \label{t44} For all positive integers
$k$, either
$$H^{4k}(Aut(F_\infty); \Q) \not =0$$
or
$$H^{4k+1}(Q_m; \Z)
\hbox{ never stabilizes as } m \to \infty.$$
\end{thm}

\PF From \cite{[H-V]}, if $m \geq 8k+3,$ then the standard map
$$H^{4k}(Aut(F_{m+1}); \Z) \to H^{4k}(Aut(F_m); \Z)$$
is an isomorphism.
Observe that $H^{4k}(Aut(F_{8k+3}); \Z) =
H^{4k}(Aut(F_\infty); \Z)$ is a finitely generated
abelian group.  If it contains a torsion free summand isomorphic to
$\Z$, then we are done and $H^{4k}(Aut(F_\infty); \Q) \not =0$.
Otherwise, choose a prime $q$ such that $2q-1 \geq 8k+3$
and so that for all primes $p \geq q$ there is no $p$-torsion
in $H^{4k}(Aut(F_{8k+3}); \Z)$.
We will show that $H^{4k+1}(Q_{2p-1}; \Z)$  has $p$-torsion for all
primes $p \geq q$, which will prove the theorem.

Let $p \geq q$.  From the lemma above, if we
use the standard equivariant spectral sequence to calculate
$H^{*}(Aut(F_{2p-1}); \Z_{(p)})$, then a class
$\alpha \in E_1^{2,4k-2}$ in the $E_1$-page survives at least
until the $E_{4k-1}$-page.
Because $H^{4k}(Aut(F_{2p-1}); \Z)$ has no
$p$-torsion and $H^{4k}(Aut(F_{2p-1}); \Q) = 0$, 
we have 
$H^{4k}(Aut(F_{2p-1}); \Z_{(p)})=0$.  
Hence the class $\alpha \in E_1^{2,4k-2}$ cannot survive to the
$E_\infty$ page. 

It follows that there is $p$-torsion in
$E_{4k-1}^{4k+1,0}$. 
Recall that
$E_1^{r,0}$ corresponds to the cellular
chain complex with $\Z_{(p)}$ coefficients
for $Q_{2p-1}$.  The $p$-torsion in
$E_{4k-1}^{4k+1,0}$, therefore, would have to have been created
when going from the $E_1$ to $E_2$ pages, because any of
the torsion above the horizontal axis of the spectral sequence
could not map onto a torsion free element on the horizontal axis.
So $H^{4k+1}(Q_{2p-1}; \Z_{(p)})$ has $p$-torsion, and thus
$H^{4k+1}(Q_{2p-1}; \Z)$ has $p$-torsion. \END

From calculations in \cite{[V]} that
$H^4(Aut(F_\infty);\Q)=0$, it follows that

\begin{cor} \label{t45} The cohomology group $H^5(Q_m; \Z)$ never
stabilizes as $m \to \infty$.
\end{cor}

\chapter{Graphs without basepoints} \label{c15}

Let $\hat X_n$ be the spine of Culler-Vogtmann space
\cite{[C-V]}
for $Out(F_n)$, a moduli space of marked graphs.
Let $\hat Q_n$, a moduli space of graphs,
be the quotient of $\hat X_n$ by $Out(F_n)$.

For the case of graphs without basepoints,
substantially the same arguments can be made as
those in the previous two chapters.  For example,
we can state exact analogues of
Theorem \ref{t44} and Corollary \ref{t45}:

\begin{thm} \label{t46} For all positive integers
$k$, either
$$H^{4k}(Out(F_\infty); \Q) \not =0$$
or
$$H^{4k+1}(\hat Q_m; \Z)
\hbox{ never stabilizes as } m \to \infty.$$
\end{thm}

\PF From \cite{[H]},
$$H^{4k}(Out(F_\infty); \Q) = H^{4k}(Aut(F_\infty); \Q)$$
and
if $m \geq 4k^2+10k+1,$ then the standard map
$$H^{4k}(Aut(F_{m}); \Z) \to H^{4k}(Out(F_m); \Z)$$
is an isomorphism.
Observe that $H^{4k}(Out(F_{4k^2+10k+1}); \Z) =
H^{4k}(Out(F_\infty); \Z)$ is a finitely generated
abelian group.  If it contains a torsion free summand isomorphic to
$\Z$, then we are done and $H^{4k}(Out(F_\infty); \Q) \not =0$.
Otherwise, choose a prime $q$ such that $q + 1 \geq 4k^2+10k+1$
and so that for all primes $p \geq q$ there is no $p$-torsion
in $H^{4k}(Out(F_{4k^2+10k+1}); \Z)$.
We will show that $H^{4k+1}(\hat Q_{p+1}; \Z)$ has $p$-torsion for
infinitely many primes $p$,
namely those primes
$p \geq q$ which are congruent to $3$
modulo $4$, which will prove the theorem.

Let $p \geq q$ with $p \equiv 3 \hbox{ (mod } 4 \hbox{)}$.
Because $H^{4k}(Out(F_{p+1}); \Z)$ has no
$p$-torsion, there is also no $p$-torsion in
$H^{4k}(Out(F_{p+1}); \Z_{(p)})$.
From the calculation of Glover and Mislin
in \cite{[G-M]} of the $E_2$-page of the
equivariant spectral sequence used to
calculate $H^{*}(Out(F_{2p-1}); \Z_{(p)})$,
we know that this $E_2$-page,
in the rows $0 \leq s < 2(p-1)$,
is given by
\bigskip
$$\matrix{
\hfill E_2^{r,s} &=& &\left\{\matrix{
H^r(\hat Q_{p+1}; \Z_{(p)})  \hfill &s = 0 \hfill \cr
\Z/p \hfill &r=0 \hbox{ and } s = 4k > 0, k \in \Z^+ \hfill \cr
(n_p) \Z/p \hfill &r=1 \hbox{ and } s = 4k > 0, k \in \Z^+ \hfill \cr
0 \hfill &\hbox{otherwise} \hfill \cr} \right. \hfill \cr
}$$
\bigskip
where $n_p = (p-1)/12 - \epsilon_p$
and $\epsilon_p \in \{0,1\}$.

Hence a class
$\hat \alpha \in E_2^{0,4k}$ in the $E_2$-page survives at least
until the $E_{4k+1}$-page.
The class $\hat \alpha \in E_2^{0,4k}$ cannot survive to the
$E_\infty$ page, however, because there is no $p$-torsion in
the finite (since $H^{4k}(Out(F_{p+1}); \Q) = 0$)
additive group $H^{4k}(Out(F_{p+1}); \Z_{(p)})$.

It follows that there is $p$-torsion in
$$E_{4k+1}^{4k+1,0} = H^{4k+1}(\hat Q_{p+1}; \Z_{(p)}).$$
Thus $H^{4k+1}(\hat Q_{p+1}; \Z)$ has $p$-torsion. \END

From calculations in \cite{[V]} that
$H^4(Out(F_\infty);\Q)=0$, it follows that

\begin{cor} \label{t47} The cohomology group $H^5(\hat Q_m; \Z)$ never
stabilizes as $m \to \infty$.
\end{cor}

Note that the behavior illustrated in Theorem \ref{t46}
appears to be a fairly common phenomenon, and that
additional theorems could be stated about similar
phenomena.  For example, from looking at the
terms $E_2^{1,4k}$ in the spectral sequence
Glover and Mislin \cite{[G-M]} used to calculate
$H^*(Out(F_{p+1}); \Z_{(p)})$, one could state
an exact analogue of Theorem \ref{t46} which implies that
either
$$H^{4k+1}(Out(F_\infty); \Q) \not =0$$
or
$$H^{4k+2}(\hat Q_m; \Z)
\hbox{ never stabilizes as } m \to \infty.$$
From this in turn one could use Vogtmann's
calculation that
$H^{5}(Out(F_\infty); \Q)=0$
to get an analogue of Corollary \ref{t47}
about $H^6(\hat Q_m; \Z)$ never stabilizing.

\appendix

\chapter{On the prime $p=3$} \label{appen}

For the sake of having a concrete example, we
calculate the cohomologies of all of the
quotient spaces involved in Theorem \ref{t9} when $p=3$.
This was also done, independently, by Glover and
Henn.

\begin{thm} \label{tapp}
\medskip

$$\hat H^t(Aut(F_{4}); \Z_{(3)}) = 
\left\{\matrix{
\Z/9 \oplus 3(\Z/3) \hfill &t = 0 \hfill \cr
([{3k \over 2}] + 2)\Z/3 \hfill &|t| = 4k \not = 0 \hfill \cr
\Z/3 \hfill &t = 1 \hfill \cr
0 \hfill &|t| \equiv 1,2  \hbox{ } (\hbox{mod } 4)
\hbox{ and } t \not = 1 \hfill \cr
([{3k \over 2}] - 1)\Z/3 \hfill &|t| = 4k-1 \hfill \cr
} \right.$$

\medskip
\end{thm}

\PF 
Examining Theorem \ref{t9} reveals that we
must show that
none of the various quotient spaces
$H^r(\tilde Q_{k}; \Z/3)$, $H^r(Q_{k}; \Z/3)$,
and $H^r(Q_{2}^\omega; \Z/3)$ (where $k = 1,2$)
contribute any nonzero cohomology classes.

The groups $H^r(Q_{1}; \Z/3)$, $H^r(Q_{2}; \Z/3)$,
$H^r(\tilde Q_{1}; \Z/3)$, and $H^1(\tilde Q_{2}; \Z/3)$
are all zero by Lemma \ref{tr32},
Lemma \ref{tr33}, and Proposition \ref{tr36}.
In addition, 
Fact \ref{tr31}
gives us that $H^2(\tilde Q_{2}; \Z/3)=0$.

For $H^r(Q_{3}^\omega; \Z/3)$
(see Definition \ref{tr17}
to recall the definitions related to this space), 
we note that the relevant
marked graphs in 
$Q_{3}^{\langle \omega \rangle} \cong X_4^A / N_{Aut(F_4)}(A)$
are those listed in component $(\mathcal{A})$
of section 4 in the paper \cite{[G-M]} by Glover and Mislin,
with the additional complication that a basepoint $*$ can
be added to the graphs in various places.  However, most of these
graphs have inessential edges (see \cite{[K-V]}) under the action
of $N_{Aut(F_4)}(A)$, and are thus collapsed directly
away when we reduce from  
$Q_{3}^{\langle \omega \rangle} \cong X_4^A / N_{Aut(F_4)}(A)$
to the space $Q_{3}^\omega$.
We list the graphs
from \cite{[G-M]} that give
$X_4^A / N_{Aut(F_4)}(A)$
here:
\begin{itemize}
\item $R_4$.  The rose has no inessential edges, even when you
attach the basepoint $*$ to the middle of one of the petals.
\item $\Theta_4$.  This $\Theta$-graph has no inessential
edges, regardless of where the basepoint is attached.
\item $W_3 \vee R_1$.  In our notation, this would be
$\Xi_3 \vee R_1$. The inessential edges in the ``spokes'' of the
graph $W_3$ are collapsed and reduce this graph to $R_4$.
\item $\Theta_3 * R_1$.  This is $\Theta_3$ with a loop $R_1$
attached to the middle of one of the edges of the $\Theta$-graph.
The edge of the $\Theta$-graph that the loop is attached
to is inessential.  (It will be 2 or 3 actual edges in the
resulting graph, all of which are inessential, depending upon
where the basepoint is placed.)  Collapsing the
inessential edges yields $R_4$.
\item $\Theta_2 \lozenge Y$.  A graph in the shape of a
letter $Y$ attached to the graph $\Theta_2$, with the
top vertices of the $Y$ attached to one side of the
$\Theta$-graph, and the bottom vertex of the $Y$ attached
to the other side of the $\Theta$-graph.  The bottom edge of
the $Y$ is inessential.  Collapsing this gives $\Theta_4$.
\item $\Theta_2 ** \Theta_1$.  Two $\Theta$-graphs with
a line drawn from the left vertex of one to the left
vertex of the other, and a line drawn from the right vertex of
one to the right vertex of the other.  The new lines drawn
are inessential edges, and can be collapsed away to
yield $\Theta_4$.
\end{itemize}

So we are left with 4 basepointed graphs.  Two come
from the rose $R_4$, depending upon where we place the
basepoint, and the other two come from $\Theta_4$ in
a similar manner.  In particular, only $\Theta_4$
(with the basepoint $*$ placed in the middle of one of
its edges)
can contribute a $2$-simplex to our complex, and the
relevant marked graph only has one maximal subforest
(up to an isomorphism of the graph).  Hence it contributes
exactly two $2$-simplices, which join together
to form a square.
Consequently it is clear that $H^2(Q_{3}^\omega; \Z/3)=0$,
which is all we needed to show to prove that
$Q_{3}^\omega$ contributes nothing more to our
cohomology calculations.
  
For the final case of
considering the contributions of $H^3(\tilde Q_{2}; \Z/3)$,
we again use arguments like
those in \cite{[V]} and Proposition \ref{t16} from Part \ref{p2}.
We show that all of the $3$-simplices in 
$\tilde Q_{2}$ can be collapsed away,
so that $H^3(\tilde Q_{2}; \Z/3)$ is
necessarily zero.  The relevant graphs which can
give $3$-simplices are:

\bigskip

\input{fg6.pic}
\begin{figure}[here]
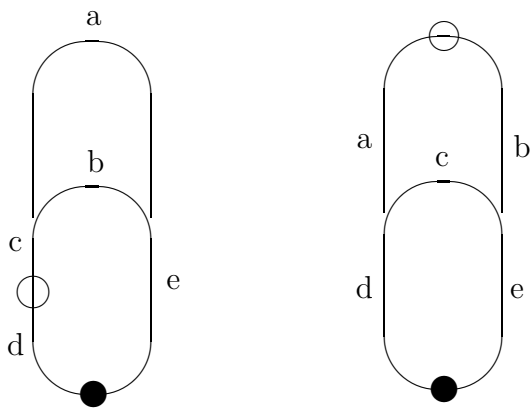

\caption{\label{fg6} Graphs giving some $3$-simplices}
\end{figure}

For the graphs in Figure \ref{fg6},
the basepoint is given by s filled-in dot. 
The open circle $\circ$
is the other ``distinguished point'' of the
graph 
and indicates where a $\Theta$-graph $\Theta_2$
should be attached.

The first of these graphs
has four subforests $\{b,c,e\}$, $\{b,c,d\}$,
$\{b,d,e\}$, and $\{c,d,e\}$,
each of which gives a
$3$-dimensional cube.
These can be collapsed in a manner
similar to that described
in Proposition \ref{t16} in Chapter \ref{c10}.
That is, the cube corresponding to the
first subforest $\{b,c,e\}$
has a free minusface obtained by collapsing $b$.
So we can collapse the interior of this cube away from that face.
Then the cube $\{b,c,d\}$ has a free
plusface corresponding to
$\{b,c\}$ and $\{b,d,e\}$ has a free
plusface corresponding to
$\{b,e\}$.  Both of these cubes can
be collapsed away from those
respective plusfaces.
This leaves the cube corresponding to
$\{c,d,e\}$ with all plusfaces free.
Thus we can disregard the first of the graphs
in Figure \ref{fg6}.
(For a more detailed
description of what plusfaces and minusfaces are,
along with several more examples, see part \ref{p2} of this
paper.)

The second of the graphs that give $3$-simplices
contributes 4 cubes, one for each of the
subforests $\{a,d,e\}$, $\{a,b,d\}$, $\{a,c,d\}$,
and  $\{a,c,e\}$.
The cube corresponding to $\{a,d,e\}$ has a
free minusface given by collapsing the
edge $a$.  So we can disregard this cube.
The remaining 3 cubes join together to
form a solid $3$-ball as described in
the proof of Proposition \ref{t16}, and so can also be
collapsed away. \END

\end{document}